\documentclass[12pt,a4paper,reqno]{amsart}
\usepackage{preemble}

\title[Invariant measures of the topological flow]
      {Invariant measures of the topological flow and measures at infinity on hyperbolic groups}
\date{\today}

\author{Stephen Cantrell}
\address{Department of Mathematics, 
University of Chicago,
Chicago, Illinois 60637, USA}
\email{scantrell@uchicago.edu}

\author{Ryokichi Tanaka}
\address{Department of Mathematics, 
Kyoto University,
Kyoto 606-8502, JAPAN}
\email{rtanaka@math.kyoto-u.ac.jp}

\date{\today. 
\\
2020 \textit{Mathematics Subject Classification}. Primary 20F67; Secondary 37D35, 37D40. 
\\
\textit{Key words and phrases}.
Hyperbolic groups, Patterson-Sullivan measures and symbolic coding}

\begin{document}

\maketitle

\begin{abstract}
We show that for every non-elementary hyperbolic group, an associated topological flow space admits a coding based on a transitive subshift of finite type.
Applications include regularity results for Manhattan curves, the uniqueness of measures of maximal Hausdorff dimension with potentials, and the real analyticity of intersection numbers for families of dominated representation, thus providing a direct proof of a result established by Bridgeman, Canary, Labourie and Sambarino in 2015.
\end{abstract}


\section{Introduction}
The purpose of this paper is to show that there exists a coding based on a transitive subshift of finite type for a flow associated with every non-elementary hyperbolic group.
Applications include the real analyticity of Manhattan curves (introduced by Burger in \cite{BurgerManhattan} and generalized in \cite{CT}),
the uniqueness of measure of maximal Hausdorff dimension with potentials,
and the real analyticity of intersection numbers for perturbations of dominated representations, providing a direct proof of a result shown by Bridgeman {\it et al.}\ in \cite{BCLS}.
We also apply our results to give a multifractal analysis of harmonic measures for finite range random walks and discuss, for explicit examples, the relationships between various measure classes on the boundary of the group.

Let $\G$ be a non-elementary hyperbolic group and $\partial \G$ be the ideal boundary equipped with a quasi-metric.
We consider the boundary square $\partial^2 \G:=(\partial \G)^2\setminus\{\text{diagonal}\}$ endowed with the diagonal action by $\G$.
Given a H\"older continuous cocycle $\k: \G \times \partial^2 \G \to \R$ coming from the Busemann function for a strongly hyperbolic metric,
we define a compact space $\Fc_\k:=\G \backslash (\partial^2 \G \times \R)$ equipped with an $\R$-action which we call a \textit{topological flow} (see Section \ref{Sec:TF} for the precise definitions).
The space $\Fc_\k$ is the one constructed by Mineyev \cite{MineyevFlow}.
Our main coding theorem reads the following:

\begin{theorem}\label{Thm:enhanced-intro}
For every non-elementary hyperbolic group $\G$, 
let $\Fc_\k$ be the topological flow space associated with a H\"older continuous cocycle $\k$ defined by a strongly hyperbolic metric. There exists a topologically transitive subshift of finite type $(\SS_0, \s)$ with a positive H\"older continuous function $r_0$
such that a continuous map from the associated suspension flow based on $(\SS_0, \s, r_0)$ to $\Fc_\k$
is surjective, equivariant with the flows and bounded-to-one.
\end{theorem}

The proof is given in Section \ref{Sec:enhanced} (see Theorem \ref{Thm:enhanced}).
A main contribution of this work is that, in the above theorem, the coding of our flow is based over a shift space $(\SS_0, \s)$ which is topologically transitive and we can find such a coding for every non-elementary hyperbolic group.
We call it an {\it enhanced coding}.
This enables us to apply thermodynamic formalism to the flow in a strength that has been missing so far, and to establish various regularity results concerning measures on the boundary.

A flow associated with hyperbolic groups was introduced by Gromov \cite{GromovHyperbolic}.
Since then, there have been a measure theoretical construction due to Furman (\cite{FurmanCoarse}, \cite{BaderFurman}) and a topological construction due to Mineyev \cite{MineyevFlow}.
The latter machinery is a key tool in the theory of Anosov representations; see e.g., \cite{BCLS}, a survey \cite{BCSsurvey} and references therein.
Although there are known codings on such flows e.g., via horofunctions (Busemann functions) \cite{CoornaertPapadopoulos}, or via Cannon's automatic structures \cite{BourdonThese} and \cite{TopFlows},
they are constructed over subshifts of finite type which may not be transitive (except for special cases such as Fuchsian groups associated with specific sets of generators).
For the case of Fuchsian groups, see the classical results by Bowen and Series \cite{BowenSeriesMarkov}. 
Recently, Constantine, Lafont and Thompson have established Bowen-Series type coding for flows on locally $\mathrm{CAT}(-1)$ metric spaces \cite{CLTstrong}.
Their construction is based on showing that the flows admit a metric Anosov structure in the sense of Pollicott \cite{PollicottFlow}.
Note that it is not known as to whether the Mineyev flow carries such a structure (see the discussion in \cite[Introduction]{BCSsurvey}).
Their result also applies to the flow associated with a projective Anosov representation \cite[Theorem C]{CLTstrong} and thus provides a self-contained proof on the existence of a coding in the setting of \cite{BCLS}.
We note, however, that there exists a non-elementary hyperbolic group for which any representation to a general linear group (over any field) factors through a finite group \cite[Theorem 8.1]{KapovichRep} as discussed in \cite[Remark 3.6]{BochiPotrieSambarino}.
Our result in Theorem \ref{Thm:enhanced-intro} applies to such hyperbolic groups without imposing any additional assumptions.

We construct a coding of the flow based on Cannon's automatic structure, which we know every hyperbolic group admits (for the proof, see e.g., \cite[Theorem 3.2.2]{Calegari}).
The underlying finite directed graph may have several distinct components
and consequently known codings (e.g., the one used in \cite{TopFlows}) 
are not necessarily based over a transitive subshifts of finite type.
Despite this issue it is still possible to use such a coding to prove some results.
Indeed, an approach due to Calegari and Fujiwara \cite{CalegariFujiwara2010} is to analyze each maximal (in the sense of spectrum) component in Cannon's automatic structure and to show that certain values of interest agree across components.
We show, by employing ergodic theory on the flow space, that
any one of the maximal components covers the flow space; see more precise discussions in Section \ref{Sec:enhanced}. (In a special setting where the group acts on a $\mathrm{CAT}(-1)$-space, an analogous argument has been carried out by Bourdon \cite[Lemme 5.1.2]{BourdonThese}.)
This part involves studying invariant measures under the flow whereas the conclusion is a purely topological one.
We note that there can exist more than one maximal component and the set of such components depend on the choice of the roof function (or, the cocycle).
Therefore, at first sight our coding is not fully constructive;
however, in fact, one may take \textit{any} component whose adjacency matrix has the maximal modulus $e^{v_S}$, where $v_S$ is the exponential growth rate for a chosen finite symmetric set of generators in the automatic structure (see Proposition \ref{Prop:component}).
That component again may not be unique and depend on the choice of a set of generators and the automaton, but one can find and use it for computation whenever any reasonably small group presentation is available, e.g., \cite[Section 5]{CT}.
We discuss some concrete examples where one has access to explicit coding structures in Section \ref{Sec:example}.

In the course of our discussion, we develop a general theory on invariant Radon measures on the  boundary square without any reference to the topological flow (see Section \ref{Sec:invariant}).
The approach was initiated by Kaimanovich \cite{KaimanovichInvariant} in the setting of negatively curved manifold.
In the same paper, it was suggested to carry out the argument entirely on the group itself (ibid, Remark 2 in Section 2.1).
Our framework unifies known approaches in the theory of (generalized) geodesic flows, random walks and dominated (or, projective Anosov) representations for hyperbolic groups.
Let us briefly describe the case when we have a dominated representation $\rho:\G \to \GL(m, \R)$ for an integer $m \ge 2$ (see  Section \ref{Sec:main-rep} for the precise definitions in the following discussion).
Let $\Dc_\G$ be the set of hyperbolic metrics which are left-invariant, quasi-isometric to some (equivalently, any) word metric in $\G$. (For more on the structure of $\Dc_\G$, see recent results in \cite{OregonReyes}.)
If we define 
\[
\psi_\rho(x, y):=\log \|\rho(x^{-1}y)\| \quad \text{for $x, y \in \G$},
\]
where $\|\cdot\|$ stands for a fixed Euclidean norm in $\R^m$.
then $\psi_\rho$ defines a $\G$-invariant tempered potential relative to a metric $d$ in $\Dc_\G$ (Definition \ref{Def:tempered}); informally speaking, it behaves as a rough geodesic metric in $\G$ in the same quasi-isometric class although it is not necessarily symmetric in $x$ and $y$.
Given any ergodic $\G$-invariant Radon measure $\L$ on $\partial^2 \G$,
the {\it local intersection number} of $\psi_\rho$ relative to $d$ in $\Dc_\G$
is defined by the following constant
\[
\t(\rho/d:\L):=\lim_{t \to \infty}\frac{\p_\rho(\g(0), \g(t))}{d(\g(0), \g(t))},
\]
for $\L$-almost every $(\x_-, \x_+) \in \partial^2 \G$,
where $\g$ is any quasi-geodesic converging to $\x_+$ and the limit exists independently of the choice of $\g$.

Considering the contragredient representation $\check \rho:\G \to \GL(m,\R)$, 
we define $\psi_{\check \rho}(x, y)$ for $\check \rho$ similarly to the case of $\rho$,
and a pair of cocycle $c_\rho(x, \x)$ and $c_{\check \rho}(x, \x)$ analogous to Busemann cocycles in the case of hyperbolic metrics.
The pair
admits a Gromov product $[\x | \y]_{\rho}$ for $(\x, \y) \in \partial^2 \G$ in terms of an associated pair of equivariant representations with $\rho$ on the boundary.
In the case of dominated representations,
$[\x|\y]_\rho$ defines a locally H\"older continuous function on $\partial^2 \G$.
If we consider the series
\[
\sum_{x \in \G}\exp(-s\log\|\rho(x)\|) \quad \text{for $s \in \R$},
\]
then the divergence exponent $v_\rho$ in $s$ is positive and finite and 
we call $v_\rho$ the {\it exponential growth rate} for $\rho$.
By definition, $v_{\check \rho}=v_\rho$, and
associated with the pair of cocycles,
finite measures $\m_\rho$ and $\m_{\check \rho}$ on $\partial \G$ are defined by the Patterson-Sullivan procedure.
Combining with the Gromov product for the pair of cocycles, 
we obtain a $\G$-invariant Radon measure on $\partial^2 \G$
given by
\[
\Lambda_\rho = \exp(v_\rho [\x|\y]_\rho)\m_{\check \rho}\otimes \m_{\rho}.
\]
The measure $\L_{\rho}$ is ergodic with respect to the $\G$-action on $\partial^2 \G$ by a generalized Hopf argument (see Theorem \ref{Thm:Hopf}).
Furthermore, $\L_\rho$ is an analogue of a Bowen-Margulis measure in the case of hyperbolic metrics and we have the following uniqueness result with respect to the maximality of Hausdorff dimension.
For the flip-involution $\i(\x, \y):=(\y, \x)$ on $\partial^2 \G$, let $\check \L:=\L\circ \i$.

\begin{theorem}\label{Thm:UMMD_ergodic-flip_intro}
Let $d$ be a hyperbolic metric in $\Dc_\G$, and $\rho:\G \to \GL(m, \R)$ be a dominated representation for $m \ge 2$. For any $\G$-invariant Radon measure $\L$ which is ergodic on $\partial^2 \G$,
if
\[
\dim_H(\L, q_\times) \ge v_\rho(\t(\check \rho/d:\check \L)+\t(\rho/d:\L)),
\]
then $\L$ is a constant multiple of $\L_\rho$, where $q_\times$ is defined in $\partial^2 \G$ by the maximum of quasi-metric $q(\x_1, \x_2)=e^{-(\x_1|\x_2)}$ associated with $d$.
\end{theorem}

We prove the above result in a more general form in Theorem \ref{Thm:UMMD_ergodic-flip}.
For example, if $\G$ is the fundamental group of closed Riemann surface of genus $2$ with the translation generators $S$, and further the representation $\rho$ induces isometries admitting a regular octagon as a Dirichlet domain on the hyperbolic plane,
then $v_\rho=2$ and a numerical experiment indicates that
\[
\t(\rho/d_S:\L_\rho)=1.13837...,
\]
where $d_S$ is the word metric relative to $S$ (see Section \ref{Sec:surface}).

\bigskip

The organization of this paper is the following:
In Section \ref{Sec:pre}, we provide preliminaries on hyperbolic groups, the Patterson-Sullivan construction with potentials and invariant Radon measures on the boundary square. 
In Section \ref{Sec:TF}, we construct the topological flow space due to Mineyev \cite{MineyevFlow} based on \cite{TopFlows}, whose results are generalized for our purposes.
In Section \ref{Sec:symbolic}, we introduce symbolic dynamics including thermodynamic formalism, and prove Theorem \ref{Thm:enhanced-intro} in Theorem \ref{Thm:enhanced}.
In Section \ref{Sec:application-coding}, we provide some applications of Theorem \ref{Thm:enhanced-intro} to the real analyticity of Manhattan curves for pairs of strongly hyperbolic metrics (Theorem \ref{Thm:Manhattan}) and word metrics (Theorem \ref{Thm:Manhattan_word}) improving upon results in \cite[Theorems 1.3 and 4.16]{CT} and the uniqueness of measure of maximal Hausdorff dimension with potentials (Theorems \ref{Thm:UMMD} and \ref{Thm:UMMD_ergodic-flip}) generalizing \cite[Theorem 1.1]{TopFlows}.
In Section \ref{Sec:RW}, we discuss finite range random walks on hyperbolic groups and their harmonic measures, showing that a multifractal profile function is real analytic (Theorem \ref{Thm:Manhattan_harm} and Remark \ref{Rem:multifractal}).
In Section \ref{Sec:main-rep}, we discuss dominated representations in the sense of Bochi, Potrie and Sambarino \cite{BochiPotrieSambarino}, prove Theorem \ref{Thm:UMMD_ergodic-flip_intro} as a special case of Theorem \ref{Thm:UMMD_ergodic-flip}, define the intersection number for a pair of dominated representations via a large deviation result (Theorem \ref{Thm:intersection_rep}) without appealing to Bowen's equidistribution theorem in the setting of metric Anosov flows, which may be of independent interest, and give a direct proof of the real analyticity under variations parameterized by real analytic manifolds.
In Section \ref{Sec:example}, we discuss a free product group and a surface group with the translation generators and compare various finite measures defined on the boundary such as Patterson-Sullivan measures associated with different hyperbolic metrics and a harmonic measure, indicating a possible numerical approach to understand their measure classes based on the local intersection numbers.

\subsection*{Notations}
We write $C, C', C'', \dots$ for constants whose exact values may vary from line to line,
and $C_\d, C_{\d, R}, \dots$ for constants depending on parameters $\d, R, \dots$.
For two real-valued functions $f(t)$ and $g(t)$ defined on some set,
we write $f(t)=g(t)+O(1)$ if there exists a constant $C$ such that $|f(t)-g(t)| \le C$ on the set, and replace $O(1)$ by $O_\d, O_R, \dots$ to indicate that the implied constant $C$ depends on some parameters $\d, R, \dots$.
We write $\Leb_{[0, T)}$ the Lebesgue measure $dt$ restricted on the interval $[0, T)$ where $T$ is a positive real number,
$\Z_{\ge 1}:=\{1, 2, \dots\}$ the set of positive integers, and $\#A$ the cardinality of a set $A$.

\section{Preliminaries}\label{Sec:pre}

\subsection{Hyperbolic groups}\label{Sec.hyp}

Let $(X, d)$ be a metric space.
The Gromov product is defined by
\[
(x|y)_w:=\frac{1}{2}(d(w, x)+d(w, y)-d(x, y)) \quad \text{for $x, y, w \in X$}.
\]
For a non-negative real number $\d \ge 0$, we say that a metric space $(X, d)$ is $\d$-{\it hyperbolic} if
\begin{equation}\label{Eq:hyperbolic}
(x|y)_w \ge \min\left\{(x|z)_w, (y|z)_w\right\}-\d \quad \text{for all $x, y, z, w \in X$}.
\end{equation}
A metric space is called {\it hyperbolic} if it is $\d$-hyperbolic for some $\d \ge 0$.

Let $\G$ be a finitely generated group, i.e., there exists a finite set of generators $S$ such that each group element can be represented as a word with letters in $S$.
We assume that $S$ is {\it symmetric}, i.e., $s^{-1} \in S$ whenever $s \in S$, enlarging a set of generators if necessary.
For any finite symmetric set of generators $S$, we define the word norm
\[
|x|_S:=\min\{ k \ge 0 \ : \ x=s_1 \dots s_k \ \text{for $s_1, \dots, s_k \in S$}\} \quad \text{for $x \in \G$},
\]
where $|o|_S=0$ (for the identity element $o:=\id$), and a word metric $d_S(x, y):=|x^{-1}y|_S$ for $x, y \in \G$.
Note that $d_S$ is invariant under left-translation.
We say that $\G$ is a hyperbolic group if $(\G, d_S)$ is a hyperbolic metric space for some (equivalently, any) finite symmetric set of generators $S$.
A hyperbolic group is called \textit{non-elementary} if it is non-amenable, and \textit{elementary} otherwise.
Elementary hyperbolic groups are either finite groups or contain an infinite cyclic group as a finite index subgroup.
Basic examples of non-elementary hyperbolic groups include finite rank free groups of rank at least $2$,
cocompact lattices of $\PSL(2, \R)$,
and
the fundamental groups of closed Riemannian manifolds with negative sectional curvature.

We say that two metrics $d, d_\ast$ in $\G$ are {\it quasi-isometric}
if there exist constants $L>0$ and $C \ge 0$ such that
\[
L^{-1} d(x, y)-C \le d_\ast(x, y) \le L d(x, y)+C \quad \text{for all $x, y \in \G$}.
\]
Let $\Dc_\G$ be the set of metrics which are left-invariant, hyperbolic and quasi-isometric to some (equivalently, any) word metrics in $\G$.

For example, if $\G$ is a cocompact lattice of $\PSL(2, \R)$ and torsion-free,
then a natural isometric action of $\G$ on the hyperbolic plane $\H^2$ yields a metric in $\G$ by restricting the Poincar\'e metric $d_{\H^2}$ on an orbit $\G w$ for some fixed base point $w$ in $\H^2$.
If $\G$ is not torsion-free and the stabilizer at $w$ is non-trivial,
then defining $d_w(x, y):=d_{\H^2}(xw, yw)+1$ if $x \neq y$ and $d_w(x, y):=0$ if $x=y$,
we obtain a metric $d_w$ in $\Dc_\G$.

For a metric $d$ in $\G$,
we say that a map $\g:I \to \G$ for an interval $I$ in $\R$ (where $I$ is possibly bounded, semi-infinite or bi-infinite) is
an \textit{$(L, C)$-quasi-geodesic} for constants $L>0$ and $C \ge 0$ if 
\[
L^{-1}|s-t|-C \le d(\g(s), \g(t)) \le L|s-t|+C \quad \text{for all $s, t \in I$},
\]
and a \textit{$C$-rough geodesic} for a constant $C \ge 0$ if
\[
|s-t|-C \le d(\g(s), \g(t))\le |s-t|+C \quad \text{for all $s, t \in I$}.
\]
We say that $\g$ is a \textit{geodesic} if it is a $0$-rough geodesic.
Furthermore a metric space $(\G, d)$ is called $C$-{\it roughly geodesic}
if for all $x, y \in \G$ there exists a $C$-rough geodesic $\g:[0, T] \to \G$ such that $\g(0)=x$ and $\g(T)=y$,
and {\it roughly geodesic} if it is $C$-roughly geodesic for some $C \ge 0$.
A geodesic metric space is a $0$-roughly geodesic metric space.
Every metric in $\Dc_\G$ is roughly geodesic;
this follows from a ``rough reparameterization'' of quasi-geodesics in a hyperbolic metric space-----but note that it is not necessarily geodesic.
Moreover if $d \in\Dc_\G$ is $C_0$-roughly geodesic,
then the {\it Morse lemma} holds:
for every $(L, C)$-quasi-geodesic $\g$ in $(\G, d)$
there exists a $C_0$-rough geodesic $\g_0$ such that $\g$ and $\g_0$ are within Hausdorff distance $C(L, C, C_0, \d)$,
where $\d$ is a hyperbolic constant of the metric $d$
\cite[Proposition 5.6]{BonkSchramm}.

\begin{definition}\label{Def:SH}
A hyperbolic metric $d$ in $\G$ is called {\it strongly hyperbolic} if it satisfies the following:
There exist positive constants $c, R_0>0$ such that for all $R \ge R_0$, and all $x, x', y, y' \in \G$,
if
\[
d(x, y)-d(x, x')+d(x', y')-d(y, y') \ge R,
\]
then
\begin{equation}\label{Eq:SH}
|d(x, y)-d(x', y)-d(x, y')+d(x', y')| \le \exp(-c R).
\end{equation}
\end{definition}

Mineyev \cite{MineyevFlow} has shown that a strongly hyperbolic metric exists in $\Dc_\G$ for any hyperbolic group, see also \cite{NicaSpakula} for another construction based on the Green metric.
Note that strongly hyperbolic metrics which are in $\Dc_\G$ are roughly geodesic.
The strongly hyperbolic condition means that
if $x, x'$ are fixed and $y, y'$ are far away,
then \eqref{Eq:SH} gives a fine control on the difference between two Busemann functions (see below for the definition)
depending on two asymptotic directions in which $y$ and $y'$ follow respectively.
We use the notation $\hat d$ for a strongly hyperbolic metric following Mineyev's hat metric.

Let us define the ideal boundary $\partial \G$ of $\G$.
Fix a hyperbolic metric in $\Dc_\G$ and define the Gromov product associated with it.
We say that a sequence $\{x_n\}_{n=0}^\infty$ {\it diverges}
if $(x_n|x_m)_w$ diverges as $\min\{n, m\}$ tends to infinity for some (equivalently, any) $w \in \G$.
Let us fix the identity $o=\id$ as a base point.
Two divergent sequences $\{x_n\}_{n=0}^\infty$ and $\{y_n\}_{n=0}^\infty$ are called {\it equivalent} if 
$(x_n|y_m)_o$ diverges as $\min\{n, m\}$ tends to infinity.
For a divergent sequence,
let us denote by $\x:=[\{x_n\}_{n=0}^\infty]$ the equivalence class, and write $x_n \to \x$ as $n \to \infty$.
We define $\partial \G$ as the set of equivalence classes of divergent sequences in $\G$
and call $\partial \G$ the {\it ideal boundary} of $\G$.
If $d$ is in $\Dc_\G$ and $C$-roughly geodesic,
then for each $\x$ in $\partial \G$ there exists a $C$-rough geodesic $\g:[0, \infty) \to \G$ 
such that $\g(0)=o$ and $\g(n) \to \x$ as $n \to \infty$.
Furthermore for each pair $\x, \y$ in $\partial \G$ with $\x \neq \y$
there exists a $C$-rough geodesic $\g:(-\infty, \infty) \to \G$
such that $\g(-n) \to \x$ and $\g(n) \to \y$ respectively as $n \to \infty$ \cite[Proposition 5.2]{BonkSchramm}.

We extend the Gromov product to points in $\G \cup \partial \G$.
Let
\[
(\x|\y)_o:=\sup\left\{\liminf_{n \to \infty}(x_n|y_n)_o \ : \ \x=[\{x_n\}_{n=0}^\infty], \  \y=[\{y_n\}_{n=0}^\infty]\right\},
\]
for $\x, \y \in \G \cup \partial \G$,
where if $\x$ or $\y$ is in $\G$, then we understand that $\{x_n\}_{n=0}^\infty$ is a constant sequence, i.e., $x_n=\x$ for all $n\ge 0$.
If $\{x_n\}_{n=0}^\infty$ and $\{x_n'\}_{n=0}^\infty$ (resp.\ $\{y_n\}_{n=0}^\infty$ and $\{y_n'\}_{n=0}^\infty$) are divergent and equivalent,
then by \eqref{Eq:hyperbolic},
\[
\liminf_{n \to \infty}(x_n'|y_n')_o \ge \limsup_{n \to \infty}(x_n|y_n)_o-2\d,
\]
where $\d$ is a hyperbolic constant for the fixed metric.
Therefore it holds that
\[
(\x|\y)_o \ge \min\left\{(\x|\z)_o, (\z|\y)_o\right\}-3\d \quad \text{for $\x, \y, \z \in \G \cup \partial \G$}.
\]
The \textit{quasi-metric} in $\partial \G$ associated to a metric in $\Dc_\G$ is defined by
\[
q(\x, \y):=\exp\(-(\x|\y)_o\) \quad \text{for $\x, \y \in \partial \G$},
\]
where $q(\x, \x)=0$.
A quasi-metric satisfies that $q(\x, \y)=0$ if and only if $\x=\y$, $q(\x, \y)=q(\y, \x)$,
and 
$q(\x, \y) \le e^{2\d}\max\left\{q(\x, \z), q(\z, \y)\right\}$
for $\x, \y, \z \in \partial \G$.
We endow $\partial \G$ with a topology by declaring that sets $\{\y \in \partial \G \ : \ q(\x, \y)<\e\}$ for $\x \in \partial \G$ and $\e>0$ form an open basis.
It is known that $q^{\e_0}$ for some $\e_0>0$ is bi-Lipschitz to a genuine metric.
The Morse lemma implies that the ideal boundaries defined in terms of metrics from $\Dc_\G$
are mutually homeomorphic.

Letting $d$ be a hyperbolic metric in $\Dc_\G$ and $w$ be in $\G$,
we define the Busemann function 
$\b_w: \G \times (\G \cup \partial \G) \to \R$ based at $w$
for $x \in \G$,
\[
\b_w(x, \x):=\sup\left\{\limsup_{n \to \infty}(d(x, x_n)-d(w, x_n)) \ : \ \x=[\{x_n\}_{n=0}^\infty]\right\} \quad
\text{for $\x \in \partial \G$},
\]
and 
$\b_w(x, y)=d(x, y)-d(w, y)=d(w, x)-2(x|y)_w$
for $y \in \G$.
Note that we have
\[
\b_w(x, \x)=d(w, x)-2(x|\x)_w +O_\d,
\]
and further,
\[
\b_w(xy, \x)=\b_w(y, x^{-1}\x)+\b_w(x w, \x)+O_\d,
\]
for $w, x, y \in \G$ and $\x \in \partial \G$.

If we fix a strongly hyperbolic metric $\hat d$ in $\Dc_\G$,
then
we denote the corresponding (extended) Gromov product and Busemann functions by
$\langle \x|\y\rangle_w$ and $\hat \b_w(x, \x)$ respectively. 
The definition of a strongly hyperbolic metric (Definition \ref{Def:SH}) implies that
for each $w, x$ in $\G$,
the Busemann function $\hat \b_w(x, \x)$ based at $w$ is obtained as a genuine limit 
\[
\hat \b_w(x, \x)=\lim_{n \to \infty}\(\hat d(x, x_n)-\hat d(w, x_n)\)=\hat d(w, x)-2\langle x | \x \rangle_w \quad \text{for $(x, \x) \in \G \times \partial \G$},
\]
where $\x=[\{x_n\}_{n=0}^\infty]$,
and it is continuous with respect to $\x$ in $\partial \G$.
Moreover, we have an identity
\[
\hat \b_w(xy, \x)=\hat \b_w(y, x^{-1}\x)+\hat \b_w(x w, \x) \quad \text{for $w, x, y \in \G$ and $\x \in \G\cup \partial \G$}.
\]
This exact identity (without ``quasification") is used in later sections.

\subsection{Patterson-Sullivan construction with tempered potentials}\label{Sec:PStempered}

Let us consider the quasi-metric $q(\x, \y)=\exp(-(\x|\y)_o)$ in $\partial \G$ associated with a metric $d$ in $\Dc_\G$
and define
\[
B(\x, r):=\{\y \in \partial \G \ : \ q(\x, \y) < r\} \quad \text{for $\x \in \partial \G$ and real $r \ge 0$}.
\]
We define the {\it shadow}
\[
\Oc(x, R):=\{\x \in \partial \G \ : \ (\x|x)_o \ge d(o, x)-R\} \quad \text{for $x \in \G$ and real $R\ge 0$}.
\]
If the metric is $\d$-hyperbolic for $\d \ge 0$,
then for each $r \ge 0$ and each $R \ge r+3\d$,
the following holds: 
For all $\x \in \partial \G$ and all $x \in \G$ such that $(o|\x)_x \le r$, we have that
\begin{equation}\label{Eq:shadow-ball}
B\(\x, e^{-3\d+R-d(o, x)}\) \subset \Oc\(x, R\) \subset B\(\x, e^{3\d+R-d(o, x)}\),
\end{equation}
(cf.\ \cite[Proposition 2.1]{BHM11}).
Let $\Oc(x, R)$ (resp.\ $\Oc_\ast(x, R)$) be shadows associated with $d$ (resp.\ $d_\ast$), where $d$ and $d_\ast$ are in $\Dc_\G$.
Since these two metrics are quasi-isometric,
for each $R \ge 0$,
\[
\Oc(x, R) \subset \Oc_\ast(x, R') \quad \text{for every $x \in \G$},
\]
where $R'$ is a positive constant depending only on $R$, the quasi-isometric constants and the hyperbolic constants.
This follows from the Morse lemma.
Hence for each fixed $x \in \G$, up to changing $R$, the shadows $\Oc(x, R)$ are including each other for any metrics in $\Dc_\G$.

For each hyperbolic metric $d$ in $\Dc_\G$, let
\[
B(x, r):=\{y \in \G \ : \ d(x, y) \le r\} \quad \text{for $x \in \G$ and real $r \ge 0$}.
\]
The \textit{exponential growth rate} relative to the metric is defined by
\[
v:=\limsup_{r \to \infty}\frac{1}{r}\log \#B(o, r),
\]
where $\#A$ stands for the cardinality of a set $A$.
If $\G$ is a non-elementary hyperbolic group, then $v$ is finite and strictly positive.

Recall that for each $d$ in $\Dc_\G$,
we have a finite Borel regular measure $\m$ called a \textit{Patterson-Sullivan measure} on $\partial \G$
satisfying that
\[
\exp(-4\d v) \le \exp(v \b_o(x, \x))\cdot \frac{d x_\ast\m}{d\m}(\x) \le \exp(4\d v),
\]
for $x \in \G$ and for $\m$-almost every $\x \in \partial \G$,
where $\d$ is a hyperbolic constant of the metric
(cf.\ \cite[Th\'eor\`eme 5.4]{Coornaert1993}).
We apply the construction to more general functions.

For a function $\p: \G \times \G \to \R$,
let
\[
(x|y)^\p_w:=\frac{1}{2}(\p(x, w)+\p(w, y)-\p(x, y)) \quad \text{for $x, y, w \in \G$},
\]
where we {\em do not} assume that $\p(x, y)=\p(y, x)$.
We say that $\p$ satisfies (QE) (``quasi-extension'')
if there exist a function $(\cdot\ |\ \cdot)_w^\p: \G \times (\G \cup \partial \G) \to \R$ and a constant $C \ge 0$
such that
\begin{equation}\label{Eq:QE}\tag{QE}
\limsup_{n \to \infty}(x|x_n)^\p_w -C \le (x|\x)^\p_w \le \liminf_{n \to \infty}(x|x_n)^\p_w+C
\end{equation}
for every $(x, \x) \in \G \times (\G \cup \partial \G)$ and for any $\{x_n\}_{n=0}^\infty$ in the class $\x$.
Let 
\[
\b^\p_w(x, \x):=\sup\left\{\limsup_{n \to \infty}(\p(x, x_n)-\p(w, x_n)) \ : \ \x=[\{x_n\}_{n=0}^\infty]\right\},
\]
for $(x, \x) \in \G \times \partial \G$.
If $\p$ satisfies (QE), then we have
\[
\b^\p_w(x, \x)=\p(x, w)-2(x|\x)^\p_w+O(1) \quad \text{for all $(x, \x) \in \G \times \partial \G$}.
\]
Moreover, if $\p$ is invariant under left-translations, i.e., $\p(g x, g y)=\p(x, y)$ for $g, x, y \in \G$,
then 
\[
\b^\p_w(xy, \x)=\b^\p_{w}(y, x^{-1}\x)+\b^\p_w(x w, \x)+O(1) \quad \text{for $x, y \in \G$ and $\x \in \partial \G$}.
\]
We say that $\p$ satisfies (RG) (``roughly geodesic'') relative to $d$ in $\Dc_\G$
if for all large enough constants $C, R \ge 0$ there exists a constant $C' \ge 0$
such that for all $C$-rough geodesics $\g$ connecting $w$ and $y$,
and for all $x$ in the $R$-neighborhood of $\g$,
\begin{equation}\label{Eq:RG}\tag{RG}
|\p(w, x)+\p(x, y)-\p(w, y)| \le C'.
\end{equation}
The Morse lemma shows that if $\p$ satisfies \eqref{Eq:RG} relative to $d$,
then it satisfies \eqref{Eq:RG} relative to $d'$ for $d, d' \in \Dc_\G$,
furthermore, if $\p$ satisfies \eqref{Eq:RG} for some $d \in \Dc_\G$, 
then there exist $C, R \ge 0$ such that
\begin{equation*}
|\b^\p_o(x, \x)+\p(o, x)| \le C \quad \text{for all $\x \in \Oc(x, R)$ and all $x \in \G$}.
\end{equation*}

\begin{definition}\label{Def:tempered}
A function $\p: \G \times \G \to \R$ is called a {\it tempered potential} (relative to $d$ in $\Dc_\G$)
if $\p$ satisfies \eqref{Eq:QE} and \eqref{Eq:RG} for some (equivalently, any) metric $d$ in $\Dc_\G$.
Moreover, $\p$ is called {\it $\G$-invariant} if $\p$ is invariant under left-translations.
\end{definition}

\begin{example}\label{Ex:Manhattan}
If $d$ and $d_\ast$ are hyperbolic metrics in $\Dc_\G$,
then $d_\ast$ defines a $\G$-invariant tempered potential relative to $d$, which follows from the Morse lemma.
Moreover, for each $s \in \R$, the function $\p_s=s d_\ast$ defines a $\G$-invariant tempered potential relative to $d$.
\end{example}

Let $c(x, \x)$ be a real-valued (additive) quasi-cocycle on $\G \times \partial \G$, i.e.,
\[
c(xy, \x)=c(y, x^{-1}\x)+c(x, \x)+O(1) \quad \text{for all $x, y \in \G$ and $\x \in \partial \G$}.
\]
We say that a finite Borel measure $\m$ on $\partial \G$ satisfies (QC) (``quasi-conformal'') with a quasi-cocycle $c(x, \x)$ if
there exists a positive constant $C$ depending only on the quasi-cocycle such that
\begin{equation}\label{Eq:QC}\tag{QC}
C^{-1} \le \exp\(c(x, \x)\)\cdot\frac{dx_\ast \m}{d\m}(\x)\le C,
\end{equation}
for $x \in \G$ and for $\m$-almost all $\x \in \partial \G$.
For example,
if $d \in \Dc_\G$, then any Patterson-Sullivan measure (based at $o$) associated with $d$ satisfies (QC) with the corresponding Busemann quasi-cocycle $\b_o(x, \x)$.
More generally, we define the following.

\begin{definition}\label{Def:tempered-potential}
For a real value $\th$,
we say that $\p$ is a \textit{$\G$-invariant tempered potential with exponent $\th$ relative to $d$} if 
$\th$ is the abscissa of a convergence of the following series in $s$,
\[
\sum_{x \in \G}\exp\(-\p(o, x)-s d(o, x)\).
\]
\end{definition}

\begin{proposition}[Proposition 2.7 in \cite{CT}]\label{Prop:PSpotential}
Let $d$ be in $\Dc_\G$ and $\b_o(x, \x)$ be the corresponding Busemann quasi-cocycle.
If $\p$ is a $\G$-invariant tempered potential relative to $d$ with exponent $\th$,
then there exists a finite Borel measure $\m$ on $\partial \G$ satisfying \eqref{Eq:QC} with the quasi-cocycle $\b^\p_o(x, \x)+\th \b_o(x, \x)$.
Furthermore, for any such finite Borel measure $\m$,
there exist positive constants $C, R>0$ such that
\[
C^{-1} \le \exp(\p(o, x)+\th d(o, x)) \cdot \m\(\Oc(x, R)\) \le C \quad \text{for all $x \in \G$}.
\]
\end{proposition}

\begin{remark}\label{Rem:doubling}
Let us say that a Borel measure $\m$ on $\partial \G$ is {\it doubling} relative to a quasi-metric $q$
if $\m$ has a positive measure on each non-empty open set and there exists a positive constant $C>0$
such that
\[
\m\(B(\x, 2r)\)\le C\m\(B(\x, r)\) \quad \text{for all $\x \in \partial \G$ and all $r> 0$},
\]
where $B(\x, r)$ are (open) balls defined in terms of $q$.
It has been shown that 
for every $\G$-invariant tempered potential relative to a metric $d$ which is in $\Dc_\G$,
any finite Borel measure satisfying \eqref{Eq:QC} is doubling relative to $q$,
where $q$ is associated to $d$.
Furthermore the measure is ergodic with respect to the $\G$-action on $\partial \G$,
i.e., all $\G$-invariant Borel set in $\partial \G$ has either $0$ or the full measure \cite[Lemma 2.9]{CT}.
\end{remark}

\begin{remark}\label{Rem:exponent0}
If a $\G$-invariant tempered potential has the exponent $0$ relative to some $d$ in $\Dc_\G$,
then it has the exponent $0$ relative to every $d$ in $\Dc_\G$.
Indeed, the exponent relative to $d$ is given by
\[
\th=\limsup_{n \to \infty}\frac{1}{n}\log \sum_{x \in S(n, R_0)}e^{-\p(o, x)},
\quad \text{where} \quad S(n, R_0):=\{x \in \G \ : \ |d(o, x)-n|\le R_0\}
\]
for a large enough $R_0>0$,
and if $\th=0$, then Proposition \ref{Prop:PSpotential} implies that there exists a finite Borel measure $\m$ on $\partial \G$ such that
\[
C^{-1} \le \exp(\p(x, y))\cdot \m(\Oc(x, R))\le C \quad \text{for all $x \in \G$}.
\]
Further the shadows $\Oc(x, R)$ and $\Oc_\ast(x, R)$ associated with $d$ and $d_\ast$ in $\Dc_\G$, respectively, are including each other up to changing $R$ uniformly for $x \in \G$, 
and shadows $\Oc(x, R)$ for $x \in S(n, R_0)$ cover the boundary with a bounded multiplicity for all large enough integers $n$.
This shows that the exponent $0$ is relative to any other metric in $\Dc_\G$.
\end{remark}

\subsection{Invariant measures on the boundary square}\label{Sec:invariant}

Let us define
the boundary square $\partial^2 \G:=(\partial \G)^2\setminus \{\text{\rm diagonal}\}$ endowed with the restriction of the product topology.
We define the diagonal action of $\G$ on $\partial^2 \G$, where the action is continuous.

Let $c(x, \x)$ be a real-valued quasi-cocycle defined on $\G \times \partial \G$.
For any pair of such quasi-cocycles $c(x, \x)$ and $\check c(x, \x)$,
we say that the pair admits a \textit{Gromov product} on $\partial^2 \G$
if there exists a family of locally bounded (i.e., uniformly bounded from above and from below on each compact set) Borel measurable functions $[\,\cdot\,|\,\cdot\,]_w:\partial^2 \G \to \R$ for $w \in \G$ such that
\[
[x \x|x \y]_{xw}=[\x|\y]_w \quad \text{for $x \in \G$ and for $(\x, \y) \in \partial^2 \G$},
\]
and
\[
\check c(x, \x)+c(x, \y)=[\x|\y]_x-[\x|\y]_o +O(1),
\]
for all $x \in \G$ and for all $(\x, \y) \in \partial^2 \G$.
For example, if $d \in \Dc_\G$, 
then the pair of two identical Busemann quasi-cocycles associated with $d$ admits a Gromov product as twice of the one in the original sense. 
This generalized notion of Gromov product has already appeared in \cite{Sambarino}.

\begin{proposition}\label{Prop:Radon}
Let $c(x, \x)$ and $\check c (x, \x)$ be real-valued quasi-cocycles on $\G \times \partial \G$.
If there exist finite Borel measures $\m$ and $\check \m$ on $\partial \G$ satisfying \eqref{Eq:QC} with $c(x, \x)$ and $\check c(x, \x)$ respectively,
and the pair admits a Gromov product $[\,\cdot\,|\,\cdot\,]_w$ on $\partial^2 \G$ for $w \in \G$,
then there exists a $\G$-invariant Radon measure $\L$ on $\partial^2 \G$ equivalent to
\[
\exp([\x|\y]_o)\check \m \otimes \m.
\]
Moreover, if $\m$ and $\check \m$ have a positive measure on each open set in $\partial \G$,
then $\L$ has a positive measure on each open set in $\partial^2 \G$.
In the case when $c(x, \x)$ and $\check c(x, \x)$ are cocycles with the equality
\[
\check c(x, \x)+c(x, \y)=[\x|\y]_x-[\x|\y]_o \quad \text{for all $x \in \G$, $(\x, \y) \in \partial^2 \G$},
\]
then we may take 
\[
\L=\exp([\x|\y]_o)\check \m \otimes \m.
\]
\end{proposition}

\proof
This is a straightforward generalization of known cases.
If we define
\[
\n:=\exp([\x|\y]_o)\check \m \otimes \m,
\]
then for all $x \in \G$, by using $[x^{-1}\x|x^{-1}\y]_o=[\x|\y]_x$, we have that
\[
\frac{d x_\ast \n}{d\n}(\x, \y)=\exp([\x|\y]_x-[\x|\y]_o)\frac{d x_\ast\check \m}{d\check \m}(\x)\frac{d x_\ast \m}{d \m}(\y),
\]
for $\n$-almost every $(\x, \y) \in \partial^2 \G$.
Since $\m$ and $\check \m$ satisfy \eqref{Eq:QC} with $c(x, \x)$ and $\check c(x, \x)$ respectively, and the pair admits the Gromov product $[\,\cdot\,|\,\cdot\,]_w$ for $w \in \G$,
the function 
\[
\f(\x, \y):=\sup_{x \in \G}\frac{d x_\ast \n}{d\n}(\x, \y)
\]
is essentially uniformly bounded from above and uniformly bounded away from $0$.
If we define $\L:=\f(\x, \y)\n$, then $\L$ is $\G$-invariant.
Moreover $\L$ is Radon since the Gromov product is assumed to be locally bounded. 
For the second statement we note that, with the additional assumptions, $c(x,\x)$ and $\check{c}(x,\x)$ are genuine cocycles and the condition (\ref{Eq:QC}) holds with $C =1$ for both $\mu$ and $\check{\mu}$.
It follows that we can take $\f(\x, \y)\equiv 1$ for $\check \m \otimes \m$-almost every $(\x, \y)$ (see \cite[Remark 2.12]{TopFlows}) and the result follows.
\qed

\medskip

For a $\G$-invariant tempered potential $\p$ relative to some $d$ in $\Dc_\G$, 
let
\[
\check \p(x,y):=\p(y,x) \quad \text{for $x, y \in \G$}.
\]
Note that $\check \p$ is $\G$-invariant.
If $\check \p$ is also a tempered potential relative to some metric in $\Dc_\G$, then the above results are applied to
\[
\check c(x, \x):=\b_o^{\check \p}(x, \x) \quad \text{and} \quad c(x, \x):=\b_o^\p(x, \x) \quad \text{for $(x, \x)\in \G \times \partial \G$},
\]
for which the pair of the quasi-cocycles admits a Gromov product.

\begin{example}\label{Ex:interpolation}
For any pair of hyperbolic metrics $d$ and $d_\ast$ from $\Dc_\G$,
we recall that $\p_s=s d_\ast$ is a $\G$-invariant tempered potential relative to $d$ for each $s \in \R$ (Example \ref{Ex:Manhattan}).
Therefore it applies to Proposition \ref{Prop:PSpotential}, which yields the following:
For each $s \in \R$ there exists a constant $t:=\th(s)$ and a finite Borel regular measure $\m_{s, t}$ on $\partial \G$ satisfying \eqref{Eq:QC} with the quasi-cocycle $s\b_{\ast o}+t\b_o$,
where $\b_{\ast o}$ and $\b_o$ are Busemann quasi-cocycles for $d_\ast$ and for $d$ respectively.
Note that the quasi-cocycle $s \b_{\ast o}+t\b_o$ admits a Gromov product
defined by $2s (\x|\y)_{\ast o}+2t (\x|\y)_o$ for $\x, \y \in \partial \G$,
where $(\,\cdot\,|\,\cdot\,)_{\ast o}$ and $(\,\cdot\,|\,\cdot\,)_o$ are the Gromov products for $d_\ast$ and $d$ respectively.
Hence
there exists a $\G$-invariant Radon measure $\L_{s, t}$ on $\partial^2 \G$ equivalent to
\[
\exp\(2s (\x|\y)_{\ast o}+2 t(\x|\y)_o\)\m_{s, t}\otimes \m_{s, t},
\]
by Proposition \ref{Prop:Radon}.
\end{example}

\section{Topological flows}\label{Sec:TF}

\subsection{Topological flows: construction}\label{Sec:TopFlows}

Recall that $\partial^2 \G:=(\partial \G)^2\setminus \{\text{\rm diagonal}\}$ endowed with the diagonal action of $\G$.
If we fix a strongly hyperbolic metric $\hat d$ in $\Dc_\G$,
then for a fixed constant $C \ge 0$, for every $(\x, \y) \in \partial^2 \G$ there exists a $C$-rough geodesic $\g_{\x, \y}:\R \to \G$ 
such that $\g_{\x, \y}(-t) \to \x$ and $\g_{\x, \y}(t) \to \y$ respectively as $t \to \infty$.
Let us parameterize $\g_{\x, \y}$ by shifting $t \mapsto t+T$ if necessary so that
\[
\hat d(o, \g_{\x, \y}(0))=\inf\left\{\hat d(o, \g_{\x, \y}(t)) \ : \ t \in \R\right\}.
\]
We define the evaluation map by
\[
\ev: \partial^2 \G \times \R \to \G, \quad (\x, \y, t) \mapsto \g_{\x, \y}(t),
\]
where the space $\partial^2 \G \times \R$ equipped with the product topology.
Note that although $\ev$ depends on the choice of $C$-rough geodesics,
any other choice yields a map whose image lies within a distance at most some constant depending only on the metric.
Furthermore we may define $\ev$ as a Borel measurable map in the following way:
Let us define the set of $C$-rough geodesics $\g:\R \to \G$ and endow it with a pointwise convergence topology.
Fixing a finite set of generators equipped with a total order in $\G$,
we assign each $(\x, \y) \in \partial^2 \G$ to a $C$-rough geodesic $\g_{\x, \y}^{\min}$ 
which is lexicographically minimal as a sequence of group elements $\{\g(n)\}_{n \in \Z}$,
and define $\g_{\x, \y}(t):=\g_{\x, \y}^{\min}(\lfloor t\rfloor)$ for $t \in \R$
where $\lfloor t\rfloor$ is the largest integer at most $t$.
Defining the evaluation map by this assignment yields a Borel measurable map $\ev$. 

Let us define the cocycle
\begin{equation}\label{Eq:top-cocycle}
\k(x, \x, \y):=\frac{1}{2}\(\hat \b_o(x^{-1}, \x)-\hat \b_o(x^{-1}, \y)\) \quad \text{for $(x, \x, \y) \in \G \times \partial^2 \G$}.
\end{equation}
The group $\G$ acts on $\partial^2 \G \times \R$ continuously via the \textit{$(\G, \k)$-action},
\[
x\cdot\(\x, \y, t\):=(x\x, x\y, t-\k(x, \x, \y)) \quad \text{for $x \in \G$ and $(\x, \y, t) \in \partial^2 \G \times \R$}.
\]
It is known that the $(\G, \k)$-action on $\partial^2 \G \times \R$ is properly discontinuous and co-compact,
i.e., the quotient topological space is compact \cite[Lemma 3.2]{TopFlows}.
Let us denote the compact space by
\[
\Fc_\k:=\G \backslash (\partial^2 \G \times \R).
\]
Note that $\Fc_\k$ is Hausdorff and in fact, metrizable (cf.\ \cite[Theorem 60(d)]{MineyevFlow} and \cite[Remark 3.3]{TopFlows}).
We define a continuous $\R$-action, namely, a flow on $\Fc_\k$ induced from the translation on the $\R$-component in $\partial^2 \G \times \R$.
More precisely,
let 
\[
\widetilde \F_t(\x, \y, s):=(\x, \y, t+s) \quad \text{for $t \in \R$ and $(\x, \y, s) \in \partial^2 \G \times \R$}.
\]
Since $\{\widetilde \F_t\}_{t \in \R}$ and the $(\G, \k)$-action commute,
it induces an $\R$-action on $\Fc_\k$,
\[
\F_t: \Fc_\k \to \Fc_\k, \quad [\x, \y, s]\mapsto [\x, \y, t+s], \quad \text{for all $t \in \R$}.
\]
Let us call $\{\F_t\}_{t \in \R}$ the {\it topological flow} on $\Fc_\k$.

\subsection{Measures and ergodicity}\label{Sec:measures-ergodicity}

Let us consider a Radon measure $\L$ on $\partial^2 \G$,
i.e., $\L$ is Borel regular and finite on every compact set.
Further 
we consider the product measure
$\L \otimes dt$ on $\partial^2 \G \times \R$
where $dt$ is the Lebesgue measure on $\R$ normalized so that the interval $[0, 1]$ has a unit mass.
For any $\G$-invariant Radon measure $\L$ on $\partial^2 \G$,
the measure $\L\otimes dt$ is also Radon and invariant under the $(\G, \k)$-action on $\partial^2 \G \times \R$.
Moreover, $\L \otimes dt$ is invariant under the flow $\{\widetilde \F_t\}_{t \in \R}$,
and we naturally associate with it a flow-invariant finite measure $m$ on $\Fc_\k$ in the following way.

Let us denote by $C_c(\partial^2 \G \times \R)$ the space of real-valued continuous functions with compact supports on $\partial^2 \G \times \R$.
For each $f \in C_c(\partial^2 \G \times \R)$, 
let
\[
\wbar f(\x, \y, t):=\sum_{x \in \G}f(x\cdot (\x, \y, t)),
\]
and since $\wbar f$ is well-defined and $(\G, \k)$-invariant, 
it is identified with a function (which we denote by the same symbol) on $\Fc_\k$.
It has been shown that
for every $\G$-invariant Radon measure $\L$ on $\partial^2 \G$,
there exists a unique finite Borel measure $m$ on $\Fc_\k$ such that
\begin{equation*}
\int_{\partial^2 \G \times \R}f\,d\L\otimes dt=\int_{\Fc_\k}\wbar f\,dm \quad \text{for all $f \in C_c(\partial^2 \G \times \R)$},
\end{equation*}
and furthermore, $m$ is invariant under the topological flow $\{\F_t\}_{t \in \R}$ on $\Fc_\k$ \cite[Lemma 3.4]{TopFlows}.
Note that $m$ has a positive mass on every non-empty open set if $\L \otimes dt$ does.
We normalize $\L$ so that $m$ is a probability measure.

The following theorem is known in \cite[Theorem 2.5]{KaimanovichCrelle} through a generalization of the classical Hopf argument (see also \cite{BaderFurman}).
Let us state it in the form we use in the present discussion.

\begin{theorem}\label{Thm:Hopf}
Assume that there are finite Borel measures $\m$ and $\check \m$ on $\partial \G$ and a $\G$-invariant Radon measure $\L$ on $\partial^2 \G$ in the measure class of $\check \m \otimes \m$.
Let $m$ be the Borel probability measure $m$ on $\Fc_\k$ invariant under the topological flow $\{\Phi_t\}_{t \in \R}$ associated with $\L$.
For every $f \in L^1(\Fc_\k, m)$, we have
\[
\lim_{T \to \infty}\frac{1}{T}\int_0^T f\circ \F_t(x)\,dt=\int_{\Fc_\k}f\,dm, \quad \text{for $m$-almost every $x \in \Fc_\k$}.
\]
In particular, the measure $m$ is ergodic under the topological flow $\{\F_t\}_{t \in \R}$.
Moreover, $\L$ is ergodic with respect to the $\G$-action on $\partial^2 \G$, i.e., for any $\G$-invariant Borel set in $\partial^2 \G$ either it or the complement has $\L$-measure null.
\end{theorem}

\proof
The proof proceeds as an adaptation to \cite[Theorem 3.6 and Corollary 3.7]{TopFlows};
we reproduce the main argument for the sake of convenience.
For every $f$ in $L^1(\Fc_k, m)$,
by the Birkhoff ergodic theorem,
there exists the limit
\[
f_\infty(x):=\lim_{T \to \infty}\frac{1}{T}\int_0^T f\circ \Phi_t(x)\,dt,
\]
for $m$-almost every $x$ in $\Fc_\k$.
The convergence takes place also in $L^1(\Fc_\k, m)$.
We show that $f_\infty$ is constant $m$-almost everywhere.
Since $\|f_\infty\|_1 \le \|f\|_1$ by the Fatou lemma, and $\Fc_\k$ is a compact metrizable space,
it suffices to prove the claim for continuous functions.
Letting $\pi: \partial^2 \G \times \R \to \Fc_\k$ be the quotient map,
we define $\tilde f_\infty:=f_\infty \circ \pi$.
Note that $\tilde f$ is $\{\wt \F_t\}_{t \in \R}$-invariant, whence $\tilde f_\infty$ is defined on $\partial^2 \G$.
By the assumption, $\L$ and $\check \m \otimes \m$ are mutually absolutely continuous on $\partial^2 \G$,
and we have that for $\check \m \otimes \m$-almost every $(\x, \y) \in \partial^2 \G$ and for all $s \in \R$,
\[
\lim_{T \to \infty}\frac{1}{T}\int_0^T (f\circ \pi)(\wt \F_t(\x, \y, s))\,dt=\tilde f_\infty(\x, \y).
\]
Let us take any compact set $K$ in $\partial^2 \G$.
There exists a constant $c_K$ such that $\langle \x|\y\rangle_o \le c_K$ for all $(\x, \y) \in K$.
For $\x \in \partial \G$, letting $x$ be a point on a rough geodesic from $o$ toward $\x$ at distance $T$,
we have 
\[
\hat q(x^{-1}\y, x^{-1}\y') \le C e^{-T}\hat q(\y, \y') \quad \text{for all $(\x, \y), (\x, \y') \in K$},
\]
where $C$ is some positive constant and $\hat q$ denotes the quasi-metric associated to $\hat d$ (cf.\ Section \ref{Sec.hyp}).
Since $f\circ \pi$ is $\G$-invariant and uniformly continuous,
for $\check \m \otimes \m$-almost every $(\x, \y), (\x, \y') \in K$,
we have $\tilde f_\infty(\x, \y)=\tilde f_\infty(\x, \y')$.
Noting that $\partial^2 \G$ is $\s$-compact, we deduce that for $\check \m$-almost every $\x$ in $\partial\G$,
the function $\tilde f_\infty(\x, \cdot\,)$ is constant $\m$-almost everywhere.
Similarly for $\m$-almost every $\y$ in $\partial \G$,
the function $\tilde f_\infty(\,\cdot, \y)$ is constant $\check \m$-almost everywhere.
Therefore $\tilde f_\infty$ is constant $\check \m \otimes \m$-almost everywhere, implying that $f_\infty$ is constant $m$-almost everywhere.
\qed

\subsection{Local intersection numbers}\label{Sec:localint}

Let $\p$ be a $\G$-invariant function on $\G \times \G$ satisfying \eqref{Eq:RG}. We then have that
\[
\p(x, z)=\p(x, y)+\p(y, z)+O(1),
\]
for any $x, y, z$ aligned on a rough geodesic in this order in $\G$.
For $d \in \Dc_\G$ and for $\x \in \partial \G$,
let
\[
\tau_{\inf}^\p(\x):=\liminf_{t \to \infty}\frac{\p(\g(0), \g(t))}{d(\g(0), \g(t))}
\quad \text{and} \quad
\tau_{\sup}^\p(\x):=\limsup_{t \to \infty}\frac{\p(\g(0), \g(t))}{d(\g(0), \g(t))},
\]
where $\g$ is a quasi-geodesic converging to $\x$ in $\partial \G$.
Note that since $\p$ is $\G$-invariant, for each $\x$ in $\partial \G$ these two values are independent of the choice of quasi-geodesics toward it, or of the starting points by the Morse lemma.
Moreover, they are $\G$-invariant.
If $\tau_{\inf}^\p(\x)=\tau_{\sup}^\p(\x)$,
then we define the value $\tau^\p(\x)$ and call the {\it local intersection number} of $\p$ relative to $d$ at $\x$.

\begin{lemma}\label{Lem:localint}
Let $\p$ be a $\G$-invariant function on $\G \times \G$ satisfying \eqref{Eq:RG} and $d$ be a hyperbolic metric in $\Dc_\G$.
If $\L$ is a $\G$-invariant Radon measure which is ergodic on $\partial^2 \G$,
then the local intersection number $\t^\p(\x_+)$ of $\p$ relative to $d$ exists and is constant for $\L$-almost every $(\x_-, \x_+) \in \partial^2\G$.

In particular,
if there exist finite Borel measures $\m$ and $\check \m$ on $\partial \G$ and
a $\G$-invariant Radon measure $\L$ on $\partial^2 \G$ in the measure class of $\check \m \otimes \m$,
then the local intersection number $\t^\p(\x)$ of $\p$ relative to $d$ exists at $\m$-almost every $\x$, i.e., 
\[
\t_{\inf}^\p(\x)=\t_{\sup}^\p(\x) \quad \text{for $\m$-almost every $\x$ in $\partial \G$},
\]
and further $\t^\p(\x)$ is constant $\m$-almost everywhere on $\partial \G$.
\end{lemma}

\proof
Taking a Borel fundamental domain $F$ in $\partial^2 \G \times \R$ for the $(\G, \k)$-action and 
a measurable section $\i: \Fc_\k \to F$,
let $\tilde w:=\i(w)$ for $w \in \Fc_\k$, and $\tilde w_t:=\wt \F_t(\tilde w)$ for $t \in \R$.
Define
\[
c_\ast^\p(s, t):=\p(\ev(\tilde w_s), \ev(\tilde w_t)) \quad \text{for $s, t \in \R$}.
\]
Let $m$ denote the Borel probability measure on $\Fc_\k$ associated with $\L$ invariant under $\{\F_t\}_{t \in \R}$.
Note that $m$ is ergodic under the flow by the assumption,
and $c_\ast^\p(s, t)$ is additive modulo a uniform constant such that $\sup_{|s-t| \le 1}|c_\ast^\p(s, t)| \le C$ for a constant $C>0$ by the $\G$-invariance of $\p$.
Thus the Kingman subadditive ergodic theorem shows that
there exists a constant $\t_\ast^\p$ such that
\[
\lim_{t \to \infty}\frac{1}{t}c_\ast^\p(0, t)=\t_\ast^\p \quad \text{for $m$-almost everywhere on $\Fc_\k$}.
\]
Applying the above argument to the case when $\p=d$,
we show that there exists a positive constant $\t_\ast^d>0$ such that
\[
\lim_{t \to \infty}\frac{1}{t}c_\ast^d(0, t)=\t_\ast^d \quad \text{for $m$-almost everywhere on $\Fc_\k$},
\]
where we have that $\t_\ast^d>0$ since $d$ and $\hat d$ are quasi-isometric.
Therefore for $m$-almost every point $w$ in $\Fc_\k$,
we have that
\[
\lim_{t \to \infty}\frac{c_\ast^\p(0, t)}{c_\ast^d(0, t)}=\frac{\t_\ast^\p}{\t_\ast^d}.
\]
Noting that $\L\otimes dt=\sum_{x \in \G}x_\ast \i_\ast m$,
we conclude the first claim.
Concerning the second claim, since $\L$ is in the measure class of $\check \m \otimes \m$ by the assumption,
$\L$ is ergodic by Theorem \ref{Thm:Hopf}, and thus
we conclude that the local intersection number $\t^\p(\x)$ of $\p$ relative to $d$ exists and is a constant $\m$-almost everywhere in $\partial^2 \G$.
\qed

\section{Symbolic coding}\label{Sec:symbolic}

\subsection{Automatic structures and shift spaces}

Fix a finite symmetric set of generators $S$ of $\G$, i.e., $S=S^{-1}$.
Associated with $S$,
there exists a finite state automaton $\Ac=(\Gc, w, S)$
where $\Gc=(V, E, s_\ast)$ is a finite directed graph with a distinguished vertex $s_\ast$ (the \textit{initial state}),
endowed with a labeling $w: E \to S$.
For a directed edge path $\o=(\o_0, \o_1, \dots, \o_{n-1})$ (where the terminus of $\o_i$ coincides with the origin of $\o_{i+1}$)
there is an associated path in the Cayley graph $\Cay(\G, S)$ beginning at the identity: the path corresponds to
\[
w(\o):=(o, w(\o_0), w(\o_0)w(\o_1), \dots, w(\o_0)\cdots w(\o_{n-1})).
\]
Let 
\[
w_\ast(\o):=w(\o_0)\cdots w(\o_{n-1}).
\]
We say that $\Ac=(\Gc, w, S)$ is a {\it strongly Markov automatic structure}
if
\begin{itemize}
\item[(1)] for each vertex $v$ in $\Gc$ there exists a directed path from $s_\ast$ to $v$,
\item[(2)] for each directed path $\o$ in $\Gc$ the associated path $w(\o)$ is a geodesic in $\Cay(\G, S)$, and
\item[(3)] the map $w_\ast$ defines a bijection between the set of directed paths from $s_\ast$ in $\Gc$ and $\G$.
\end{itemize}

It is known that every hyperbolic group admits a strongly Markov automatic structure for any finite symmetric set of generators (cf.\ \cite[Section 3.2]{Calegari}).
Let us fix an automaton $\Ac=(\Gc, w, S)$ for $(\G, S)$.
We denote by $\SS^\ast$ (resp.\ by $\SS^+$) the set of finite directed paths (resp.\ the set of semi-infinite paths $(\o_i)_{i=0}^\infty$) in $\Gc$,
where paths are not necessarily starting from $s_\ast$ in both cases.
Letting $\wbar \SS^+:=\SS^\ast \cup \SS^+$,
we extend and define
\[
w_\ast: \wbar \SS^+ \to \G \cup \partial \G, 
\]
by assigning to each element of $\wbar \SS^+$ the terminus of the corresponding geodesic segment or ray issuing from $o$ in $\Cay(\G, S)$.

Let us define the space of bilateral directed paths in $\Gc$ by
\[
\SS:=\{(\o_i)_{i \in \Z} \ : \ A(\o_i, \o_{i+1})=1 \ \text{for all $i \in \Z$}\},
\]
where $A=(A(e, e'))_{e, e' \in E}$ is the adjacency matrix on edges of $\Gc$,
i.e., $A(e, e')=1$ if the terminus of $e$ is the origin of $e'$ and $0$ otherwise.
We define the shift 
\[
\s: \SS \to \SS, \quad \s(\o_i)_{i \in \Z}:=(\o_{i+1})_{i \in \Z}.
\]
Note that $(\SS, \s)$ is a subshift of finite type over alphabets $E$.
We define a metric in $\SS$ by setting
\[
d_\SS(\o, \o'):=
\begin{cases}
\exp(-n) &\text{where $n:=\sup\{m \ge 0 \ : \ \o_i=\o'_i \ \text{for all $|i| \le m$}\}$ if $\o \neq \o'$},\\
0 &\text{if $\o=\o'$}.
\end{cases}
\]
Each bilateral path $\o$ yields a geodesic line passing through $o$ in $\Cay(\G, S)$:
\[
(\x_n(\o))_{n \in \Z}:=(\dots, w(\o_{-1})^{-1}w(\o_{-2})^{-1}, w(\o_{-1})^{-1}, o, w(\o_0), w(\o_0)w(\o_1), \dots),
\]
where $\x_0(\o):=o$.
Abusing the notation, let us define the map which assigns the pair of two (distinct) extreme points in $\partial \G$ to each bilateral path by
\[
w_\ast: \SS \to \partial^2 \G, \quad \o \mapsto (\x_-(\o),\x_+(\o)),
\]
where $\x_-(\o)$ and $\x_+(\o)$ are $[\{\x_{-n}(\o)\}_{n=0}^\infty]$ and $[\{\x_n(\o)\}_{n=0}^\infty]$, respectively.
Note that the map $w_\ast$ is Lipschitz continuous with respect to the associated quasi-metric by the hyperbolicity of $\Cay(\G, S)$.

\subsection{Suspension flows}\label{Sec:suspension}

Given the cocycle $\k: \G \times \partial^2 \G \to \R$ defined by \eqref{Eq:top-cocycle},
let
\[
\tilde \k: \SS \to \R, \quad \o\mapsto \k(w(\o_0)^{-1}, w_\ast(\o)).
\]
We will write
\[
S_n \tilde \k(\o):=\sum_{i=0}^{n-1}\tilde \k(\s^i \o) \quad \text{for $n \in \Z_{\ge 1}$}.
\]
There exist an $N \in \Z_{\ge 1}$ and positive constants $c_1, c_2>0$ such that
$c_1\le S_N\tilde \k(\o)\le c_2$ for all $\o \in \SS$ \cite[Lemma 4.3]{TopFlows}.
Fix such an $N$ and let
\[
r_N:=S_N \tilde \k.
\]
Let us endow $\SS \times \R$ with the product topology.
We define
\[
\Sus(\SS, r_N):=\(\SS \times \R\)/\sim,
\]
where 
$(\s^N \o, t) \sim (\o, t+r_N(\o))$ for $(\o, t) \in \SS \times \R$,
and endow it with the quotient topology.
Further defining the continuous $\R$-action on $\Sus(\SS, r_N)$ by
\[
\s_t[\o, s]:=[\o, t+s] \qquad \text{for $t \in \R$},
\]
we obtain the flow $\{\s_t\}_{t \in \R}$ on $\Sus(\SS, r_N)$ as a {\it suspension flow} over $(\SS, \s^N)$.

The natural projection
\[
\widetilde \Pi: \SS \times \R \to \partial^2\G \times \R, \quad (\o, t) \mapsto (w_\ast(\o), t)
\]
induces the map
\[
\Pi: \Sus(\SS, r_N) \to \Fc_\k, \quad [\o, t]\mapsto [w_\ast(\o), t],
\]
as a consequence of the cocycle relation for $\k$.
It has been shown that $\Pi$ is continuous and equivariant with flows, i.e.,
$\Pi \circ \s_t=\F_t \circ \Pi$ for all $t \in \R$,
furthermore, $\Pi$ is surjective and the cardinality of each fiber is uniformly bounded \cite[Proposition 4.4]{TopFlows}.

\begin{lemma}\label{Lem:coboundary}
For the cocycle $\k: \G \times \partial^2 \G \to \R$ defined in \eqref{Eq:top-cocycle} by the Busemann cocycle $\hat \b_o$,
if we define 
\[
\tilde \k: \SS \to \R,\quad
\o \mapsto \k(w(\o_0)^{-1}, w_\ast(\o)),
\]
and further
\[
\P_B: \SS \to \R, \quad \o \mapsto \hat \b_o(w(\o_0), \xi_+(\o))  
\quad
\text{and}
\quad
U: \SS \to \R, \quad \o \mapsto \langle \xi_-(\o)|\xi_+(\o)\rangle_o,
\]
then
the functions $\tilde \k$, $\P_B$ and $U$ are H\"older continuous on $(\SS, d_\SS)$, i.e.,
there exist constants $L \ge 0$ and $\a>0$ such that
\[
|\tilde \k(\o)-\tilde \k(\o')| \le L d_\SS(\o, \o')^\a \quad \text{for $\o, \o' \in \SS$},
\]
and the same inequalities hold for $\P_B$ and $U$.
Furthermore we have
\[
\tilde \k=-\P_B+U\circ \s-U.
\]
\end{lemma}

\proof
The H\"older continuity follows since $\hat \b_o$ is defined by a strongly hyperbolic metric in $\Dc_\G$.
The second claim follows the definition of the cocycle $\k$; for the details see \cite[Lemma 5.1]{TopFlows}.
\qed

\subsection{Thermodynamic formalism}

Let us consider the set of unilateral (finite or infinite) paths $\wbar \SS^+$ (including the empty path $\emptyset$) and define the shift $\s: \wbar \SS^+ \to \wbar \SS^+$ by deleting the initial edge of each path.
Here we understand that $\s(\o)=\emptyset$ for paths $\o$ of length $1$ (i.e., edges).
Further we define the metric
\[
d_{\wbar \SS^+}(\o, \o'):=
\begin{cases}
\exp(-n)	&\text{where $n:=\sup\{m \ge 0 \ : \ \o_i=\o'_i \ \text{for all $0 \le i \le m$}\}$ if $\o \neq \o'$},\\
0 			&\text{if $\o=\o'$},
\end{cases}
\]
similarly to the case of bilateral paths $\SS$.

For each real-valued H\"older continuous function (potential) $\P:\wbar \SS^+ \to \R$ relative to $d_{\wbar \SS^+}$,
the transfer operator is defined on continuous functions $f$ on $\wbar \SS^+$ by
\[
\Lc_\P f(\o)=\sum_{\s(\o')=\o}e^{\P(\o')}f(\o').
\]
We say that a finite directed graph is \textit{recurrent} if there exists a directed path between any pair of vertices.
A {\it component} of a finite directed graph is a maximal subgraph which is recurrent.
The underlying graph $\Gc$ of an automaton $\Ac$ can admit several distinct components which are not singletons.
This causes difficulties when we analyze the spectral properties of transfer operators on the associated shift space.

Each component $\Cc$ of $\Gc$ has a period $p_\Cc \ge 1$,
i.e., the length of any closed path in $\Cc$ is a multiple of $p_\Cc$.
If the period $p_\Cc$ is $1$, then $\Cc$ is {\it topologically mixing}, i.e., there exists a directed path between any pair of vertices with length $n$
for every $n \ge N$ for some integer $N$.
If $p_\Cc>1$,
then the set of vertices $V(\Cc)$ of $\Cc$ admits a {\it cyclic decomposition}
$V(\Cc)=\bigsqcup_{i \in \Z/p_\Cc\Z}V_i$ where $V_i$ are disjoint subsets of vertices such that any edge whose origin is in $V_i$ has the terminus in $V_{i+1}$.
If we denote by $\wbar \SS^+_i$ the set of paths starting from a vertex in $V_i$ and the empty path,
then $\s: \wbar \SS^+_i \to \wbar \SS^+_{i+1}$
and the restriction of $\s^{p_\Cc}$ on $\wbar \SS^+_i$ is a topological mixing subshift of finite type for each $i \in \Z/p_\Cc\Z$.

\def\Pr{{\rm Pr}}
For each component $\Cc$, define the transfer operator $\Lc_\Cc$ by restricting $\P$ on the set of paths in $\Cc$.
Since $\Cc$ is recurrent, $\Lc_\Cc$ has finitely many eigenvalues of maximal modulus $e^{\Pr_\Cc(\P)}$
where $\Pr_\Cc(\P)$ is called the \textit{pressure}.
Let
\[
\Pr(\P):=\max_\Cc \Pr_\Cc(\P),
\]
where $\Cc$ runs over all components in $\Gc$.
A component $\Cc$ is called {\it maximal} for $\P$ if $\Pr(\P)=\Pr_\Cc(\P)$.

\begin{lemma}\label{Lem:Pr}
For $d$ in $\Dc_\G$, let $\p$ be a $\G$-invariant tempered potential with exponent $0$ relative to $d$; see Definition \ref{Def:tempered-potential}.
If for a strongly Markov automatic structure $\Ac=(\Gc, w, S)$,
the corresponding unilateral shift space $(\wbar \SS^+, \s)$ admits a H\"older continuous potential $\P$ such that
\[
S_n\P(\o)=\sum_{i=0}^{n-1}\P(\s^i(\o))=-\p(o, w_\ast (\o))+O(1) \quad \text{for all $\o=(\o_0, \dots, \o_{n-1}) \in \SS^*$},
\]
then $\Pr(\P)=0$.
\end{lemma}

\proof
The proof is based on the following estimate: For all $n$,
\begin{equation}\label{Eq:L}
\Lc_\P^n \1_{[E_\ast]}(\emptyset)=\sum_{|x|_S=n}e^{-\p(o, x)+O(1)},
\end{equation}
where $[E_\ast]$ denotes the set of paths in $\wbar \SS^+$ starting at the initial state.
Since $\p$ has exponent $0$ relative to $d$,
for an associated finite Borel measure $\m$ satisfying \eqref{Eq:QC} with respect to $\b_o^\p$,
there exist positive constants $C, R>0$ such that
\[
C^{-1} \le \exp(\p(o, x))\cdot \m(\Oc(x, R)) \le C \quad \text{for all $x \in \G$}.
\]
This shows that \eqref{Eq:L} is uniformly bounded from below and from above
(cf.\ \cite[Lemmas 2.8 and 4.7]{CT}). This yields the claim since $\Pr(\P)$ is the exponential growth rate
of \eqref{Eq:L} in $n$.
\qed

\begin{example}\label{Ex:hatBpotential}
For the potential associated with the Busemann function $\hat \b_o$ by
$\P_B:\wbar \SS^+ \to \R$, $\o \mapsto \hat \b_o(w(\o_0), w_\ast(\o))$,
we note that $\P_B$ is H\"older continuous on $\wbar \SS^+$, 
and for $\o \in \SS^+$ (semi-infinite paths),
\begin{align*}
S_n\P_B(\o)=\sum_{i=0}^{n-1}\P_B\circ \s^i(\o)&=\hat \b_o(w_\ast(\o_0, \dots, \o_{n-1}), w_\ast(\o))\\
&=-\hat d(o, w_\ast(\o_0, \dots, \o_{n-1}))+O(1),
\end{align*}
where the implicit constant depends only on the metrics $\hat d$ and $d_S$ by \eqref{Eq:RG}.
For the exponential growth rate $\hat v$ for $\hat d$,
noting that $\hat v \hat d$ is a $\G$-invariant tempered potential with exponent $0$ relative to $\hat d$,
we have that $\Pr(\hat v\P_B)=0$ by Lemma \ref{Lem:Pr} (cf.\ \cite[Lemma 5.4]{TopFlows}).
\end{example}

For each recurrent component $\Cc$, let $\wbar \SS^+_\Cc$ be the set of unilateral paths staying in $\Cc$ all the time,
and denote by $(\wbar \SS^+_\Cc, \s)$ the associated subshift.
Let $C^\a(\wbar \SS^+_\Cc)$ be the space of real-valued H\"older continuous functions on $\wbar \SS^+_\Cc$ with a fixed exponent $\a \in (0, 1)$,
and $\|\cdot\|_{\a}$ denotes the corresponding H\"older norm.

\begin{theorem}\label{Thm:TF}
If $\P$ is a function in $C^\a(\wbar \SS^+_\Cc)$
and a recurrent component $\Cc$ has the period $p_\Cc$ with a cyclic decomposition $\Cc=\bigsqcup_{i \in \Z/p_\Cc\Z}\Cc_i$,
then there exist positive constants $C, \e_0>0$, functions $h_i$ in $C^\a(\wbar \SS^+_\Cc)$ and finite measures $\lambda_i$ such that the following hold:
\begin{itemize}
\item[(1)] 
For all $f$ in $C^\a(\wbar \SS^+_\Cc)$ and all $n \ge 0$,
\[
\left\|\Lc_\P^n f -e^{\Pr_\Cc(\P)n}\sum_{i \in \Z/p_\Cc \Z}\(\int_{\wbar \SS^+_\Cc}f\,d\lambda_{i-n \,{\rm mod}\, p_\Cc}\)h_i\right\|_{\a}\le C\|f\|_{\a}\cdot e^{(\Pr_\Cc(\P)-\e_0)n},
\]
and the measure $\m_\Cc=\sum_{i \in \Z/p_\Cc\Z}h_i\lambda_i$ is shift-invariant and ergodic.
\item[(2)] For all $\f$ in $C^\a(\wbar \SS^+_\Cc)$ with small enough norm,
there exist functions $h_i^\f$ in $C^\a(\wbar \SS^+_\Cc)$ and finite measures $\lambda_i^\f$ with the same supports as $h_i$ and $\lambda_i$ respectively such that
for all $f$ in $C^\a(\wbar \SS^+_\Cc)$ and all $n \ge 0$,
\[
\left\|\Lc_{\P+\f}^n f-e^{\Pr_\Cc(\P+\f)n}\sum_{i \in \Z/p_\Cc \Z}\(\int_{\wbar \SS^+_\Cc}f\,d\lambda_{i-n \,{\rm mod}\,p_\Cc}^\f\)h_i^\f\right\|_{\a} \le C\|f\|_{\a}\cdot e^{(\Pr_\Cc(\P)-\e_0)n}.
\]
Furthermore, the maps $\f \mapsto \Pr_\Cc(\P+\f)$, $\f \mapsto h_i^\f$ and $\f \mapsto \lambda_i^\f$
from $C^\a(\wbar \SS^+_\Cc)$ to $\R$, $C^\a(\wbar \SS^+_\Cc)$ and the dual of $C^\a(\wbar \SS^+_\Cc)$ respectively,
are real analytic in a small neighborhood of $0$ in $C^\a(\wbar \SS^+_\Cc)$, and we have that
\[
\frac{\partial}{\partial s}\Big|_{s=0}\Pr_\Cc(\Psi+s\f)=\int_{\wbar \SS_\Cc^+}\f\,d\m_{\Cc},
\]
where $\m_\Cc$ is determined by $\Psi$ in the first statement (1) and normalized as the probability measure.
\end{itemize}
\end{theorem}

The proof can be reduced to the case when $\Cc$ is topologically mixing; see \cite[Theorem 2.2, Propositions 4.6 and 4.10]{ParryPollicott}.

Let $\SS_\Cc$ be the set of bilateral paths in $\Cc$ and $(\SS_\Cc, \s)$ be the associated shift space.
Since $\Cc$ is recurrent, this shift space is {\it topologically transitive}, i.e., for any two non-empty open sets $U$ and $V$, there exists $n \in \Z$ such that $U\cap \s^n V \neq \emptyset$.
Furthermore, for a function $\P$ on $\wbar \SS^+$ we define (and use the same symbol) the induced potential $\P$ on $\SS_\Cc$ by setting $\P(\o):=\P(\o_0, \o_1, \dots)$,
where the value depends only on the non-negative indices.
Let us denote by $\Mcc(\s, \SS_\Cc)$ the set of all $\s$-invariant Borel probability measures on $\SS_\Cc$, by $h(\s, \lambda)$ the measure theoretical entropy of $(\SS_\Cc, \s, \lambda)$ (see Section \ref{Sec:entropy}).
Let $[\o_0, \dots, \o_{n-1}]$ be cylinder sets in $\SS_\Cc$.

\begin{proposition}[Uniqueness of equilibrium state]\label{Prop:VP}
If $\P$ is a H\"older continuous function on $\SS_\Cc$,
then 
\[
\Pr_\Cc(\P)=\sup_{\lambda \in \Mcc(\s, \SS_\Cc)}\left\{h(\s, \lambda)+\int_{\SS_\Cc}\P\,d\lambda\right\},
\]
and the supremum is attained by a unique equilibrium state, i.e., a unique $\s$-invariant Borel probability measure $\m_\Cc$ on $\SS_\Cc$, and further there exist positive constants $c_1, c_2>0$ such that
\begin{equation}\label{Eq:Gibbs}
c_1 \le \frac{\m_\Cc[\o_0, \dots, \o_{n-1}]}{\exp\(-n\Pr_\Cc(\P)+S_n \P(\o)\)} \le c_2,
\end{equation}
for all $\o \in [\o_1, \dots, \o_{n-1}]$ and for all $n \in \Z_{\ge 1}$,
where $S_n\P:=\sum_{i=0}^{n-1}\P\circ \s^i$.
\end{proposition}

For the proof, see \cite[Theorem 3.5]{ParryPollicott} (cf.\ \cite[Proposition 5.7]{TopFlows}).

Recall that $\hat d$ is a strongly hyperbolic metric in $\Dc_\G$ and $\hat \b_o$ is the associated Busemann cocycle.
Let us assume that $\p$ is a $\G$-invariant tempered potential relative to $\hat d$ with exponent $\th$. 
If $\check \p$ is also a tempered potential, where $\check \p(x, y):=\p(y, x)$ for $x, y \in \G$,
then $\check \p$ is a $\G$-invariant tempered potential relative to $\hat d$ with the same exponent $\th$ (see Proposition \ref{Prop:PSpotential}).
Let $\m$ and $\check \m$ be finite Borel measures on $\partial \G$ satisfying \eqref{Eq:QC} with the quasi-cocycles $\b_o^\p +\th \hat \b_o$ and $\b_o^{\check \p}+\th \hat \b_o$ respectively,
i.e., for some positive constant $C>0$, we have
\[
C^{-1} \le \exp(\b_o^\p(x, \x) +\th \hat \b_o(x, \x))\cdot \frac{d x_\ast \m}{d \m}(\x) \le C,
\]
for $(x, \x) \in \G \times \partial \G$
and similarly for $\check \m$.

\begin{lemma}\label{Lem:codingPSpotential}
Let $\p$ be a $\G$-invariant tempered potential relative to $\hat d$ with exponent $\th$ and $\check \p(x, y)=\p(y, x)$ for $x, y \in \G$.
Assume that $\check \p$ is also a tempered potential and
there exists a H\"older continuous function $\P$ on $(\wbar \SS^+, \s)$ such that
\[
S_n \P(\o)=\sum_{i=0}^{n-1}\P(\s^i(\o))=(\b_o^\p+\th \hat\b_o)(w_\ast(\o_0, \dots, \o_{n-1}), w_\ast(\o)) \quad \text{for all $\o \in \wbar \SS^+$}.
\]
If for a recurrent component $\Cc$,
we have
\[
\Pr_\Cc(\P)=0,
\] 
then for the induced potential $\P$ on $\SS_\Cc$,
the corresponding equilibrium state $\m_\Cc$ satisfies the following:
\begin{itemize}
\item[(1)]
For finite Borel measures $\m$ and $\check \m$ on $\partial \G$ satisfying \eqref{Eq:QC} with quasi-cocycles $\b_o^\p+\th \hat \b_o$ and $\b_o^{\check \p}+\th \hat \b_o$, respectively,
there exists a constant $C>0$ such that
\[
w_\ast \m_\Cc \le C \check \m \otimes \m.
\]
\item[(2)]
If we further assume that there exists a $\G$-invariant Radon measure $\L_\p$ in the same measure class of $\check \m \otimes \m$ on $\partial^2 \G$,
then 
letting $m_\p$ be a finite Borel measure on $\Fc_\k$ invariant under the topological flow $\{\F_t\}_{t \in \R}$ associated with $\L_\p$ and
$\pi: \partial^2 \G\times \R \to \Fc_\k$ be the quotient map,
we have that
for any $T>0$, there exists a constant $C>0$ such that
\[
\pi_\ast(w_\ast \m_\Cc \otimes \Leb_{[0, T)}) \le C m_\p.
\]
\end{itemize}
\end{lemma}

\proof
First we show (1).
Let us consider shadows $\Oc(x, R)$ on $\Cay(\G, S)$ for some $R>0$.
For each $x, y \in \G$, and for positive integers $m, n >0$,
let $[\o_{-m}, \dots, \o_{n-1}]$ be any cylinder set such that
for all $\o$ in the cylinder set, we have
\begin{equation}\label{Eq:cylinder}
\x_{-m}(\o) \in B_S(x, R) \quad \text{and} \quad \x_n(\o) \in B_S(y, R),
\end{equation}
where $B_S(x, R)$ stands for the ball of radius $R$ centered at $x$ in $\Cay(\G, S)$.
By assumption, for all $\o \in [\o_{-m}, \dots, \o_{n-1}]$,
we have
\begin{align*}
S_n \P(\o)	&=-(\p+\th \hat d)(o, w_\ast(\o_0, \dots, \o_{n-1}))+O(1)\\
			&=-(\p+\th \hat d)(o, y)+O(1),
\end{align*}
since $\p$ and $\hat d$ satisfy \eqref{Eq:RG}, and further
\begin{align*}
S_m \P(\s^{-m}\o)	&=-(\p+\th \hat d)(o, w_\ast(\o_{-m}, \dots, \o_{-1}))+O(1) \\
					&=-(\p+\th \hat d)(o, x^{-1})+O(1)=-(\check\p+\th \hat d)(o, x)+O(1),
\end{align*}
where the last equality holds by the $\G$-invariance of $\p$ and the definition of $\check \p$.
By Proposition \ref{Prop:PSpotential}, we have
\[
\m(\Oc(y, R))=\exp(-\p(o, y)-\th \hat d(o, y)+O(1))
\]
and
\[
\check \m(\Oc(x, R))=\exp(-\check \p(o, x)-\th \hat d(o, x)+O(1)).
\]
For the unique equilibrium state $\m_\Cc$ on $(\SS_\Cc, \s)$ with the potential $\P$, we have that
\begin{align*}
w_\ast \m_\Cc(\Oc(x, R) \times \Oc(y, R))	
&=\m_\Cc\(\{\o \in \SS \ : \ w_\ast(\o) \in \Oc(x, R) \times \Oc(y, R)\}\)\\
& \le \sum \m_\Cc[\o_{-m}, \dots, \o_{n-1}],
\end{align*}
where the summation runs over all cylinder sets satisfying \eqref{Eq:cylinder} (and $m, n$ may vary).
Note that the number of those cylinder sets is at most $\# B_S(x, R) \cdot \# B_S(y, R)$.
Since $\m_\Cc$ satisfies \eqref{Eq:Gibbs}, for all $\o$ in such cylinder set $[\o_{-m}, \dots, \o_{n-1}]$,
we have
\[
\m_\Cc[\o_{-m}, \dots, \o_{n-1}]
=\exp(-(n+m)\Pr_\Cc(\P)+S_{n+m}\P(\s^{-m}\o)+O(1)).
\]
Hence if $\Pr_\Cc(\P)=0$, then
\begin{align*}
\m_\Cc[\o_{-m}, \dots, \o_{n-1}]	
&=\exp(S_m\P(\s^{-m}\o)+S_n\P(\o) +O(1))\\
& \le C\check \m(\Oc(x, R))\cdot \m(\Oc(y, R)),
\end{align*}
where $C$ is a positive constant independent of $x$ and $y$,
and thus
\[
w_\ast \m_\Cc(\Oc(x, R) \times \Oc(y, R)) \le C(\# B_S(o, R))^2 \check \m(\Oc(x, R))\cdot \m(\Oc(y, R)),
\]
where we have used $\#B_S(o,R) = \#B_S(x,R)$ due to the $\Gamma$-invariance of the word metric.
Replacing shadows by balls in the boundary (see \eqref{Eq:shadow-ball}), we obtain for some positive constants $C, C'>0$,
\[
w_\ast \m_\Cc(B(\x, r) \times B(\y, s)) \le C \check \m(B(\x, C' r))\cdot \m(B(\y, C' s))
\]
for all $(\x, \y) \in \partial^2 \G$ and all $r, s >0$.
We note that $\check \m$ and $\m$ are doubling (Remark \ref{Rem:doubling}).
The rest follows as in \cite[Lemma 5.8]{TopFlows}.

Second we show (2).
Note that $w_\ast \m_\Cc \otimes \Leb_{[0, T)}$ is compactly supported.
Letting
\[
M_0:=\sum_{x \in \G}x.
\(w_\ast \m_\Cc\otimes\Leb_{[0, T)}\),
\]
we have that $M_0$ is $(\G, \k)$-invariant and moreover, it is Radon since 
the $(\G, \k)$-action is properly discontinuous and cocompact.
Further since $w_\ast \m_\Cc \le C \check \m\otimes \m$ for a constant $C>0$ on $\partial^2 \G$ by the first claim,
and restricting on the compact support of $w_\ast \m_\Cc$,
the measures $\check \m\otimes \m$ and $\L_\p$ mutually dominate each others up to multiplicative constants,
we have a constant $C'>0$ such that 
$M_0 \le C' \L_\p\otimes dt$.

If we take a Borel fundamental domain $F$ and a measurable section $\i:\Fc_\k \to F \subset \partial^2 \G \times \R$,
then
\[
\sum_{x \in \G}x.\i_\ast \pi_\ast\(w_\ast \m_\Cc \otimes \Leb_{[0, T)}\)=\sum_{x \in \G}x.\(w_\ast \m_\Cc\otimes \Leb_{[0, T)}\)=M_0.
\]
Hence we have
\[
\int_{\Fc_\k}\wbar f\,d\(\pi_\ast(w_\ast \m_\Cc \otimes \Leb_{[0, T)})\)=\int_{\partial^2 \G \times \R}f\,dM_0 \qquad \text{for all $f \in C_c(\partial^2 \G \times \R)$},
\]
and thus for all non-negative valued function $f \in C_c(\partial^2 \G \times \R)$,
\[
\int_{\Fc_\k}\wbar f\,d\(\pi_\ast(w_\ast \m_\Cc\otimes \Leb_{[0, T)})\) \le C' \int_{\Fc_\k}\wbar f\,dm_\p.
\]
For any continuous functions $\f$ on $\Fc_\k$, there exists $f \in C_c(\partial^2 \G \times \R)$
such that $\f=\wbar f$, further if $\f \ge 0$, then one may take $f \ge 0$ by invoking the Urysohn lemma.
Therefore for all continuous function $\f\ge 0$ on $\Fc_\k$,
we have that
\[
\int_{\Fc_\k}\f\,d\(\pi_\ast(w_\ast \m_\Cc \otimes \Leb_{[0, T)})\) \le C' \int_{\Fc_\k}\f\,dm_\p,
\]
and noting that $\Fc_\k$ is a compact metrizable space (Section \ref{Sec:TopFlows}), 
we obtain the claim up to changing constants if required.
\qed

\begin{proposition}\label{Prop:ac}
Fix a strongly hyperbolic metric $\hat d$ in $\Dc_\G$ with the exponential growth rate $\hat v$.
Letting $\P_B(\o):=\hat \b_o(w(\o_0), w_\ast(\o))$ for $\o \in \wbar \SS^+$,
we take a maximal component $\Cc$ for $\hat v\P_B$.
For a finite Borel measure $m_{\hat d}$ on $\Fc_\k$ invariant under the topological flow $\{\F_t\}_{t \in \R}$ associated with $\hat d$,
we have that for any $T>0$, there exists a constant $C>0$ such that
\[
\pi_\ast(w_\ast \m_\Cc \otimes \Leb_{[0, T)}) \le C m_{\hat d}
\]
where $\m_\Cc$ is the unique equilibrium state for the induced potential $\hat v\P_B$ on $(\SS_\Cc, \s)$ and $\pi: \partial^2 \G\times \R \to \Fc_\k$ is the quotient map.
\end{proposition}

\proof
Note that $\Pr(\hat v \P_B)=0$
 (cf.\ Example \ref{Ex:hatBpotential}),
 whence for any maximal component $\Cc$ for $\hat v \P_B$, we have $\Pr_\Cc(\hat v \P_B)=0$.
 Applying Lemma \ref{Lem:codingPSpotential} to the case when $\p\equiv 0$ and $\th=\hat v$,
 we obtain the claim.
\qed

\subsection{An enhanced coding}\label{Sec:enhanced}

Fix a strongly hyperbolic metric $\hat d$ in $\Dc_\G$ with the exponential growth rate $\hat v$
and fix an arbitrary maximal component $\Cc_0$ for $\hat v\P_B$.
The following proposition shows that one may take for $\Cc_0$ any component whose adjacency matrix has an eigenvalue of maximal modulus.

\begin{proposition}\label{Prop:component}
For any finite symmetric set of generators $S$, let $v_S$ be the exponential growth rate for the word metric, and $\Ac=(\Gc, w, S)$ be an arbitrary strongly Markov automatic structure for $(\G, S)$.
Let us consider the constant function $\P_S\equiv 1$ and $\P_B$ defined by the Busemann cocycle for a strongly hyperbolic metric $\hat d$ on the associated shift space $(\wbar \SS^+, \s)$.
A component $\Cc$ in $\Gc$ is maximal for $\hat v \P_B$ if $\Cc$ is maximal for $\P_S$.
In particular, a component whose adjacency matrix has an eigenvalue of maximal modulus $e^{v_S}$ is maximal for $\hat v \P_B$.
\end{proposition}

\proof
On the one hand, we recall that $\P_B$ is a H\"older continuous potential of the form:
\[
\P_B(\o)=\hat \b_o(w(\o_0), w_\ast(\o)) \quad \text{for $\o \in \wbar\SS^+$},
\]
where $\hat \b_o$ is the Busemann function associated with $\hat d$.
For all $n$, we have that
\[
\Lc_{\hat v\P_B}^n\1_{[E_\ast]}(\emptyset)=\sum_{x : \ |x|_S=n} \exp(-\hat v\hat d(o, x)+O(1)).
\]
If we define the series in $t$ by
\[
\hat\Pc(t):=\sum_{x \in \G}\exp(-\hat v\hat d(o, x)-t|x|_S) \quad \text{for $t \in \R$},
\]
then the divergence exponent $\hat \th_S$ is given by
$\hat \th_{S}=\Pr(\hat v\P_B)$ by Theorem \ref{Thm:TF}(1) (in fact, $\hat \th_S=0$, but we do not use this fact).
On the other hand, for any maximal component $\Cc$ for $\P_S$, which depends only on $\Ac=(\Gc, w, S)$,
let $\Gamma_\Cc$ be the set of group elements defined by the words corresponding to finite paths in $\Cc$.
Then there exists a finite set $B$ in $\G$ such that $B\G_\Cc B=\G$ (for the proof, see \cite[Lemma 4.19]{CT} and the references therein, where the main argument relies on \cite[Theorem 3]{arzlys}).
This implies that if we define
\[
\hat\Pc_\Cc(t):=\sum_{x \in \G_\Cc}\exp(-\hat v\hat d(o, x)-t|x|_S) \quad \text{for $t \in \R$},
\]
then 
\[
\hat \Pc_\Cc(t) \le \hat \Pc(t) \le e^{c(t)}(\#B)^2\hat \Pc_\Cc(t),
\]
where $c(t):=2\hat v\max_{x \in B}\hat d(o, x)+2|t|\max_{x \in B} |x|_S$,
showing that the divergence exponents for $\hat \Pc(t)$ and for $\hat \Pc_\Cc(t)$ in $t$ coincide.
Since 
\[
\hat \Pc_\Cc(t)=\sum_{n=1}^\infty\sum_{\o=(\o_0, \dots, \o_{n-1}) \ \text{in $\Cc$}}\exp(-\hat v S_n\P_B(\o)-t n),
\]
we have that $\hat \th_{S}=\Pr_\Cc(\hat v\P_B)$.
Therefore,
we have that $\Pr(\hat v\P_B)=\Pr_\Cc(\hat v\P_B)$,
i.e., any maximal component $\Cc$ for the constant potential $\P_S$ is maximal for $\hat v \P_B$.
The last assertion follows since a maximal component for the constant $1$ potential $\P_S$ has the adjacency matrix with the maximal eigenvalue in modulus $e^{v_S}$ by the Perron-Frobenius theorem.
\qed

\medskip

Let $\SS_0:=\SS_{\Cc_0}$.
Since $(\SS_0, \s)$ defines a subshift of $(\SS, \s)$,
we may define the suspension flow space $\Sus(\SS_0, r_0)$ where $r_0:=r_N|_{\SS_0}$ is the restriction of $r_N$ on $\SS_0$
as a subsystem of $\Sus(\SS, r_N)$.
We define 
\[
\Pi_0: \Sus(\SS_0, r_0) \to \Fc_\k,
\]
as the restriction of $\Pi: \Sus(\SS, r_N) \to \Fc_\k$.
Note that $\Pi_0$ is continuous, equivariant with the flows and the cardinality of each fiber is uniformly bounded (cf.\ Section \ref{Sec:suspension}).
Let us show that 
the map $\Pi_0: \Sus(\SS_0, r_0) \to \Fc_\k$ is surjective, and thus establish the following.

\begin{theorem}[Main coding theorem]\label{Thm:enhanced}
For any non-elementary hyperbolic group $\G$, 
let $\Fc_\k$ be the topological flow space associated with a cocycle $\k$.
There exists a topologically transitive subshift of finite type $(\SS_0, \s)$ with a positive H\"older continuous function $r_0$
such that a natural continuous map $\Pi_0$ from the associated suspension flow $\Sus(\SS_0, r_0)$ to $\Fc_\k$
is surjective, equivariant with the flows and the cardinality of each fiber is uniformly bounded.
\end{theorem}

\proof
It suffices to show that the map $\Pi_0$ is surjective.
Let $\Xc_0:=\Sus(\SS_0, r_0)$ and $\Xc:=\Sus(\SS, r_N)$ for the simplicity of notations.
Note that $\Xc_0$ is a compact subset of $\Xc$
and the image $\Pi_0(\Xc_0)$ is a compact subset of the compact Hausdorff space $\Fc_\k$.
Hence $\Pi_0(\Xc_0)$ is closed.
Since $\Pi_0$ is equivariant with the flows, the closed set $\Pi_0(\Xc_0)$ is invariant under the topological flow $\{\F_t\}_{t \in \R}$ in $\Fc_\k$.
Note that there is a natural bijection between the space of $\s^N$-invariant probability measures $\Mcc(\s^N, \SS_0)$ (where $N$ is as in Section \ref{Sec:suspension} and we use the notations therein in the following)
and the space of flow-invariant Borel probability measures $\Mcc(\{\s_t\}_{t \in \R}, \Xc_0)$ by
\[
\Mcc(\s^N, \SS_0) \to \Mcc(\{\s_t\}_{t \in \R}, \Xc_0), \quad \lambda \mapsto \n_\lambda=\frac{1}{\lambda\otimes dt(Z)}\lambda\otimes dt|_Z,
\]
where $Z:=\{(\o, t) \in \SS \times \R \ : \ 0\le t < r_0(\o)\}$ is identified with $\Xc_0$.
For the unique equilibrium state $\m_0$ on $\SS_0$ for $\hat v \P_B$,
we consider the corresponding flow-invariant probability measure $\n_0$ on $\Xc_0$ and define
\[
m_0:=\Pi_{0 \ast}\n_0.
\]
Note that $m_0$ is invariant under the topological flow $\{\F_t\}_{t \in \R}$ on $\Fc_\k$.

Let us show that $m_0$ is absolutely continuous with respect to $m_{\hat d}$.
If we define $T_0:=\sup_{\o \in \SS_0}r_0(\o)$,
then
\[
\widetilde \Pi_\ast \n_0 \le C w_\ast \m_0 \otimes \Leb_{[0, T_0)},
\]
for the positive constant $C:=1/\m_0\otimes dt(Z)>0$.
Further recalling that $\pi: \partial^2 \G \times \R \to \Fc_\k$ is the quotient map,
we have
\[
m_0=\Pi_{0 \ast} \n_0=\pi_\ast \widetilde \Pi_\ast \n_0 \le C \pi_\ast\(w_\ast \m_0\otimes \Leb_{[0, T_0)}\) \le C' m_{\hat d},
\]
for some positive constant $C'>0$ by Proposition \ref{Prop:ac}.
Therefore $m_0$ is absolutely continuous with respect to $m_{\hat d}$.

This implies that $m_0=m_{\hat d}$ since $m_0(\Fc_\k)=1$ and $m_{\hat d}$ is ergodic with respect to $\{\F_t\}_{t \in \R}$ (Theorem \ref{Thm:Hopf}).
Note that 
$m_{\hat d}(U)>0$ for all non-empty open sets $U$ in $\Fc_\k$ (see the second paragraph of Section \ref{Sec:measures-ergodicity}).
We also have that $m_0$ is supported on $\Pi_0(\Xc_0)$.
Suppose that $\Pi_0(\Xc_0)$ is a proper subset of $\Xc$;
then since $\Pi_0(\Xc_0)$ is closed in $\Fc_\k$ there exists a non-empty open set $U$ disjoint from $\Pi_0(\Xc_0)$, implying that $m_0(U)=0$,
which leads a contradiction to that $m_0(U)=m_{\hat d}(U)>0$.
Therefore $\Pi_0(\Xc_0)=\Fc_\k$, i.e., $\Pi_0$ is surjective.
\qed

\section{Applications through the enhanced coding}\label{Sec:application-coding}

Let us fix $(\SS_0, \s)$ and $\Pi_0: \Sus(\SS_0, r_0) \to \Fc_\k$ in Theorem \ref{Thm:enhanced}.

\subsection{Invariant measures via the enhanced coding}

Let us denote by $\{\s_t\}_{t \in \R}$ the suspension flow on $\Sus(\SS_0, r_0)$.
\begin{proposition}\label{Prop:code-measure}
There exist positive constants $T_0, c_1, c_2>0$ such that the following hold:
\begin{itemize}
\item[(1)] For every flow $\{\F_t\}_{t \in \R}$-invariant Borel probability measure $m$ on $\Fc_\k$,
there exists a $\{\s_t\}_{t \in \R}$-invariant Borel probability measure $\tilde m$ on $\Sus(\SS_0, r_0)$ such that
\[
\Pi_{0\ast}\tilde m=m,
\]
and further
there exists a $\s$-invariant Borel probability measure $\lambda$ on $\SS_0$ such that
\[
c_1 m \le \pi_\ast \(w_\ast \lambda \otimes \Leb_{[0, T_0)}\) \le c_2 m,
\] 
where $\pi:\partial^2 \G \times \R \to \Fc_\k$ is the quotient map.
\item[(2)] For every $\G$-invariant Radon measure $\L$ on $\partial^2 \G$,
there exists a $\s$-invariant Borel probability measure $\lambda$ on $\SS_0$ such that
\[
c_1 \L\otimes dt \le \sum_{x \in \G}x.\(w_\ast\lambda\otimes \Leb_{[0, T_0)}\) \le c_2 \L\otimes dt \quad \text{on $\partial^2 \G \times \R$}.
\]
\end{itemize}
\end{proposition}

\proof
The claims have been shown for $\Sus(\SS, r_N)$ in \cite[Lemma 4.5 and Proposition 4.7]{TopFlows},
and the same proofs work for $\Sus(\SS_0, r_0)$ by Theorem \ref{Thm:enhanced}.
\qed

\subsection{Entropy}\label{Sec:entropy}

For a compact topological space $\Xc$, and a continuous map $f: \Xc \to \Xc$,
let $\Mcc(f, \Xc)$ be the space of all $f$-invariant Borel probability measures $\n$ on $\Xc$.
For any finite Borel partition $\Pc=\{A_1, \dots, A_k\}$ of $\Xc$,
let
\[
H(\n, \Pc):=-\sum_{i=1}^k \n(A_i)\log \n(A_i) \quad \text{for $\n \in \Mcc(f, \Xc)$},
\] 
where $0\log 0:=0$.
Letting
\[
h(f, \n, \Pc):=\lim_{n \to \infty}\frac{1}{n}H\Big(\n, \bigvee_{i=0}^{n-1}f^{-i}\Pc\Big)=\inf_{n \ge 1}\frac{1}{n}H\Big(\n, \bigvee_{i=0}^{n-1}f^{-i}\Pc\Big),
\]
where $\bigvee_{i=0}^{n-1}f^{-i}\Pc:=\Pc\vee f^{-1}\Pc\vee\cdots \vee f^{-(n-1)}\Pc$ denotes the partition consisting of
\[
A_{i_1}\cap f^{-1}A_{i_2}\cap \cdots \cap f^{-(n-1)}A_{i_n} \quad \text{for $i_1, \dots, i_n \in \{1, \dots, k\}$},
\]
and the limit exists since the function $n \mapsto H(\n, \bigvee_{i=0}^{n-1}f^{-i}\Pc)$ is subadditive,
we define the \textit{measure theoretical entropy} of $f$ for $(\Xc, \n)$ by
\[
h(f, \n):=\sup_\Pc h(f, \n, \Pc),
\]
where the supremum is taken over all finite Borel partitions $\Pc$.
In the case of a flow $\{\s_t\}_{t \in \R}$ on $\Xc$ with a flow invariant probability measure $\n$, 
the entropy is defined to be $h(\s_1, \n)$ for the time-one map $\s_1$. 

If $\L$ is a Radon measure on $\partial^2 \G$,
then we define
\[
D_q(\x_-, \x_+: \L):=\liminf_{r \to 0}\frac{\log \L(B(\x_-, r)\times B(\x_+, r))}{\log r} \quad \text{for $(\x_-, \x_+) \in \partial^2 \G$}, 
\]
where (open) balls $B(\x, r)$ are defined in terms of a quasi-metric $q$ in $\partial \G$.

\begin{lemma}\label{Lem:dim-ent}
Suppose that $\L$ is a Radon measure on $\partial^2 \G$ and $\lambda$ is a $\s$-invariant Borel probability measure on $\SS_0$ such that 
\[
w_\ast \lambda \le C \L \quad \text{on $\partial^2 \G$}
\]
for some positive constant $C>0$.
If for the strongly hyperbolic metric $\hat d$ defining $r_0$ in $\Sus(\SS_0, r_0)$,
\[
D_{\hat q}(\x_-, \x_+: \L) \ge D \quad \text{for $\L$-almost every $(\x_-, \x_+)$ in $\partial^2 \G$},
\]
where $\hat q$ is the associated quasi-metric in $\partial \G$ with $\hat d$,
then
we have that
\[
h(\s, \lambda)+\frac{D}{2}\int_{\SS_0}\P_B\,d\lambda \ge 0,
\]
where $\P_B(\o)=\hat \b_o(w(\o_0), w_\ast(\o))$ for $\o \in \SS_0$.
\end{lemma}

\proof
Let us identify $\Sus(\SS_0, r_0)$ with $Z=\{(\o, t) \in \SS_0 \times \R \ : \ 0 \le t < r_0(\o)\}$.
Considering the family of cylinder sets $[\o_0, \dots, \o_{N-1}]$ of length $N$ in $\SS_0$ (where $N$ is the fixed natural number appearing in the definition of $r_0$),
we define the partition $\Pc$ of $Z$ by the sets
\[
\{(\o, t) \in \SS_0 \times \R \ : \ \o \in [\o_0, \dots, \o_{N-1}], \ 0 \le t<r_0(\o)\}
\]
over the cylinder sets.
Let $\e_0:=(1/2)\min_{\o \in \SS_0}r_0(\o)$ and $\s_\ast:=\s_{\e_0}$ for the suspension flow $\{\s_t\}_{t \in \R}$ in $\Sus(\SS_0, r_0)$, and we analyze the partitions
\[
\s_\ast^n \Pc \vee \cdots \vee \s_\ast^{-n}\Pc \qquad \text{for $n \in \Z_{\ge 1}$}.
\]
For each $n \ge 1$,
let $\Pc_{[-n, n]}([\o, t])$ denote the set in this partition containing $[\o, t]$ in $\Sus(\SS_0, r_0)$ (or, in $Z$).
We fix a rough geodesic $\g_\o$ from $\R$ to the hyperbolic metric space $(\G, \hat d)$ 
with extreme points $w_\ast(\o)=(\x_-(\o), \x_+(\o))$ such that the distance between $o$ and $\g_\o$ is realized at $0$.
Note that for a large enough $R>0$, for every $[\o', t'] \in \Pc_{[-n, n]}([\o, t])$ we have
\[
w_\ast(\o') \in \Oc(\g_\o(t-\e_0 n), R) \times \Oc(\g_\o(t+\e_0 n), R).
\]
By the assumption, we have $w_\ast \lambda \otimes \Leb_{[0, T_0)}\le C \L\otimes dt$,
where $T_0:=\max_{\o \in \SS_0}r_0(\o)$.
Letting $\n_\lambda:=\frac{1}{\lambda\otimes \Leb(Z)}\lambda\otimes dt|_Z$,
we obtain
for $C':=C/\lambda\otimes \Leb(Z)$,
\begin{align*}
\n_\lambda(\Pc_{[-n, n]}([\o, t])) 
&\le \n_\lambda\(\{w_\ast(\o') \in \Oc(\g_\o(t-\e_0 n), R) \times \Oc(\g_\o(t+\e_0 n), R), \ |t-t'| \le T_0\}\)\\
&\le (2T_0 C') \cdot \L\(\Oc(\g_\o(t-\e_0 n), R) \times \Oc(\g_\o(t+\e_0 n), R)\).
\end{align*}
Hence taking the logarithms on both sides and dividing by $-2n$,
we have that 
\begin{align*}
\liminf_{n\to\infty} -\frac{1}{2n}\log \n_\lambda\(\Pc_{[-n, n]}([\o, t])\)
&\ge\liminf_{n \to \infty} -\frac{1}{2n}\L\(\Oc(\g_\o(t-\e_0 n), R) \times \Oc(\g_\o(t+\e_0 n), R)\)\\
&= \frac{\e_0}{2} D_{\hat q}(\x_-(\o), \x_+(\o): \L),
\end{align*}
where we have compared shadows with balls relative to $\hat q$ (see \eqref{Eq:shadow-ball}).
Furthermore, $h(\s_\ast, \n_\lambda, \Pc)$ is equal to:
\begin{align*}
\lim_{n \to \infty}\frac{1}{2n}H\big(\n_\lambda, \bigvee_{i=-n}^n \s_\ast^{-i}\Pc\big)
&=\liminf_{n \to \infty}-\frac{1}{2n}\int_{\Sus(\SS_0, r_0)}\log \n_\lambda\(\Pc_{[-n, n]}([\o, t])\)\,d\n_\lambda\\
&\ge \frac{\e_0}{2}\int_{\Sus(\SS_0, r_0)}D_{\hat q}(\x_-(\o), \x_+(\o): \L)\,d\n_\lambda,
\end{align*}
where the last inequality follows from the Fatou lemma.
Therefore if $D_{\hat q}(\x_-, \x_+: \L) \ge D$ for $\L$-almost every $(\x_-, \x_+)$ in $\partial^2 \G$,
then since $w_\ast \lambda \le C \L$ by the assumption and
$h(\s_\ast, \n_\lambda) \ge h(\s_\ast, \n_\lambda, \Pc)$ by the definition of the measure theoretical entropy,
we have
\[
h(\s_\ast, \n_\lambda) \ge \frac{\e_0}{2}D.
\]
The Abramov formula implies that
\[
h(\s_\ast, \n_\lambda)=\e_0 h(\s_1, \n_\lambda) \quad \text{and} \quad h(\s_1, \n_\lambda)=h(\s^N, \lambda)/\int_{\SS_0}r_0\,d\lambda,
\]
(e.g., \cite[Theorem 4.1.4 and Corollary 4.1.10]{FisherHasselblatt}).
Noting that $h(\s^N, \lambda)=N h(\s, \lambda)$ and $\int_{\SS_0}r_0\,d\lambda=N \int_{\SS_0}\tilde \k\,d\lambda$,
and further $\tilde \k=-\P_B+U\circ \s-U$ by Lemma \ref{Lem:coboundary},
we obtain
\[
h(\s_1, \n_\lambda) \ge \frac{D}{2} \quad \text{and} \quad h(\s, \lambda) \ge \frac{D}{2}\int_{\SS_0}\tilde \k\,d\lambda \ge \frac{D}{2}\(-\int_{\SS_0}\P_B\,d\lambda\),
\]
concluding the claim (cf.\ \cite[Lemma 6.3]{TopFlows}).
\qed

\begin{lemma}\label{Lem:localint_potential}
Let $\p$ be a $\G$-invariant tempered potential on $\G$,
and let $\m$ and $\check \m$ be finite Borel measures on $\partial \G$
such that there exists a $\G$-invariant Radon measure on $\partial^2 \G$ equivalent to $\check \m \otimes \m$.
If there exists a bounded measurable function $\P$ on $(\SS_0, \s)$ such that
\[
S_n\P(\o)=\sum_{i=0}^{n-1}\P(\s^i(\o))=-\p(o, \x_n(\o))+O(1) \quad \text{for all $\o \in \SS_0$},
\]
and a $\s$-invariant Borel probability measure $\lambda$ on $(\SS_0, \s)$ 
satisfying that
\[
w_\ast \lambda \le C \check \m \otimes \m,
\]
then for the local intersection number $\tau$ of $\p$ relative to $\hat d$, i.e., $\tau^\p(\x)=\tau$ for $\m$-almost every $\x \in \partial \G$, we have that
\[
\int_{\SS_0}\P\,d\lambda=\tau \int_{\SS_0}\P_B\,d\lambda,
\]
where $\P_B(\o)=\hat \b_o(w(\o_0), w_\ast(\o))$ for $\o \in \SS_0$.
\end{lemma}

\proof
Lemma \ref{Lem:localint} shows that $\tau^\p(\x)=\tau$ for $\m$-almost every $\x \in \partial \G$ for a constant $\tau$ relative to $\hat d$.
Since $S_n\P(\o)=-\p(o, \x_n(\o))+O(1)$ by the assumption, and $S_n\P_B(\o)=-\hat d(o, \x_n(\o))+O(1)$ for $\o \in \SS_0$ (see Example \ref{Ex:hatBpotential}), 
we have that by the definition of $\tau^\p(\x)$ relative to $\hat d$ and by the Birkhoff ergodic theorem,
\[
\lim_{n \to \infty}\frac{1}{n}S_n \P(\o)=\tau^\p(w_\ast(\o))\lim_{n \to \infty}\frac{1}{n}S_n\P_B(\o) \quad \text{for $\lambda$-almost every $\o \in \SS_0$ and in $L^1$}.
\]
Since by assumption $w_\ast \lambda \le C \check \m\otimes \m$,
and $\t^\p(\x)$ is constant $\m$-almost everywhere,
we have that $\tau^\p(w_\ast(\o))=\tau$ for $\lambda$-almost every $\o \in \SS_0$,
and thus conclude the claim.
\qed

\subsection{Real analyticity of Manhattan curves}\label{Sec:Manhattan}

Let us consider a pair of hyperbolic metrics $d$ and $d_\ast$ from $\Dc_\G$.
For $d$ (resp.\ $d_\ast$),
we define the stable translation length
\[
\ell[x]:=\lim_{n \to \infty}\frac{1}{n}d(o, x^n) \quad \text{for $x \in \G$},
\]
(resp.\ $\ell_\ast[x]$), where the limit exists since $d(o, x^n)$ is subadditive in $n$.
Since $\ell[x]=\ell[yxy^{-1}]$, it defines a function on the set of conjugacy classes $\conj_\G$ in $\G$.
Let
\[
\Qc(s, t):=\sum_{[x] \in \conj_\G}\exp(-s \ell_\ast[x]-t\ell[x]) \quad \text{for $s, t \in \R$},
\]
where $[x]$ is the conjugacy class of $x \in \G$.
The {\it Manhattan curve} $\Cc_M$ for a pair of metrics $(d, d_\ast)$ in $\Dc_\G$ is defined by the boundary of the convex domain where 
$\Qc(s, t)<\infty$ in $\R^2$.
If we define
\[
\Pc(s, t):=\sum_{x \in \G}\exp(-s d_\ast(o, x)-t d(o, x)),
\]
then
it has been shown that $\Qc(s, t)<\infty$ if and only if $\Pc(s, t)<\infty$ \cite[Proposition 3.1]{CT}.
For each fixed $s$, let $\th(s)$ be the abscissa of convergence in $t$ for $\Pc(s, t)$.
We call $\th$ the {\it pressure curve}
and identify the Manhattan curve $\Cc_M$ with the graph of the function $\th$.
It has been shown that $\Cc_M$ is $C^1$ for any pair of metrics from $\Dc_\G$ \cite[Theorem 1.1]{CT},
and $C^2$ for any pair of strongly hyperbolic metrics, or any pair of word metrics \cite[Theorems 1.3 and 4.14]{CT}.

\begin{theorem}\label{Thm:Manhattan}
Let $\G$ be a non-elementary hyperbolic group.
For any pair of strongly hyperbolic metrics $d$ and $d_\ast$ in $\Dc_\G$,
the Manhattan curve $\Cc_M$ is real analytic.
Moreover, $\th'(0)=-\t(d_\ast/d)$ where $\t(d_\ast/d)$ is the local intersection number of $d_\ast$ relative to $d$.
\end{theorem}

\proof
Let $\b_{\ast o}$ and $\b_o$ be the Busemann cocycles for $d_\ast$ and for $d$ respectively,
and let
\[
\P_\ast(\o):=\b_{\ast o}(w(\o_0), w_\ast(\o)) \quad \text{and} \quad \P(\o):=\b_o(w(\o_0), w_\ast(\o)) \quad \text{for $\o \in \SS_0$}.
\]
Since $d_\ast$ and $d$ are strongly hyperbolic metrics,
$\P_\ast$ and $\P$ define H\"older continuous potentials on $(\SS_0, \s)$.
For any $(s, t) \in \R^2$, let us define
\[
\P_{s, t}:=s \P_\ast+t \P.
\]
If $t=\th(s)$, then $s d_\ast +t d$ has exponent $0$ with respect to any hyperbolic metric in $\Dc_\G$, and thus by Lemma \ref{Lem:Pr}, we have $\Pr(\P_{s, t})=0$.
Therefore for each $(s, t) \in \Cc_M$, we have 
$\Pr_{\Cc_0}(\P_{s, t}) \le 0$.
If $t=\th(s)$, then there exists a $\G$-invariant Radon measure $\L_{s, t}$ on $\partial^2 \G$ equivalent to $\m_{s, t}\otimes \m_{s, t}$ (Example \ref{Ex:interpolation}).
For this measure $\L_{s, t}$ on $\partial^2 \G$, by Proposition \ref{Prop:code-measure}(2),
there exists a $\s$-invariant Borel probability measure $\lambda_{s, t}$ on $\SS_0$ such that
\[
w_\ast \lambda_{s, t}\otimes \Leb_{[0, T_0)} \le c\L_{s, t}\otimes dt \quad \text{on $\partial^2 \G \times \R$},
\]
for some positive constants $T_0, c>0$, and thus $w_\ast \lambda_{s, t} \le C \m_{s, t}\otimes \m_{s, t}$ on $\partial^2 \G$ for some positive constant $C>0$.
For the fixed strongly hyperbolic metric $\hat d$ which is used to define $\Sus(\SS_0, r_0)$ in Theorem \ref{Thm:enhanced} and the corresponding quasi-metric $\hat q$ in $\partial \G$,
we have that
\[
D_{\hat q}(\x_-, \x_+: \L_{s, t})=2\tau_{s, t} \quad \text{for $\L_{s, t}$-almost all $(\x_-, \x_+) \in \partial^2 \G$},
\]
where $\tau_{s, t}$ is the local intersection number of $s d_\ast+t d$ relative to $\hat d$ for $\m_{s, t}$-almost everywhere
(by Lemma \ref{Lem:localint}).
Moreover, Lemma \ref{Lem:localint_potential} implies that
\[
\int_{\SS_0}\P_{s, t}\,d\lambda_{s, t} = \tau_{s, t}\int_{\SS_0}\P_B\,d\lambda_{s, t},
\]
where $\P_B(\o)=\hat \b_o(w(\o_0), w_\ast(\o))$.
By Lemma \ref{Lem:dim-ent} we have that
\[
h(\s, \lambda_{s, t})+\int_{\SS_0}\P_{s, t}\,d\lambda_{s,t} =h(\s, \lambda_{s, t})+\tau_{s, t}\int_{\SS_0}\P_B\,d\lambda_{s, t} \ge 0.
\]
The variational principle (Proposition \ref{Prop:VP}) shows that $\Pr_{\Cc_0}(\P_{s, t}) \ge 0$.
Therefore for $s \in \R$, if $t=\th(s)$, then
\[
\Pr_{\Cc_0}(\P_{s, t})=0,
\]
and Proposition \ref{Prop:VP}
implies that $\lambda_{s, t}$ is the unique equilibrium state for $\P_{s, t}$.
Note that the function $\R^2 \to \R$, $(s, t) \mapsto \Pr_{\Cc_0}(\P_{s, t})$ is real analytic by Theorem \ref{Thm:TF}(2).
Moreover if $t=\th(s)$, then
\[
\frac{\partial}{\partial t}\Big|_{t=\th(s)}\Pr_{\Cc_0}(\P_{s, t})=\int_{\SS_0}\P\,d\lambda_{s,t}<0,
\]
where the last strict inequality follows since $S_n \P<0$ for all large enough $n$.
This holds for each $s \in \R$.
The implicit function theorem implies that $\th$ is real analytic on the whole of $\R$, concluding the first claim.
The second claim holds since
\[
\th'(0)=-\frac{\int_{\SS_0}\P_\ast\,d\lambda_{0,\th(0)}}{\int_{\SS_0}\P\,d\lambda_{0,\th(0)}},
\]
and the right hand side coincides with $-\t(d_\ast/d)$ by Lemmas \ref{Lem:localint} and \ref{Lem:localint_potential}.
\qed

In fact, we have a short proof that applies to Manhattan curves
for pairs of word metrics.
Since it is of interest on its own right, we include it below.

\begin{theorem}\label{Thm:Manhattan_word}
Let $\G$ be a non-elementary hyperbolic group.
For any pair of word metrics in $\G$,
the Manhattan curve $\Cc_M$ is real analytic. 
\end{theorem}

\proof
For any pair of finite symmetric sets of generators $S, S_\ast$,
there exist a strongly Markov automatic structure $\Ac=(\Gc, w, S)$ and a function $d\f_\ast$ on the edges such that
\[
|w_\ast(\o)|_{S_\ast}=\sum_{i=0}^{n-1}d\f_\ast(\o_i),
\]
for all paths $\o=(\o_0, \dots, \o_{n-1})$ with $n=|w_\ast(\o)|_S$ from the initial state in the underlying directed graph \cite[Lemma 3.8]{CalegariFujiwara2010}.
Let $\P_{S_\ast}(\o):=d\f_\ast(\o_0)$ for $\o \in \wbar \SS^+$.
Since $\P_{S_\ast}$ depends only on the first coordinate of $\o$, it defines a H\"older continuous potential on $(\wbar \SS^+, \s)$.

For $s, t \in \R$, we consider the series $\Pc(s, t)$ for the pair $d_S, d_{S_\ast}$.
Defining the potential $\P_S:=1$, we take a maximal component $\Cc$ for $\P_S$.
If we define $\G_\Cc$ as the set of group elements which are obtained as images of the words corresponding to finite paths in $\Cc$,
then there exists a finite set $B$ in $\G$ such that $B \G_\Cc B=\G$ (cf.\ the proof of Proposition \ref{Prop:component}; see \cite[the proof of Lemma 4.19]{CT}).
For any $s, t \in \R$, letting
\[
\Pc_{\Cc}(s, t):=\sum_{x \in \G_\Cc}\exp(-s |x|_{S_\ast}-t|x|_S),
\]
we have that
\[
\Pc_\Cc(s, t) \le \Pc(s, t)=\sum_{x \in B\G_\Cc B}\exp(-s |x|_{S_\ast}-t|x|_S) \le e^{c(s, t)}(\# B)^2\, \Pc_\Cc(s, t),
\]
where $c(s, t):=2|s|\max_{x \in B}|x|_{S_\ast}+2|t| \max_{x \in B} |x|_{S}$.
Therefore for each fixed $s \in \R$, the divergence exponents of $\Pc_\Cc(s, t)$ and $\Pc(s, t)$ coincide.
Furthermore this exponent given by $\Pr_\Cc(-s\P_\ast)$ since for $t>\th(s)$,
\[
\Pc_\Cc(s, t)=\sum_{n=1}^\infty\sum_{\o=(\o_0, \dots, \o_{n-1}) \ \text{in $\Cc$}} \exp(-s S_n \P_\ast(\o)-t n).
\]
By \eqref{Eq:Gibbs} in Proposition \ref{Prop:VP}, we have
$\th(s)=\Pr_\Cc(-s\P_\ast)$ for every $s \in \R$.
Since $\Pr_\Cc(-s\P_\ast)$ is real analytic in $s$ by Theorem \ref{Thm:TF}(2),
we conclude the claim.
\qed

\subsection{Uniqueness of measure of maximal Hausdorff dimension}\label{Sec:UMMD}

For a hyperbolic metric $d$ in $\Dc_\G$,
let $q$ be the corresponding quasi-metric in $\partial \G$.
Let
\[
q_\times((\x_1, \y_1), (\x_2, \y_2)):=\max\{q(\x_1, \x_2), q(\y_1, \y_2)\}
\]
for $(\x_i, \y_i) \in \partial^2 \G$ and $i=1, 2$.
For a set $E$ in $\partial \G$, we denote
the Hausdorff dimension of $E$ relative to $q$ by $\dim_H(E, q)$.
For a set $E$ in $\partial^2 \G$, we write $\dim_H(E, q_\times)$
for the Hausdorff dimension of $E$ relative to $q_\times$.
For a Borel measure $\m$ on $\partial \G$,
let us define the \textit{lower Hausdorff dimension} of $\m$ by
\[
\underline{\dim}_H(\m, q):=\inf\{\dim_H(E, q) \ : \ \text{$\m(E)>0$ and $E$ is Borel}\},
\]
and the \textit{upper Hausdorff dimension} of $\m$ by
\[
\overline{\dim}_H(\m, q):=\inf\{\dim_H(E, q) \ : \ \text{$\m(\partial \G \setminus E)=0$ and $E$ is Borel}\}.
\]
If $\underline{\dim}_H(\m, q)=\overline{\dim}_H(\m, q)$,
then the common value is called the \textit{Hausdorff dimension} of $\m$ relative to $q$ and is denoted by $\dim_H(\m, q)$.
For a Radon measure $\L$ on $(\partial^2 \G, q_\times)$,
we similarly define $\underline{\dim}_H(\L, q_\times)$, $\overline{\dim}_H(\L, q_\times)$ and $\dim_H(\L, q_\times)$ relative to $q_\times$.
Furthermore $\underline{\dim}_H(\L, q_\times)$ coincides with the essential infimum of $D_q(\x_-, \x_+: \L)$ with respect to $\L$, i.e., 
\[
\underline{\dim}_H(\L, q_\times)={\rm essinf}_\L\, D_q(\x_-, \x_+: \L),
\]
(cf.\ \cite[Section 2.2]{TopFlows} and the references therein). 

For $d \in \Dc_\G$ and a $\G$-invariant tempered potential $\p$ on $\G$,
if we have a $\G$-invariant Radon measure $\L$ on $\partial^2 \G$,
then
let $\t(\p/d: \L)$ be the essential supremum of $\tau^\p_{\inf}(\x_+)$ relative to $\L$ for $(\x_-, \x_+) \in \partial^2 \G$, i.e.,
\[
\tau(\p/d: \L):={\rm esssup}_\L\, \tau^\p_{\inf}(\x_+),
\]
(cf.\ Section \ref{Sec:localint}).
Let us define the {\it flip-involution} $\i(\x, \y):=(\y, \x)$ on $\partial^2 \G$.
Note that this is isometric in $q_\times$ and commutes with the $\G$-action, and thus for any $\G$-invariant Radon measure $\L$ on $\partial^2 \G$, if we define $\check \L:=\L\circ \i$,
then $\check \L$ is $\G$-invariant and Radon, furthermore,
\[
\t(\p/d:\check \L)={\rm esssup}_\L\, \tau^\p_{\inf}(\x_-).
\]
If we have a pair of $\G$-invariant tempered potentials $\p$ and $\check \p$ admitting a Gromov product,
then let $\L_\p$ be an associated $\G$-invariant Radon measure in the measure class of $\check \m \otimes \m$ where $\m$ and $\check \m$ are finite Borel measures satisfying \eqref{Eq:QC} for $\b_o^\p$ and $\b_o^{\check \p}$ respectively.
Let us fix the strongly hyperbolic metric $\hat d$ in Theorem \ref{Thm:enhanced}, and define $\hat q$ and $\hat q_\times$ the corresponding quasi-metrics in $\partial \G$ and $\partial^2 \G$ respectively.

\begin{theorem}\label{Thm:UMMD}
Let $\hat d$ be the strongly hyperbolic metric, and $\p, \check \p$ be a pair of $\G$-invariant tempered potentials admitting a Gromov product and H\"older continuous cocycles $\b_o^\p$, $\b_o^{\check\p}$ such that $\p$ has the exponent $0$ relative to $\hat d$.
For any $\G$-invariant Radon measure $\L$ on $\partial^2 \G$,
if 
\[
\underline{\dim}_H(\L, \hat q_\times) \ge \t(\check \p/\hat d:\check \L)+\t(\p/\hat d:\L),
\]
then $\L$ is a constant multiple of $\L_\p$. 
\end{theorem}

\proof
Fix $(\SS_0, \s)$ corresponding to a component $\Cc_0$ whose adjacency matrix has an eigenvalue of maximal modulus (see Section \ref{Sec:enhanced}).
Let 
\[
\P(\o):=\b_o^\p(w(\o_0), w_\ast(\o)) \quad \text{for $\o \in \SS_0$}.
\]
Note that $\P$ is H\"older continuous on $(\SS_0, \s)$ since by assumption $\b_o^\p$ is a H\"older continuous cocycle.
Moreover by assumption, $\p$ has the exponent $0$ relative to $\hat d$.
By Lemma \ref{Lem:Pr} this implies that $\Pr(\P)=0$, hence $\Pr_{\Cc_0}(\P)\le 0$.
For any $\G$-invariant Radon measure $\L$ on $\partial^2 \G$,
by Proposition \ref{Prop:code-measure}(2) there exist a $\s$-invariant Borel probability measure $\lambda$ on $\SS_0$ and a constant $T_0>0$ such that
$\sum_{x \in \G}x \cdot (w_\ast\lambda \otimes \Leb_{[0, T_0)})$
and $\L\otimes dt$ are dominated by each other up to multiplicative constants.
Hence there exists a constant $C>0$ such that $w_\ast \lambda \le C \L$ on $\partial^2 \G$.
Let
\[
\chi(\p/\hat d:\L):=\frac{1}{2}(\t(\check \p/\hat d:\check \L)+\t(\p/\hat d:\L)).
\]
Since by the assumption that $\underline{\dim}_H(\L, \hat q_\times) \ge 2\chi(\p/\hat d:\L)$,
we have
\[
D_{\hat q}(\x_-, \x_+: \L) \ge 2\chi(\p/\hat d: \L) \quad \text{for $\L$-almost every $(\x_-, \x_+)$ in $\partial^2 \G$},
\]
Lemma \ref{Lem:dim-ent} shows that
\begin{equation}\label{Eq:dim-ent}
h(\s, \lambda)+\chi(\p/\hat d: \L)\int_{\SS_0}\P_B\,d\lambda \ge 0.
\end{equation}
Noting that 
\[
S_n \P(\o)=-\p(o, \x_n(\o))+O(1) \quad \text{and} \quad S_n \P_B(\o)=-\hat d(o, \x_n(\o))+O(1)
\]
for $\o \in \SS_0$ and for integers $n \ge 0$,
we have by the Birkhoff ergodic theorem,
\[
\frac{S_n \P(\o)}{S_n \P_B(\o)} \to \t^\p(\x_+(\o)) \quad \text{for $\lambda$-almost every $\o \in \SS_0$},
\]
where $\t^\p(\eta)$ stands for the local intersection number of $\p$ relative to $\hat d$ at $\eta$ in $\partial \G$.
Hence we have 
\[
\int_{\SS_0}\P\,d\lambda\ge \t(\p/\hat d: \L) \int_{\SS_0}\P_B\,d\lambda
\]
since $S_n\P_B<0$ for all large enough $n$ and $\t^\p(\eta_+) \le \t(\p/\hat d: \L)$ for $\L$-almost every $(\eta_-, \eta_+)$ in $\partial^2 \G$.
Furthermore, we have that for $\o \in \SS_0$ and integers $n>1$,
\begin{align*}
S_{-n}\P(\o)
&=-\p(o, w_\ast(\o_{-(n-1)},\dots, \o_0))+O(1)\\
&=-\check \p(o, w(\o_0)^{-1}w(\o_{-1})^{-1}\cdots w(\o_{-(n-1)})^{-1})+O(1)\\
&=-\check \p(o, \x_{-n}(\o))+O(1),
\end{align*}
where we have used the $\G$-invariance in the second equality and \eqref{Eq:RG} in the last equality,
and similarly,
\[
S_{-n} \P_B(\o)=-\hat d(o, \x_{-n}(\o))+O(1).
\]
By the Birkhoff ergodic theorem,
we have
\[
\frac{S_{-n} \P(\o)}{S_{-n} \P_B(\o)} \to \t^{\check \p}(\x_-(\o)) \quad \text{for $\lambda$-almost every $\o \in \SS_0$},
\]
and further, by the same reasoning as above,
\[
\int_{\SS_0}\P\,d\lambda\ge \t(\check \p/\hat d: \check \L) \int_{\SS_0}\P_B\,d\lambda.
\]
Summarizing the above inequalities together with \eqref{Eq:dim-ent} and the variational principle (Proposition \ref{Prop:VP}) on $(\SS_0, \s)$ yields
\[
\Pr_{\Cc_0}(\P) \ge h(\s, \lambda)+\int_{\SS_0}\P\,d\lambda 
\ge h(\s, \lambda)+\frac{1}{2}(\t(\check \p/\hat d:\check \L)+\t(\p/\hat d:\L))\int_{\SS_0}\P_B\,d\lambda 
\ge 0.
\]
Therefore $\Pr_{\Cc_0}(\P)=0$, and furthermore $\lambda$ satisfies that
\[
\Pr_{\Cc_0}(\P)=h(\s, \lambda)+\int_{\SS_0}\P\,d\lambda=0.
\]
This implies that $\lambda$ is the unique equilibrium state $\lambda_{\Cc_0}$ for $\P$ on $(\SS_0, \s)$ by Proposition \ref{Prop:VP}.
Since $\Pr_{\Cc_0}(\P)=0$, by Lemma \ref{Lem:codingPSpotential}(1) we have $w_\ast \lambda_{\Cc_0} \le C \L_\p$ for some positive constant $C>0$,
implying that $\L\otimes dt$ is absolutely continuous with respect to $\L_\p \otimes dt$ on $\partial^2 \G \times \R$.
Hence $\L$ is absolutely continuous with respect to $\L_\p$ on $\partial^2 \G$,
and since $\L_\p$ is ergodic with respect to the $\G$-action on $\partial^2 \G$ by Theorem \ref{Thm:Hopf},
the $\G$-invariant Radon measure $\L$ is a constant multiple of $\L_\p$, as required.
\qed

If we further assume that $\L$ be a $\G$-invariant Radon measure which is ergodic,
then we have the following valid for any quasi-metric associated with $d \in \Dc_\G$.

\begin{theorem}\label{Thm:UMMD_ergodic-flip}
Let $d$ be a hyperbolic metric in $\Dc_\G$, and $\p, \check \p$ be a pair of $\G$-invariant tempered potentials admitting a Gromov product and H\"older continuous cocycles $\b_o^\p$, $\b_o^{\check \p}$ such that $\p$ has the exponent $0$ relative to $d$.
If $\L_\p$ is an associated $\G$-invariant Radon measure with $\p$, then
\[
\dim_H(\L_\p, q_\times)=\t(\check \p/d:\check \L_\p)+\t(\p/d:\L_\p).
\]
Furthermore,
for any $\G$-invariant Radon measure $\L$ which is ergodic on $\partial^2 \G$,
if
\[
\dim_H(\L, q_\times) \ge \t(\check \p/d:\check \L)+\t(\p/d:\L),
\]
then $\L$ is a constant multiple of $\L_\p$, where $q_\times$ is the quasi-metric associated with $d$.
\end{theorem}

\proof
Since $\p$ is a $\G$-invariant tempered potential and $\L$ is ergodic,
the local intersection number $\t^\p(\x_+)$ of $\p$ relative to $d$ exists at $\L$-almost every $\x_+$ as a genuine limit and is constant for $\L$-almost every $(\x_-, \x_+) \in \partial^2 \G$ by Lemma \ref{Lem:localint}.
Let us denote this constant by $\t(\p/d:\L)$.
Similarly, $\t^{\check \p}(\x_-)=\t(\check \p/d:\check \L)$ for $\L$-almost every $(\x_-, \x_+) \in \partial^2 \G$.
Let us consider 
\[
\chi^\p(\x_-, \x_+:d):=\liminf_{t \to \infty}\frac{\p(\g(-t), \g(t))}{d(\g(-t), \g(t))} \qquad \text{for $(\x_-, \x_+) \in \partial^2 \G$},
\]
where $\g$ is a rough geodesic in $(\G, d)$ such that $\g(t) \to \x_+$ and $\g(-t) \to \x_-$ as $t\to \infty$, respectively.
We note that the definition does not depend on the choice of $\g$ by the Morse lemma.
Since $\p$ satisfies \eqref{Eq:RG} relative to $d$,
\begin{align*}
\p(\g(-t), \g(t))
&=\p(\g(-t), \g(0))+\p(\g(0), \g(t))+O(1)\\
&=\check \p(\g(0), \g(-t))+\p(\g(0), \g(t))+O(1),
\end{align*}
we have that
\[
\chi^\p(\x_-, \x_+:d)=\lim_{t\to\infty}\frac{\p(\g(-t), \g(t))}{d(\g(-t), \g(t))}=\frac{1}{2}(\t(\check \p/d:\check \L)+\t(\p/d:\L)),
\]
for $\L$-almost every $(\x_-, \x_+) \in \partial^2 \G$.
Let 
\[
\chi(\p/d:\L):=\frac{1}{2}(\t(\check \p/d:\check \L)+\t(\p/d:\L)).
\]
By definition, we have that
\begin{equation}\label{Eq:D}
D_q(\x_-, \x_+:\L)=2\liminf_{t \to \infty}\frac{\log \L(\Oc(\g(-t), R)\times \Oc(\g(t), R))}{d(\g(-t), \g(t))},
\end{equation}
for a large enough $R>0$, comparing shadows and balls (see \eqref{Eq:shadow-ball}).
By the definition of $\L_\p$, we have by \eqref{Eq:D},
\[
D_q(\x_-, \x_+:\L_\p)=2\chi^\p(\x_-,\x_+:d)=\t(\check \p/d:\check \L_\p)+\t(\p/d:\L_\p),
\] 
for $\L_\p$-almost every $(\x_-, \x_+) \in \partial^2 \G$,
and the Frostman-type lemma (see \cite[Section 2.2]{TopFlows}) shows the Hausdorff dimension formula for $\L_\p$.

Applying $\p$ to the strongly hyperbolic metric $\hat d$,
we have $\chi^{\hat d}(\x_-, \x_+:d)$ obtained as the limit for $\L$-almost every $(\x_-, \x_+) \in \partial^2 \G$,
and 
\[
\chi^{\hat d}(\x_-, \x_+:d)=\frac{1}{2}(\t(\hat d/d:\check \L)+\t(\hat d/d:\L)) \quad \text{for $\L$-almost every $(\x_-, \x_+) \in \partial^2 \G$},
\]
whence by \eqref{Eq:D},
\[
D_q(\x_-, \x_+:\L)=D_{\hat q}(\x_-, \x_+:\L)\cdot\chi^{\hat d}(\x_-, \x_+:d) \quad \text{for $\L$-almost every $(\x_-, \x_+) \in \partial^2 \G$}.
\]
Moreover, we have that
\[
\chi^\p(\x_-, \x_+:d)=\chi^\p(\x_-, \x_+:\hat d)\cdot \chi^{\hat d}(\x_-, \x_+:d)\quad \text{for $\L$-almost every $(\x_-, \x_+) \in \partial^2 \G$}.
\]

Therefore, if $\dim_H(\L, q_\times) \ge 2\chi(\p/d:\L)$,
then $D_q(\x_-, \x_+:\L) \ge 2\chi(\p/d:\L)$ for $\L$-almost every $(\x_-, \x_+) \in \partial^2 \G$,
and thus $D_{\hat q}(\x_-, \x_+:\L) \ge 2\chi(\p/\hat d:\L)$ for $\L$-almost every $(\x_-, \x_+)\in \partial^2 \G$,
where we have used $\t(\hat d/d:\L), \t(\hat d/d:\check \L)>0$ since $\hat d$ and $d$ are quasi-isometric.
Moreover, if $\p$ has the exponent $0$ relative to $d$, then it has the exponent $0$ relative to $\hat d$ (Remark \ref{Rem:exponent0}).
The rest follows from Theorem \ref{Thm:UMMD}.
\qed

\section{Random walks and their harmonic measures}\label{Sec:RW}

For a probability measure $p$ on $\G$ such that the support is finite and generates $\G$ as a semigroup,
let us consider the Green function 
\[
G(x, y):=\sum_{n=0}^\infty p^{\ast n}(x^{-1}y) \quad \text{for $x, y \in \G$},
\]
where $p^{\ast n}$ denotes the $n$-fold convolution of $p$.
If we define
\[
\p(x, y):=-\log G(x, y),
\]
then $\p$ is a $\G$-invariant tempered potential relative to any metric $d$ in $\Dc_\G$.
Indeed, $\p$ satisfies \eqref{Eq:RG} by the Ancona inequality: $G(x, z) \le C G(x, y)G(y, z)$ whenever $x$, $y$ and $z$ are aligned on a geodesic segment (in a word metric) where $C$ is a uniform constant.
Note that the inequality $G(x, z) \ge C'G(x, y)G(y, z)$ holds for all $x$, $y$ and $z$ where $C'=G(\id, \id)^{-1}$ by the definition of the Green function.
These inequalities 
imply that
the Green function is multiplicative up to uniform multiplicative constants along geodesic rays in a word metric (see \cite{GouezelAncona} for $p$ with an infinite support).
Furthermore under the present assumption, \eqref{Eq:QE} holds with
\[
\b_o^\p(x,\x)=\p(x,o)-2(x|\x)_o^\p \quad \text{and} \quad \b_o^\p(x,\x)=-\log K_o(x, \x) \quad \text{for $(x, \x) \in \G \times \partial \G$},
\]
where $K_o(x, \x)$ is the Martin kernel based at $o$.
It has been shown that $K_o(x, \x)$ is H\"older continuous in $\x$ (see \cite[Theorem 3.3]{INO}) relative to some (equivalently, any) quasi-metric $q$ for each $x$.
Note that $G(x, y)$ is not necessarily symmetric in $x$ and $y$ unless $p$ defines a symmetric random walk on $\G$.
See \cite[Corollaries 2.4 and 3.4]{INO} and \cite[Section 3]{BHM11}.

If we define $\check G(x, y):=G(y, x)$, then $\check G(x, y)$ is the Green function associated with the measure $\check p$,
where $\check p(x):=p(x^{-1})$ for $x \in \G$.
The harmonic measure $\n$ defined by $p$ (the $p$-stationary measure) on $\partial \G$ satisfies that
\[
\frac{d x_\ast \n}{d\n}(\x)=K_o(x, \x) \quad \text{for $(x, \x) \in \G \times \partial \G$},
\]
and similarly, for the harmonic measure $\check \n$ defined by $\check p$ with the associated Martin kernel $\check K_o(x, \x)$.
The {\it Naim kernel} $\phi(\y, \x)$ is defined by the extension of the following function on $\partial^2 \G$:
\[
\frac{G(x, y)}{G(x, o)G(o, y)} \quad \text{for $(x, y) \in \G\times \G$},
\]
and $\phi(\x, \y)$ is H\"older continuous on each compact subset in $\partial^2 \G$
(see e.g., \cite[the proof of Proposition 5.7]{THaus} and \cite[Proposition 3.17]{LedrappierFreeGroups}).
Moreover, we have that
\[
-\log \phi(\x, \y)+\log \phi(x^{-1}\x, x^{-1}\y)=-\log \check K_o(x, \x)-\log K_o(x, \y) \quad \text{for $(x, \x, \y)\in \G \times \partial^2 \G$}
\]
(for the derivation, see \cite[ibid.]{THaus, LedrappierFreeGroups}).
Therefore the pair of cocycles $-\log \check K_o(x, \x)$ and $-\log K_o(x, \y)$ admits
a Gromov product $[\x|\y]_o:=\log \phi(\x, \y)$, and
we obtain a $\G$-invariant Radon measure 
$\exp([\x|\y]_o)\check \n \otimes \n$ on $\partial^2 \G$ by Proposition \ref{Prop:Radon}.

For $d \in \Dc_\G$,
let us consider the following series
\[
\Pc_{G, d}(s, t):=\sum_{x \in \G}G(o, x)^s e^{-t d(o, x)} \quad \text{for $s, t \in \R$}.
\]
For each $s \in \R$, the divergence exponent in $t$ is finite, and if $s=1$, then it is $0$ (cf.\ \cite[Theorem 1.1(ii)]{BHM11} and \cite[Lemma 3.2]{THaus}; the proof can be adapted to apply to any metric in the class $\Dc_\G$).
Let us consider any strongly hyperbolic metric $\hat d$ in $\Dc_\G$, and let $t=\th_{G, \hat d}(s)$ be the divergence exponent of $\Pc_{G, \hat d}(s, t)$ in $t$ for each fixed $s \in \R$.

\begin{theorem}\label{Thm:Manhattan_harm}
The function $s \mapsto \th_{G, \hat d}(s)$ is real analytic on the whole $\R$.
\end{theorem}

\proof
The proof proceeds as in Theorem \ref{Thm:Manhattan}.
Let us define 
\[
\Psi_G(\o):=-\log K_o(w(\o_0), w_\ast(\o)) \quad \text{and} \quad \Psi_B(\o):=\hat \b_o(w(\o_0), w_\ast(\o)) \quad \text{for $\o \in \SS_0$},
\]
where $K_o(x, \x)$ is the Martin cocycle and $\hat \b_o(x, \x)$ is the Busemann cocycle associated with $\hat d$.
Define 
\[
\Psi_{s,t}:=s \Psi_G +t \Psi_B \quad \text{for $s, t \in \R$}.
\]
Note that $\Psi_G$ and $\Psi_B$ are H\"older continuous potentials on $(\SS_0, \s)$.
The rest follows as in Theorem \ref{Thm:Manhattan}; we omit the details.
\qed

\begin{remark}\label{Rem:multifractal}
The Legendre transform of the function $\th_{G, \hat d}(s)$ in Theorem \ref{Thm:Manhattan_harm} provides the {\it multifractal profile} for the harmonic measure $\n$ with respect to the quasi-metric $\hat q$ associated with the strongly hyperbolic metric $\hat d$ (see \cite[Theorem 1.2]{THaus} for the case of word metrics and \cite[Theorem 3.8]{CT} for a related result).
The result shows that the multifractal profile is a real analytic function.
Moreover, one can show that the domain is a bounded closed interval and the profile function extends continuously at extremes by adapting the method in \cite[Proposition 4.20]{CT}.
In particular this applies to finite range random walks on surface groups. In certain special cases, it has been recently shown that the harmonic measure is singular with respect to the Lebesgue measure on the boundary of the hyperbolic plane \cite[Theorem 1.1]{Kosenko} and \cite[Theorem 1]{KosenkoTiozzo}.
In those cases, $\th_{G, \hat d}$ is strictly convex.
\end{remark}

\section{Dominated representations}\label{Sec:main-rep}

\subsection{Dominated representations: preliminaries}\label{Sec:rep}

Let $\R^m$ be the Euclidean space of dimension $m \ge 2$ equipped with the standard inner product.
For each linear transformation $g$ of $\R^m$,
let us denote the singular values $\s_1(g) \ge \cdots \ge \s_m(g)$ repeated with multiplicity (where $\s_i(g)$ are the eigenvalues of $\sqrt{g^\ast g}$ and $g^\ast$ denotes the transpose of $g$).
Given $1 \le p < m$ we say that $g$ has a {\it gap of index} $p$ if $\s_p(g)>\s_{p+1}(g)$.

Let us consider any invertible matrix $g$ with a gap of index $p$ for some $1 \le p < m$.
Let $U_p(g)$ be the space spanned by eigenvectors of $\sqrt{gg^\ast}$ with the $p$-largest singular values,
\[
U_p(g):=\mathrm{span}\left\{w_i \ : \ \sqrt{gg^\ast}w_i=\s_i(g)w_i, \ i=1, \dots, p\right\}.
\]
If $v_i$ denotes an eigenvector of $\sqrt{g^\ast g}$ for the eigenvalue $\s_i(g)$,
then the space $U_p(g)$ is spanned by $\{g v_i\}_{i=1, \dots ,p}$.
Since $\s_{m-i+1}(g^{-1})=\s_i(g)^{-1}$ for $i=1, \dots, m$,
for any $g$ with a gap of index $p$, the inverse $g^{-1}$ has a gap of index $m-p$.
Let 
\[
S_{m-p}(g):=U_{m-p}(g^{-1}).
\] 
By definition, $S_{m-p}(g)$ consists of eigenvectors of $\sqrt{g^\ast g}$ with eigenvalues $\s_{i}(g)$ for $i=p+1, \dots, m$, and $S_{m-p}(g)$ and $g^{-1}U_p(g)$ are orthogonal.

Letting $\|\cdot\|$ be the associated norm in $\R^m$,
we define 
\[
\|g\|:=\sup\left\{\frac{\|g v\|}{\|v\|} \ : \ 0\neq v \in \R^m\right\} \quad \text{for $g \in \GL(m, \R)$}.
\]
Let us define a metric $d_\Pb$ in the projective space $\Pb(\R^m)$ (the space of one-dimensional subspace in $\R^m$) by 
\[
d_\Pb(x_1, x_2):=\inf\{\|v_1-v_2\| \ : \ v_i \in x_i, \ \|v_i\|=1, \ i=1, 2\}.
\]
For two closed sets $X_1, X_2$ in $\Pb(\R^m)$, we define
\[
d_\Pb(X_1, X_2):=\inf\{d_\Pb(x_1, x_2) \ : \ x_i \in X_i, \ i=1, 2\}.
\]
Below it is useful to recall the following inequalities: For $g, h \in \GL(m, \R)$, if $g$ and $gh$ have a gap of index $p$,
then
\begin{equation}\label{Eq:AAB}
d_\Pb(U_p(g), U_p(gh)) \le \|h\|\|h^{-1}\|\frac{\s_{p+1}(g)}{\s_p(g)},
\end{equation}
and
\begin{equation}\label{Eq:BABA}
d_\Pb(U_p(gh), g U_p(h)) \le \|g\|\|g^{-1}\|\frac{\s_{p+1}(h)}{\s_p(h)},
\end{equation}
(see \cite[Lemmas A4 and A5]{BochiPotrieSambarino}).

\subsection{Dominated representations}\label{Sec:dom_rep}
Let $\G$ be a non-elementary hyperbolic group.
For an integer $m \ge 2$, we consider a linear representation of $\G$ on $\R^m$ given by a homomorphism
\[
\rho: \G \to \GL(m, \R).
\] 
Let $\lambda(g)$ be the spectral radius of the linear transformation $g$, i.e., the modulus of the largest eigenvalue of $g$ in absolute value.
Recall that $g$ is {\it proximal} if $g$ has an eigenvalue of modulus $\lambda(g)$ with multiplicity one.
If $g$ is proximal, then its associated map on the projective space $\Pb(\R^m)$ has an attracting fixed point  given by the (one-dimensional) eigenspace and a repelling set determined by the complementary codimension one subspace in $\R^m$.
For $\e>0$, we say that $g$ is {\it $\e$-proximal} if $g$ is proximal,
moreover, for the attracting fixed point $x_{g+}$ and the repelling set $X_{g-}$,
defining
\[
b^\e:=\{x \in \Pb(\R^m) \ : \ d_\Pb(x, x_{g+}) \le \e\} \quad \text{and} \quad B^\e:=\{x \in \Pb(\R^m) \ : \ d_\Pb(x, X_{g-}) \ge \e\},
\]
we have $d_\Pb(x_{g+}, X_{g-}) \ge 2\e$, 
$g(B^\e) \subset b^\e$ and $g$ restricted on $B^\e$ is $\e$-Lipschitz \cite[Section 6]{Benoist}.

For each element $g \in \GL(m, \R)$,
let $\s_1(g)$ (resp.\ $\s_2(g)$) be the largest (resp.\ the second largest) singular value of $g$,
and $\lambda_1(g)$ (resp.\ $\lambda_2(g)$) be the modulus of the eigenvalue of $g$ with the largest (resp.\ the second largest) modulus. We have $\s_1(g)=\|g\|$ and $\lambda_1(g)=\lambda(g)$.
Note that $\lambda_1(g)$ and $\lambda_2(g)$ depend only on the conjugacy class of $g$.
We say that $\rho:\G \to \GL(m, \R)$ is a {\it dominated representation} if
there exist positive constants $c, C>0$ such that
\begin{equation}\label{Eq:dom}
\frac{\s_2(\rho(x))}{\s_1(\rho(x))} \le C e^{-c |x|} \quad \text{for all $x \in \G$},
\end{equation}
where $|x|$ stands for a word norm of $x$ with a fixed set of generators $S$.
In fact, such a representation is called a {\it $1$-dominated representation} in \cite{BochiPotrieSambarino};
we discuss only this case for the sake of brevity.
The definition does not depend on the choice of $S$, whence we omit the symbol from the notation.
Note that the condition \eqref{Eq:dom} implies that
\[
\frac{\s_{m}(\rho(x))}{\s_{m-1}(\rho(x))} \le C e^{-c|x|} \quad \text{for all $x \in \G$},
\]
since $\s_i(g^{-1})=\s_{m-i+1}(g)^{-1}$ for $g \in \GL(m, \R)$.

For any dominated representation $\rho: \G \to \GL(m, \R)$,
there exist positive constants $c_1, c_2>0$ such that
\begin{equation}\label{Eq:dom_QI}
c_1^{-1}|x|-c_2 \le \log \|\rho(x)\| \le c_1 |x| \quad \text{for $x \in \G$}.
\end{equation}
Indeed, for each $\rho$ the condition \eqref{Eq:dom} holds for the representation $x \mapsto |\det \rho(x)|^{-1/m}\rho(x)$, and applying it to this new representation we have
$\log \s_2(\rho(x))\ge -(m-1)\log \s_1(\rho(x))$;
combining with \eqref{Eq:dom} implies
\[
\log \s_1(\rho(x)) \ge -(m-1)\log \s_1(\rho(x))+c|x|-\log C.
\]
Rearranging the terms gives the first inequality in \eqref{Eq:dom_QI} since $\s_1(\rho(x))=\|\rho(x)\|$.
The second inequality of \eqref{Eq:dom_QI} holds with $c_1:=\max\{\|\rho(s)\| \ : \ s \in S\}$.

Note that taking the $n$-th roots of both sides in \eqref{Eq:dom} for $x^n \in \Gamma$ and letting $n$ tend to infinity, we have
\begin{equation}\label{Eq:dom_sp}
\frac{\lambda_2(\rho(x))}{\lambda_1(\rho(x))} \le e^{-c\ell[x]} \quad \text{for all $[x] \in \conj_\G$},
\end{equation}
where $\ell[x]$ denotes the stable length function defined by the chosen word metric of a conjugacy class $[x]$.
We observe that by \eqref{Eq:dom_sp} if $\rho(x)$ is proximal and $\ell[x]$ is large enough,
then $\rho(x)$ has a strong contraction property away from the repelling set in $\Pb(\R^m)$.
This is used to show that some proximal element whose attracting fixed point and repelling set are far enough is $\e$-proximal for a given $\e>0$.
Appealing to \eqref{Eq:dom_QI}, we obtain for a positive constant $c_1>0$,
\begin{equation}\label{Eq:dom_QI_sp}
c_1^{-1}\ell[x] \le \log \lambda(\rho(x)) \le c_1 \ell[x] \quad \text{for all $[x] \in \conj_\G$}.
\end{equation}

Associated with a dominated representation $\rho$, 
there exists a pair of $\rho$-equivariant continuous maps from the boundary $\partial \G$ to the projective space of $\R^m$ and that of the dual ${\R^m}^\ast$,
\[
E_\rho:\partial \G \to \Pb(\R^m) \quad \text{such that $\rho(x)E_\rho(\x)=E_\rho(x \x)$ for $x \in \G$ and $\x \in \partial \G$},
\]
and
\[
\check E_\rho: \partial \G \to \Pb({\R^m}^\ast) \quad \text{such that $\rho(x)\check E_\rho(\x)=\check E_\rho(x \x)$ for $x \in \G$ and $\x \in \partial \G$}
\]
where we understand that for $\theta \in \Pb({\R^m}^\ast)$  the action is the one induced by $\rho(x)\theta=\theta\circ \rho(x)^{-1}$.
Let us identify each $\theta \in \Pb({\R^m}^\ast)$ with the space of hyperplanes in $\R^m$ determined by the kernel of a representative of $\theta$ in the dual ${\R^m}^\ast$, where the kernel is independent of the choice.
 Through this identification $\check E_\rho(\x)$ and $E_\rho(\y)$ determine transversal subspaces in $\R^m$ if $\x \neq \y$.
Note that for every infinite order element $x$ in $\G$, the image $\rho(x)$ is proximal; for the attracting (resp.\ repelling) fixed point $x_+$ (resp.\ $x_-$) of $x$ in the boundary, $E_\rho(x_+)$ (resp.\ $\check E_\rho(x_-)$) yields an attracting fixed point (resp.\ a repelling invariant set) in $\Pb(\R^m)$.

More concretely, the maps $E_\rho$ and $\check E_\rho$ are obtained as follows.
For each $\x$ in $\partial \G$,
letting $\{x_n\}_{n=0}^\infty$ be a geodesic converging to $\x$,
we define
\[
E_\rho(\x):=\lim_{n \to \infty}U_1(\rho(x_n))
\qquad \text{and} \qquad
\check E_\rho(\x):=\lim_{n \to \infty}U_{m-1}(\rho(x_n)),
\]
where the limits are defined in $\Pb(\R^m)$ and $\Pb({\R^m}^\ast)$ with the fixed metrics, 
respectively;
they exist by \eqref{Eq:AAB}, see \cite[Lemma 4.7]{BochiPotrieSambarino}.
Moreover, the limits do not depend on the choice of geodesics. 
By \eqref{Eq:AAB},
we have that there exist positive constants $C, c>0$ such that
\[
d_\Pb(E_\rho(\x), E_\rho(\y)) \le C e^{-c (\x|\y)_o} \quad \text{for $\x, \y \in \partial \G$},
\]
and an analogous inequality holds for $\check E_\rho$, i.e.,
$E_\rho$ (resp.\ $\check E_\rho$) defines a H\"older continuous map from $\partial \G$ to $\Pb(\R^m)$ (resp.\ $\Pb({\R^m}^\ast)$).
The maps $E_\rho$ and $\check E_\rho$ are $\rho$-equivariant by \eqref{Eq:BABA}.
It holds that for any $\x, \y \in \partial \G$ with $\x \neq \y$,
the two subspaces $\check E_\rho(\x)$ and $E_\rho(\y)$ are transversal;
this is due to \cite[Theorem A]{BochiGourmelon} (see also \cite[Appendix A.4]{BochiPotrieSambarino}).

For each $(u, v) \in \Pb({\R^m}^\ast)\times \Pb(\R^m)$, 
let
\[
\llangle u, v\rrangle:=\frac{|u_0(v_0)|}{\|u_0\|\|v_0\|},
\]
where $u_0$ and $v_0$ represent $u$ and $v$ respectively, and the right hand side does not depend on the choices.
Note that $(u, v) \mapsto \llangle u, v\rrangle$ is continuous.
Identifying $u$ with the corresponding set in $\Pb(\R^m)$ determined by the hyperplane,
we have that
\begin{equation}\label{Eq:angle}
\llangle u, v\rrangle \le d_\Pb(u, v) \le 2 \llangle u, v\rrangle.
\end{equation}
This follows by noting that $\llangle u, v\rrangle$ is the sine of angle between $v$ and the corresponding hyperplane for $u$.

\begin{lemma}\label{Lem:rep_RG}
Fix a word metric in $\G$.
For any dominated representation $\rho:\G \to \GL(m, \R)$, and for any constants $C, R \ge 0$,
there exists a constant $c_0>0$ such that
for any $C$-rough geodesic segment between $x$ and $y$, and any point $z$ which is in the $R$-neighborhood of the $C$-rough geodesic segment,
we have that
\[
\|\rho(x^{-1}y)\| \ge c_0\|\rho(x^{-1}z)\|\|\rho(z^{-1}y)\|.
\]
\end{lemma}
\proof
Translating from the left hand side by $z$,
we consider the case when the $C$-rough geodesic passes through $o$ within distance $R$.
Fix a large enough $R>0$, and let
\[
c_R:=\inf\{\llangle \check E_\rho(\x), E_\rho(\y)\rrangle \ : \ (\x|\y)_o \le R+4 \d\}>0.
\]
For $\y \in \Oc(y, R)$, we have that
\[
d_{\Pb}(U_1(\rho(y)), E_\rho(\y)) \le C e^{-c|y|},
\]
and
an analogous estimate holds for $\check E(\x)$ and $U_{m-1}(\rho(x))$ with $\x \in \Oc(x, R)$ (cf.\ \eqref{Eq:AAB} and \eqref{Eq:dom}).
Hence there exist a constant $C_\ast$ depending on $c, C$ and $c_R$ such that for $x, y \in \G$ satisfying that a geodesic between $x$ and $y$ passes through $o$
and $|x|, |y| \ge C_\ast$,
\[
\llangle U_{m-1}(\rho(x)), U_1(\rho(y))\rrangle \ge \frac{c_R}{2}.
\]
Letting $\a$ be the angle between $U_{m-1}(\rho(x))$ and $U_1(\rho(y))$,
we have $\sin \a \ge c_R/2$.
Since $U_{m-1}(\rho(x))=S_{m-1}(\rho(x)^{-1})$,
which is perpendicular to $\rho(x)U_1(\rho(x)^{-1})$,
for any unit vector $w \in \rho(y)^{-1}U_1(\rho(y))$,
\[
\|\rho(x)^{-1}\rho(y)w\| \ge \s_1(\rho(x)^{-1})\|\rho(y)w\|\sin \a,
\]
whence $\|\rho(x)^{-1}\rho(y)\| \ge \|\rho(x)^{-1}\|\|\rho(y)\|c_R/2$.
Therefore, fixing a small enough constant $c_0>0$,
for all $x, y$ in $\G$ and all $z$ in the $R$-neighborhood of any $C$-rough geodesic between $x$ and $y$,
we obtain the claim.
\qed

For $(g, \x) \in \G \times \partial \G$,
if we define
\[
c_\rho(g, \x):=\log \frac{\|\rho(g)^{-1}v_\x\|}{\|v_\x\|} \quad \text{for $v_\x \in E_\rho(\x)\setminus\{0\}$},
\]
then it satisfies the following cocycle identity:
\[
c_\rho(gh, \x)=c_\rho(h, g^{-1}\x)+c_\rho(g, \x) \quad \text{for $g, h \in \G$ and $\x \in \partial \G$}.
\]
Since we have
\[
d_\Pb(E_\rho(\x), E_\rho(\y)) \le C e^{-c (\x|\y)_o} \quad \text{for $\x, \y \in \partial \G$},
\]
(cf.\ \eqref{Eq:AAB} and \eqref{Eq:dom})
and $\|\rho(g)^{-1}v_\x-\rho(g)^{-1}v_{\y}\| \le \|\rho(g)^{-1}\| \|v_\x-v_{\y}\|$,
there exists a constant $C_g$ depending on $g$ such that
\[
\left|c_\rho(g, \x)-c_\rho(g, \y)\right| \le C_g e^{-c (\x|\y)_o}  \quad \text{for $\x, \y \in \partial \G$}.
\]
This implies that $c_\rho(g, \x)$ is H\"older continuous in $\x$ for each fixed $g$.

For a unit vector $v_x$ in $\rho(x)^{-1}U_1(\rho(x))$,
we have $\|\rho(x)v_x\|=\|\rho(x)\|$, and thus
\begin{equation}\label{Eq:limit_cocycle}
\log\frac{\|\rho(g)^{-1}\rho(x)v_x\|}{\|\rho(x)v_x\|} \to c_\rho(g, \x),
\end{equation}
if $x$ tends to $\x$ along some (equivalently, any) geodesic ray.
Indeed, for any geodesic ray $\{x_n\}_{n=0}^\infty$ converging to $\x$,
we have that $\rho(x_n)v_{x_n} \in U_1(\rho(x_n))$ and $U_1(\rho(x_n))$ tends to $E_\rho(\x)$ as $n$ goes to infinity.

\begin{lemma}\label{Lem:rep_QE}
For every $\x \in \partial \G$,
if a sequence $\{x_n\}_{n=0}^\infty$ in $\G$ converges to $\x$,
then for each $g \in \G$,
\[
c_\rho(g, \x)=\lim_{n \to \infty}(\log \|\rho(g^{-1}x_n)\|-\log \|\rho(x_n)\|).
\]
\end{lemma}

\proof
For $g, x \in \G$,
and for a unit vector $v_x \in \rho(x)^{-1}U_1(\rho(x))$,
we have that
\[
\|\rho(g^{-1}x)v_x\|\ge \|\rho(g^{-1}x)\|\sin \a, 
\]
where $\sin \a=\llangle S_{m-1}(\rho(g^{-1}x)), \rho(x)^{-1}U_1(\rho(x))\rrangle$.
Note that 
\[
\sin (\pi/2-\a)=\llangle U_1(\rho(g^{-1}x)^\ast), U_1(\rho(x)^\ast)\rrangle,
\]
and 
\[
d_\Pb(U_1(\rho(g^{-1}x)^\ast), U_1(\rho(x)^\ast))\le \|\rho(g^{-1})\|\|\rho(g)\|\frac{\s_2(\rho(x)^\ast)}{\s_1(\rho(x)^\ast)} \le C_g e^{-c|x|},
\]
by \eqref{Eq:AAB} and \eqref{Eq:dom} since $\s_i(\rho(x)^\ast)=\s_i(\rho(x))$ for all $i$.
Hence $\sin \(\frac{\pi}{2}-\a\)$ tends to $0$ by \eqref{Eq:angle} and thus $\sin \a$ tends to $1$ exponentially fast as $|x| \to \infty$, respectively.
Noting that
\[
\frac{\|\rho(g^{-1}x)\|}{\|\rho(x)\|} \ge \frac{\|\rho(g^{-1}x) v_x\|}{\|\rho(x) v_x\|} \ge \frac{\|\rho(g^{-1}x)\|}{\|\rho(x)\|}\sin \a,
\]
where $\|\rho(x)v_x\|=\|\rho(x)\|$,
we conclude the claim.
\qed

Lemmas \ref{Lem:rep_RG} and \ref{Lem:rep_QE} show that for each dominated representation $\rho: \G \to \GL(m, \R)$, if we define 
\[
\psi(x, y):=\log \|\rho(x^{-1}y)\| \quad \text{for $x, y \in \G$},
\]
then $\psi$ satisfies \eqref{Eq:QE} and \eqref{Eq:RG} relative to each fixed $d \in \Dc_\G$ by the Morse lemma (see Section \ref{Sec:PStempered}).
Moreover, since $\psi$ is $\G$-invariant,
the function $\psi:\G \times \G \to \R$ defines a $\G$-invariant tempered potential relative to a metric $d$ in $\Dc_\G$ (Definition \ref{Def:tempered}).

Let us define the contragredient representation $\check \rho:\G \to \GL(m,\R)$, i.e., $\check \rho(x)$ is the transpose of $\rho(x^{-1})$.
Note that
\[
\check \psi(x, y)=\log \|\rho(y^{-1}x)\|=\log \|\check \rho(x^{-1}y)\| \quad \text{for $x, y \in \G$},
\]
and the pair of cocycles $c_\rho(x, \x)$ and $c_{\check \rho}(x, \x)$ admits a Gromov product
\[
[\x | \y]_{\rho}:=-\log \llangle \check E_\rho(\x), E_\rho(\y)\rrangle \quad \text{for $(\x, \y) \in \partial^2 \G$}.
\]
Indeed, we have that
\begin{align*}
[x \x|x \y]_\rho-[\x|\y]_\rho
&=-\log \frac{|u_{x \x}(w_{x \y})|}{\|u_{x \x}\|\|w_{x \y}\|}+\log \frac{|u_{\x}(w_\y)|}{\|u_\x\|\|w_\y\|}\\
&=\log \frac{\|\rho(x)u_\x\|}{\|u_\x\|}+\log \frac{\|\rho(x)w_\y\|}{\|w_\y\|}
=c_{\check \rho}(x^{-1}, \x)+c_\rho(x^{-1}, \y),
\end{align*}
where $u_\x$ and $w_\y$ represent $\check E_\rho(\x)$ and $E_\rho(\y)$ respectively,
and $\rho(x)u_\x=u_{x \x}$ and $\rho(x)w_\y=w_{x \y}$,
further we have used $u_{x \x}(w_{x \y})=u_\x(w_\y)$ in the second equality.
Note that $[\x|\y]_\rho$ defines a locally H\"older continuous function on $\partial^2 \G$.
Given results through Sections \ref{Sec:pre}-\ref{Sec:TF},
we obtain the following: 
Let us consider the series
\[
\sum_{x \in \G}\exp(-s\log\|\rho(x)\|) \quad \text{for $s \in \R$},
\]
and let $v_\rho$ be the divergence exponent in $s$.
Note that $v_\rho$ is positive and finite by \eqref{Eq:dom_QI}.
We call $v_\rho$ the {\it exponential growth rate} for $\rho$.
Furthermore, we have $v_{\check \rho}=v_\rho$.
Applying the above ideas to the cocycles $v_\rho c_\rho(x, \x)$ and $v_\rho c_{\check \rho}(x, \x)$,
we obtain finite measures $\m_\rho$ and $\m_{\check \rho}$ on $\partial \G$ 
satisfying \eqref{Eq:QC} by Proposition \ref{Prop:PSpotential},
where the exponent is $0$ for any metric $d$ in $\Dc_\G$.
Moreover, 
there exists a $\G$-invariant Radon measure
\[
\L_\rho:=\exp(v_\rho [\x|\y]_\rho)\m_{\check \rho}\otimes \m_{\rho} \quad \text{on $\partial^2 \G$},
\]
by Proposition \ref{Prop:Radon}.
Note that $\L_{\rho}$ is ergodic with respect to the $\G$-action on $\partial^2 \G$ by Theorem \ref{Thm:Hopf}.

Given the above preparations, Theorem \ref{Thm:UMMD_ergodic-flip} implies that
\[
\dim_H(\L_\rho, q_\times)=v_\rho(\t(\check \rho/d:\L_{\check \rho})+\t(\rho/d:\L_\rho)),
\]
for $d \in \Dc_\G$ where $q_\times$ is the quasi-metric associated with $d$ since $\check \L_\rho=\L_{\check \rho}$, and 
for any $\G$-invariant Radon measure $\L$ which is ergodic on $\partial^2 \G$,
if 
\[
\dim_H(\L, q_\times)\ge \t(\check \rho/d : \check \L)+\t(\rho/d : \L),
\]
then $\L$ is a constant multiple of $\L_\rho$.
This concludes 
Theorem \ref{Thm:UMMD_ergodic-flip_intro}.

\subsection{Intersections for pairs of dominated representations}

Fix an integer $m \ge 2$.
Let $\rho$ and $\rho_\ast$ be a pair of dominated representations from $\G$ to $\GL(m, \R)$.
We define the {\it Manhattan curve} for any such pair $(\rho, \rho_\ast)$
as the boundary of the domain where 
\[
\Qc_{\rho, \rho_\ast}(s, t):=\sum_{[x] \in \conj_\G}\exp(-s \log \lambda(\rho_\ast(x))-t \log \lambda(\rho(x))) \quad \text{for $(s, t) \in \R^2$},
\]
is finite.
Note that the Manhattan curve is convex by the H\"older inequality.
Let 
\[
\Pc_{\rho, \rho_\ast}(s, t):=\sum_{x \in \G}\exp(-s \log \|\rho_\ast(x)\|-t \log \|\rho(x)\|) \quad \text{for $(s, t) \in \R^2$}.
\] 
For each $s \in \R$, let $\th(s)$ be the divergence exponent of $\Pc_{\rho, \rho_\ast}(s, t)$ in $t$, and let us call $\th$ the {\it pressure curve} for the pair $(\rho ,\rho_\ast)$.
Analogous to the Manhattan curve for a pair of metrics \cite[Proposition 3.1]{CT},
we show that $\Qc_{\rho, \rho_\ast}(s, t)<\infty$ if and only if $\Pc_{\rho, \rho_\ast}(s, t)<\infty$ in Proposition \ref{Prop:Knieper} below.

\begin{lemma}\label{Lem:growth}
For large enough $L, R>0$ and
for each $s \in \R$, there exists a positive constant $C_{s, L, R}$ such that
\begin{equation*}\label{Eq:growth}
C_{s, L, R}^{-1}\,e^{\th(s)n} \le \sum_{x: \ |\log \|\rho(x)\|-n| \le L}e^{-s \log \|\rho_\ast(x)\|} \le C_{s, L, R}\,e^{\th(s)n},
\end{equation*}
for all positive integers $n$.
\end{lemma}
\proof
For each $s \in \R$, let $t:=\th(s)$.
Applying the cocycle $s c_{\rho_\ast}(x, \x)+tc_{\rho}(x, \x)$ to Proposition \ref{Prop:PSpotential},
we have a finite Borel measure $\m_{s, t}$ on $\partial \G$ satisfying \eqref{Eq:QC} 
(where the exponent is $0$ relative to any metric in $\Dc_\G$) and
\[
C^{-1} \le \exp(s \log \|\rho_\ast(x)\|+t\log \|\rho(x)\|)\cdot \m_{s, t}(\Oc(x, R)) \le C,
\]
for all $x \in \G$ and for some positive constants $C, R>0$.

Fix a large enough $L>0$ and for all positive integer $n$,
the shadows $\Oc(x, R)$ for $x$ with $|\log \|\rho(x)\|-n|\le L$ cover the boundary $\partial \G$ with a uniformly bounded multiplicity $M$;
indeed, for each $\x \in \partial \G$, there exists a $C$-rough geodesic $\g$ from $o$ converging to $\x$, and \eqref{Eq:dom_QI} shows that for a large enough $L>0$, there exists $x$ on $\g$ such that $|\log \|\rho(x)\|-n| \le L$ and $\x \in \Oc(x, R)$.
Furthermore, if $\x \in \Oc(x_1, R)\cap \Oc(x_2, R)$ and $|\log \|\rho(x_i)\|-n|\le L$ for $i=1, 2$, then $\g$ passes through the balls $B(x_i, R)$ at some $y_i$ for $i=1, 2$, respectively.
Suppose that $|y_1| \le |y_2|$, otherwise we change the role $y_1$ and $y_2$ below.
By \eqref{Eq:RG} shown in Lemma \ref{Lem:rep_RG},
we have 
\[
\log \|\rho(y_1)\|+\log \|\rho(y_1^{-1}y_2)\|=\log \|\rho(y_2)\|+O(1),
\]
and thus $\log \|\rho(y_1^{-1}y_2)\|=O(1)$, implying $|y_1^{-1}y_2|\le R_0$ for a constant $R_0$ by \eqref{Eq:dom_QI}.
Hence the multiplicity $M$ is at most $\#B(o, 2R+R_0)$.
This shows that for all $n$,
\begin{align*}
C^{-1}\m_{s, t}(\partial \G) &\le C^{-1}\sum_{x: \ |\log \|\rho(x)\|-n| \le L}\m_{s, t}(\Oc(x, R))
\le e^{|t|L}\sum_{x: \ |\log \|\rho(x)\|-n| \le L}\exp(-s \log \|\rho_\ast(x)\|-tn)\\
&\le C e^{2|t|L}\sum_{x: \ |\log \|\rho(x)\|-n| \le L}\m_{s, t}(\Oc(x, R))\le C e^{2|t|L} M\m_{s, t}(\partial \G),
\end{align*}
implying the claim.
\qed

\begin{proposition}\label{Prop:Knieper}
Let $\G$ be a non-elementary hyperbolic group.
If $\rho$ and $\rho_\ast$ are dominated representations from $\G$ to $\GL(m, \R)$ for $m \ge 2$,
then for each fixed $s \in \R$, and for a positive constant $L>0$,
letting
\[
\Cc_{\rho, \rho_\ast, s}(n, L):=\sum_{[z]: \  |\log \lambda(\rho(z))-n| \le L}e^{-s \log \lambda(\rho_\ast(z))},
\]
we have 
\[
C^{-1}n^{-2}e^{\th(s)n}\le\Cc_{\rho, \rho_\ast, s}(n, L)
\le
C e^{\th(s)n},
\]
for a constant $C>0$ and all large enough $n$.
In particular, for each $s \in \R$,
the divergence exponents in $t$ of $\Qc_{\rho, \rho_\ast}(s, t)$ and $\Pc_{\rho, \rho_\ast}(s, t)$ coincide.
\end{proposition}

\begin{remark}
We may take a pair of dominated representations into different spaces, namely,
$\rho: \G \to \GL(m_1, \R)$ and $\rho_\ast: \G \to \GL(m_2, \R)$ where $m_1$ and $m_2$ are integers at least $2$ and may not be identical.
The same result holds for this generality; see the proof of Proposition \ref{Prop:Knieper}.
However, we present the proof only when $m_1=m_2$ for the sake of brevity.
\end{remark}

If $\rho:\G \to \GL(m, \R)$ is a dominated representation,
then by the above discussion it yields a continuous equivariant map
\[
\partial^2 \G \to \Pb({\R^m}^\ast)\times \Pb(\R^m)\setminus \D, \quad (\x, \y) \mapsto (\check E_\rho(\x), E_\rho(\y)),
\] 
where $\D:=\{(u, v) \in \Pb({\R^m}^\ast)\times \Pb(\R^m) \ : \ \llangle u, v\rrangle=0\}$.
For any infinite order element $x \in \G$, let $x_+$ (resp.\ $x_-$) be the attracting  (resp.\ repelling) fixed point in $\partial \G$.

\begin{lemma}[cf.\ Lemma 5.7 in \cite{Sambarino}]\label{Lem:Benoist}
Assume that $\G$ is torsion-free. 
For every $0<\e<1/8$,
there exist a constant $c_\e\ge0$ and a finite set $F_\e$ in $\G$ such that
if $x \in \G \setminus F_\e$ and $\llangle \check E_\rho(x_-), E_\rho(x_+)\rrangle \ge 8 \e$,
then $\rho(x)$ is $\e$-proximal and 
\[
\log \lambda(\rho(x)) \le \log \|\rho(x)\|  \le \log \lambda(\rho(x))+c_\e.
\]
If $\G$ is not torsion-free, then the same statement holds with $F_\e$ replaced by the union of a finite set and the set of finite order elements.
\end{lemma}

\proof
First we assume that $\G$ is torsion-free so that every element has infinite order.
For any $0<\e<1/8$, let
\[
P_{8\e}:=\{(u, v) \in \Pb({\R^m}^\ast)\times \Pb(\R^m)\setminus \D \ : \ \llangle u, v\rrangle \ge 8\e\}.
\]
Note that $P_{8\e}$ is compact, and thus there exists a finite collection of open sets $U \times V$ in $\Pb({\R^m}^\ast) \times \Pb(\R^m)$ such that $\llangle u, v\rrangle \ge 7\e$ for all $(u, v) \in U \times V$, and any two points in $U$ (resp.\ $V$) are $\e/2$-close.
For each such $U \times V$,
consider the sets $\check E_\rho^{-1}(U)$ and $E_\rho^{-1}(V)$ in the boundary,
cover each of them by finitely many shadows $\Oc(w_1, R)$ and $\Oc(w_2, R)$, respectively for some $w_1, w_2 \in \G$ and some $R>0$, so that if $\x \in \Oc(w_1, R)$ and $\y \in \Oc(w_2, R)$, then $\llangle \check E_\rho(\x), E_\rho(\y)\rrangle \ge 6\e$.
Let
\[
U(w_1, w_2, R):=\{x \in \G \ : \ x_-\in \Oc(w_1, R), \ x_+ \in \Oc(w_2, R)\}.
\]
The $\d$-hyperbolicity implies that $(x_-|x_+)_o=(w_1|w_2)_o+O_{\d, R}(1)$ for all $x \in U(w_1, w_2, R)$,
and thus
\[
\ell[x]=|x|-2(x_-|x_+)_o=|x|+O_{\d, R, w_1, w_2}(1),
\]
(cf.\ \cite[Lemma 3.2]{CT}).
By \eqref{Eq:dom_sp}, we have that $\lambda_2(\rho(x))/\lambda_1(\rho(x)) \le e^{-c\ell[x]}$.
Since 
\[
\llangle \check E_\rho(x_-), E_\rho(x_+)\rrangle \ge 6\e
\]
 and $\check E_\rho(x_-)$ and $E_\rho(x_+)$ are $6 \e$-apart by \eqref{Eq:angle}, there exists a constant $C_\ast(\e)$ depending on $\e, w_1, w_2, \d$ and $R$ such that
if $|x| \ge C_\ast(\e)$,
then $\lambda_2(\rho(x))/\lambda_1(\rho(x))$ is uniformly small and $\rho(x)$ is $\e$-Lipschitz outside the $\e$-neighborhood of the hyperplane $\check E_\rho(x_-)$ in $\Pb(\R^m)$ and it sends that outside into the $\e$-neighborhood of $E_\rho(x_+)$.
Hence $\rho(x)$ is $2\e$-proximal by appealing to \cite[Lemma 6.2]{Benoist}.
Furthermore, since the set of $2\e$-proximal elements with norm $1$ is compact in the space of (nonzero) linear transformation \cite[Corollaire 6.3]{Benoist},
there exists a constant $c_{\e, \d, R, w_1, w_2}$ such that for all $x \in U(w_1, w_2, R)$ and $|x| \ge C_\ast(\e)$,
\[
\log \lambda(\rho(x)) \le \log \|\rho(x)\| \le \log \lambda(\rho(x))+c_{\e, \d, R, w_1, w_2}.
\]
Let $c_\e$ and $C_{\ast\ast}(\e)$ be the maxima of such $c_{\e, \d, R, w_1, w_2}$ and $C_\ast(\e)$ over finitely many possibilities of $w_1, w_2$, respectively.
Letting $F_\e$ be the ball of radius $C_{\ast\ast}(\e)$ in the word metric, we conclude the claim.

Finally, we consider the case when $\G$ is not necessarily torsion-free.
Applying the above argument to the set of infinite order elements,
we obtain the second claim.
\qed

\medskip

For a large enough $R>0$ and for any $z, w \in \G$,
let
\[
O(z, w, R):=\{x \in \G \cup \partial \G \ : \ (z|x)_w \le R\},
\]
and for $x \in \G$,
\[
U(x^{-1}, x, R):=\{z \in \G \ : \ z^{-1}\in O(o, x^{-1}, R), \ z \in O(o, x, R)\}.
\]
Note that for an $R'>R$ depending on $x$, for any infinite order element $z$ in $U(x^{-1}, x, R)$, 
we have $z_- \in O(o, x^{-1}, R')$ and $z_+ \in O(o, x, R')$.

\begin{lemma}\label{Lem:conj}
For each infinite order element $x$ in $\G$,
for all large enough $L>0$ there exists a constant $C_{x, R, L, \d}$ such that
for any positive integer $k$ and
any $z \in U(x^{-1}, x, R)$ with $|z|=k$,
the set
\[
C_k(z; x, L, R):=\{g\langle z\rangle \in \G/\langle z\rangle \ : \ g z g^{-1} \in U(x^{-1}, x, R), \ ||z|-k| \le L\}
\]
has cardinality at most $C_{x, R, L, \d}k$,
where $\langle z\rangle$ stands for the subgroup generated by $z$ in $\G$ and $\d$ is the hyperbolicity constant for a fixed word metric in $\G$.
\end{lemma}

\proof
The proof is found in the last paragraph in the proof of \cite[Proposition 3.1]{CT};
we include it here for the sake of convenience.
 We write $[o, z]$ for a geodesic segment between $o$ and $z$ (with respect to a fixed word metric) and take $g\in\G$ with $g\langle z\rangle \in C_k(z; x, L, R)$. Since by assumption
 \[
  (z^{-1}|z)_o=(x^{-1}|x)_o+O_R \quad \text{and} \quad (gz^{-1}g^{-1}|gzg^{-1})_o=(x^{-1}|x)_o+O_R,
 \]
we have that $\g(z):=\bigcup_{n \in \Z} z^n[o, z]$
 is the image of a $(c_0, c_1)$-quasi-geodesic invariant under $z$ for some $c_0, c_1>0$
 and
 similarly $\g(gzg^{-1})$ is the image of a $(c_0, c_1)$-quasi-geodesic invariant under $g z g^{-1}$.
 Note that $g \g(z)$ is also a $(c_0, c_1)$-quasi-geodesic invariant under $g z g^{-1}$,
 and thus $g \g(z)$ and $\g(g z g^{-1})$ lie within a bounded Hausdorff distance.
 Hence if $g\langle z\rangle \in C_k(z; x, L, R)$,
 then $g\g(z)$ passes through near $o$ within a bounded distance $C_{x, R, \d}$ of $o$.
 Here it is crucial that $c_0$, $c_1$ depend only on $x$, $R$ and the hyperbolicity constant.
 Since for all such $g$, the inverse $g^{-1}$ lies in the $C_{x, R, \d}$-neighborhood of $[o, z]$ up to translation by $z$,
 and $\g(z)$ is $z$-invariant, we obtain $\# C_k(z; x, L, R) \le C_{x, R, L, \d} k$,
 where we have incorporated the thickness of the annulus $L$ in the counting.
\qed

\medskip

The following lemma is due to \cite[the proof of Lemma 5.2]{CoornaertKnieper}.

\begin{lemma}\label{Lem:CoornaertKnieper}
Fix a large enough $R>0$.
There exists a pair of infinite order elements $x, y$ in $\G$ such that
the following hold:
If we define
\[
U:=O(o, x^{-1}, R), \quad V:=O(o, x, R),
\]
$U':=yU$, $V':=yV$,  $U'':=y^{-1}U$,  $V'':=y^{-1}V$, and
\[
S(U, V):=\{z \in \G \ : \ U\cap z V =\emptyset\},
\]
and define $S(U', V')$, $S(U'', V'')$ analogously,
then $U, V, U', V', U''$ and $V''$ are disjoint, 
\begin{equation}\label{Eq:x3}
x^3S(U, V)x^3 \subset U(x^{-1}, x, R), 
\end{equation}
and $S(U, V)=y^{-1}S(U', V')y=y S(U'',V'')y^{-1}$.
Furthermore
there exists a finite (possibly empty) set $F_{x, y, R}$ such that
\begin{equation}\label{Eq:SUV}
\G \setminus F_{x, y, R} \subset S(U, V) \cup S(U', V') \cup S(U'', V'').
\end{equation}
\end{lemma}

\proof

Since $\G$ is a non-elementary hyperbolic group,
there exists a pair of infinite order elements $x$ and $y$ such that
$n \mapsto x^n$ and $n \mapsto y^n$ for $n \in \Z$ give quasi-geodesics and their extremes are distinct in the boundary \cite[Section 3.1]{GromovHyperbolic}.
Replacing $x$ by its larger power if necessary,
we define
\[
\wt V:=O(x^{-3}, x^{-2}, R), \quad 
U:=O(o, x^{-1}, R), \quad V:=O(o, x, R) \quad \text{and} \quad \wt U:=O(x^3, x^2, R),
\]
satisfying that
\[
U\cap V =\emptyset, \quad (\G \cup \partial \G)\setminus U \subset \wt V \quad \text{and} \quad (\G \cup \partial \G)\setminus V \subset \wt U.
\]
Furthermore, we define
\[
U':=yU, \quad V':=yV, \quad U'':=y^{-1}U \quad \text{and} \quad V'':=y^{-1}V,
\]
in such a way that they and together with $U$ and $V$ are all disjoint,
replacing $y$ by its large power if necessary.
Note that for every $z \in S(U, V)$,
\[
(x^3 z x^3)\wt V \subset V \quad \text{and} \quad (x^3 z x^3)^{-1}\wt U \subset U,
\]
since $(x^3zx^3)\wt V \subset (x^3z) V \subset x^3((\G \cup \partial \G) \setminus U) \subset x^3 \wt  V  = V$ and similarly for the second expression.
Therefore,
we have that
\begin{equation*}
x^3 S(U, V) x^3 \subset U(x^{-1}, x, R),
\end{equation*}
since $o \in \wt V$ and $o \in \wt U$.
Moreover, since $U \times U' \times U''$ and $V \times V' \times V''$ are disjoint relatively compact sets in the space of ordered distinct triples in $\G \cup \partial \G$, denoted by $(\G \cup \partial \G)^{(3)}$.
The diagonal action of $\G$ on $(\G \cup \partial \G)^{(3)}$ is properly discontinuous \cite[Lemma 1.2 and Proposition 1.12]{BowditchConvergence}.
Hence the set of $z$ in $\G$ such that $U \times U' \times U'' \cap z(V \times V' \times V'')\neq \emptyset$
is finite, implying the claim.
\qed

\proof[Proof of Proposition \ref{Prop:Knieper}]
For each $s \in \R$ and a large enough $L>0$,
let
\[
\Pc_{\rho, \rho_\ast, s}(n, L):=\sum_{x:\ |\log\|\rho(x)\|-n|\le L}\exp\(-s\log\|\rho_\ast(x)\|\),
\]
and
\[
\Cc_{\rho, \rho_\ast, s}(n, L):=\sum_{[x] :\ |\log \lambda(\rho(x))-n|\le L}\exp\(-s \log\lambda(\rho_\ast(x))\).
\]
First we show that for every $s \in \R$, 
there exist positive constants $C, c>0$ such that
\begin{equation}\label{Eq:Knieper1}
\Cc_{\rho, \rho_\ast, s}(n, L)\le C \Pc_{\rho, \rho_\ast, s}(n, L+c),
\end{equation}
for all large enough $n$.
Observe that the left hand side depends only on the set of infinite elements in $\G$.
We apply Lemma \ref{Lem:Benoist} to both $\rho$ and $\rho_\ast$ simultaneously;
for any $0<\e<1/8$,
there exist a positive constant $c_\e>0$ and a set $\wt F_\e$ consisting of finite order elements together with some finite set in $\G$ such that
for all $x \in \G \setminus \wt F_\e$, if $\llangle \check E_\rho(x_-), E_\rho(x_+)\rrangle \ge 8 \e$, then we have
\[
\log \lambda(\rho(x)) \le \log \|\rho(x)\| \le \log \lambda(\rho(x))+c_\e,
\]
and analogous inequalities hold for $\rho_\ast(x)$.
Let 
\[
c_0:=\inf\{\llangle \check E_\rho(\x), E_\rho(\y)\rrangle \ : \ (\x|\y)_o \le 2\d\}>0,
\]
where $\d$ is a hyperbolic constant of $\Cay(\G, S)$ and
$c_0$ is positive since $(\x, \y)$ are uniformly away from the diagonal.
We fix $0<\e<c_0/8$.
Note that for any $[z] \in \conj_\G$ of an infinite order element,
there exists a representative $z$ such that $(z_-|z_+)_o \le 2\d$.
Therefore, for all large enough $n$, we have that
\[
\sum_{[z] :\ |\log \lambda(\rho(z))-n|\le L}\exp\(-s \log\lambda(\rho_\ast(z))\) \le e^{|s|c_\e}\sum_{z:\ |\log\|\rho(z)\|-n|\le L+c_\e}\exp\(-s\log\|\rho_\ast(z)\|\),
\]
implying \eqref{Eq:Knieper1}.

Second let us show that for every $s \in \R$,
there exist positive constants $C, c>0$ such that
\begin{equation}\label{Eq:Knieper2}
\Pc_{\rho, \rho_\ast, s}(n, L) \le C n^2\,\Cc_{\rho, \rho_\ast, s}(n, L+c),
\end{equation}
for all large enough $n$.
Observe that $\Pc_{\rho, \rho_\ast, s}(n, L)$ changes only by a multiplicative constant independent of $n$ under restricting the sum on a finite index subgroup of $\G$.
By the Selberg lemma applied to $\rho(\G)$, passing to a finite index subgroup in $\G$ if necessary,
we assume that the image $\rho(\G)$ consists of infinite order elements.

Let us take $x, y$ and define $U, V, U', V', U''$ and $V''$ as in Lemma \ref{Lem:CoornaertKnieper}.
Applying Lemma \ref{Lem:Benoist} to both $\rho$ and $\rho_\ast$ simultaneously, we have that
for each $0<\e<1/8$, 
there exist a positive constant $c_\e>0$ and a finite set $F_\e$ in $\G$ such that
if $z \in \G \setminus F_\e$ and $\llangle \check E_\rho(z_-), E_\rho(z_+)\rrangle \ge 8 \e$, 
then
\[
\log \lambda(\rho(z)) \le \log \|\rho(z)\| \le \log \lambda(\rho(z))+c_{\e}.
\] 
and similarly for $\rho_\ast(z)$.
Letting 
\[
r_{x,R'}:=\inf\{\llangle \check E_\rho(\x), E_\rho(\y) \rrangle \ : \ (\x, \y) \in O(o, x^{-1}, R') \times O(o, x, R') \}>0,
\]
for an $R'>R$ depending only on $x$,
we fix $0<\e<r_{x, R'}/8$.
Note that
if $z \in U(x^{-1}, x, R)$,
then $\llangle \check E_\rho(z_-), E_\rho(z_+)\rrangle \ge 8 \e$.

Given $z \in U(x^{-1}, x, R)$ with $|z|=k$,
the set
\[
C_k(z; x, L, R):=\{g\langle z\rangle \in \G/\langle z\rangle \ : \ g z g^{-1} \in U(x^{-1}, x, R), \ ||z|-k| \le L\}
\]
has cardinality at most $C_{x, R, L, \d}k$ by Lemma \ref{Lem:conj}.
Since $c_1^{-1}|z|-c_2 \le \log \|\rho(z)\| \le c_1|z|$ (see \eqref{Eq:dom_QI}) with constants $c_1, c_2$ depending only on the representation $\rho$ and the chosen word metric,
for $z \in U(x^{-1}, x, R)$ with $|\log \|\rho(z)\|-n|\le L$, 
one has
\[
\sum_{c_1^{-1}(n-L)\le k \le c_1(n+L+c_2)}\#C_k(z; x, L, R) \le C_{x, L, R, \d}\cdot n^2.
\]
Therefore, we have for all large enough $n$, 
\begin{align*}
&\sum_{z \in U(x^{-1}, x, R), \ |\log \|\rho(z)\|-n| \le L} \exp\(-s \log \|\rho_\ast(z)\|\) \\
&\qquad \qquad \qquad \qquad \qquad \qquad \le C_{x, L, R, \d}\cdot n^2\sum_{[z]: \ |\log\lambda(\rho(z))-n|\le L+c_{\e, x, R}}\exp\(-s \log \lambda(\rho_\ast(z))+|s|c_{\e}\),
\end{align*}
where $n$ is large enough so that $z \notin F_{x, y, R}\cup F_\e$ by the inequality $\log \|\rho(z)\| \le c_1 |z|$ in \eqref{Eq:dom_QI}.
By \eqref{Eq:x3} and \eqref{Eq:SUV} in Lemma \ref{Lem:CoornaertKnieper}, if we define 
\[
\wt L:=L+6\log \|\rho(x)\|+\log \|\rho(y)\|+\log \|\rho(y^{-1})\|+c_{\e, x, R},
\]
in analogous inequalities for $S(U', V')$ and $S(U'', V'')$ and combine them,
then for all large enough $n$,
\[
\sum_{z: \ |\log\|\rho(z)\|-n|\le L}\exp\(-s \log \|\rho_\ast(z)\|\)\le C_{\e, x, y, L, R, \d}\cdot n^2\sum_{[z]: \  |\log \lambda(\rho(z))-n| \le \wt L}\exp\(-s \log \lambda(\rho_\ast(z))\),
\]
implying \eqref{Eq:Knieper2}, as required.
The claim follows \eqref{Eq:Knieper1} and \eqref{Eq:Knieper2} combined with Lemma \ref{Lem:growth}.
\qed

\begin{theorem}\label{Thm:Manhattan_rep}
Fix an integer $m \ge 2$.
For any pair of dominated representations $\rho$ and $\rho_\ast$ from a non-elementary hyperbolic group $\G$ to $\GL(m, \R)$,
the Manhattan curve for the pair $(\rho, \rho_\ast)$ is real analytic.
\end{theorem}

\proof
Let $c_{\rho_\ast}(x, \x)$ and $c_{\rho}(x, \x)$ be the cocycles associated to $\rho_\ast$ and $\rho$ respectively.
If we define 
\[
\Psi_\ast(\o):=c_{\rho_\ast}(w(\o_0), w_\ast(\o)) \quad \text{and} \quad \Psi(\o):=c_{\rho}(w(\o_0), w_\ast(\o)) \quad \text{for $\o \in \SS_0$},
\]
then $\Psi_\ast$ and $\Psi$ define H\"older continuous potentials on the shift space $(\SS_0, \s)$ for an enhanced coding in Theorem \ref{Thm:enhanced}.
For each $s \in \R$, let $t:=\th(s)$ the value at $s$ of the pressure curve for the pair $(\rho, \rho_\ast)$.
Letting $\p_\ast(x, y):=\log \|\rho_\ast(x^{-1}y)\|$ and $\p(x, y):=\log \|\rho(x^{-1}y)\|$,
we define $\t_{s, t}$ as the local intersection number of $s \p_\ast+t\p$ relative to $\hat d$.
The rest follows as in the proof of Theorem \ref{Thm:Manhattan}.
\qed

\begin{proposition}\label{Prop:LDP}
Fix an integer $m \ge 2$.
Let $(\rho, \rho_\ast)$ be a pair of dominated representations from a non-elementary hyperbolic group $\G$ to $\GL(m, \R)$.
If we fix a large enough constant $L>0$ and define
\[
\Sc_{\rho, L}(n):=\{x\in \G:\ |\log \|\rho(x)\|-n|\le L\},
\]
then the large deviation principle holds for the sequence $(1/n)\log \|\rho_\ast(x)\|$
where $x$ is distributed according to the uniform counting measure on $\Sc_{\rho, L}(n)$ for $n=1, 2, \dots$.
More precisely, for any open set $U$ and any closed set $V$ in $\R$ such that $U \subset V$,
we have that
\begin{align*}
-\inf_{t \in U}I(t) 
&\le 
\liminf_{n \to \infty}\frac{1}{n}\log 
\frac{1}{\#\Sc_{\rho, L}(n)}\#\left\{x \in \Sc_{\rho, L}(n) \ : \ \frac{1}{n}\log \|\rho_\ast(x)\| \in U\right\}\\
&\le 
\limsup_{n \to \infty}\frac{1}{n}\log 
\frac{1}{\#\Sc_{\rho, L}(n)}\#\left\{x \in \Sc_{\rho, L}(n) \ : \ \frac{1}{n}\log \|\rho_\ast(x)\| \in V\right\}
\le -\inf_{t \in V}I(t),
\end{align*}
where 
\[
I(t):=\sup_{s \in \R}\{ts-\th(-s)\}+\th(0) \quad \text{for $t \in \R$},
\]
and $\th(s)$ is the pressure curve for the pair of dominated representations $(\rho, \rho_\ast)$.

Furthermore, for a large enough constant $L>0$,
if we define
\[
\Cc_{\rho, L}(n):=\{[x] \in \conj_\G \ : \ |\log \lambda(\rho(x))-n|\le L\},
\]
then the large deviation principle holds for the sequence $(1/n)\log \lambda(\rho_\ast(x))$ where $[x]$ is distributed according to the uniform counting measure on $\Cc_{\rho, L}(n)$ for $n=1, 2, \dots$ with the same rate function $I$.
\end{proposition}

\proof
By Lemma \ref{Lem:growth}, for each $s \in \R$, we have that
\[
\lim_{n \to \infty}\frac{1}{n}\log \frac{1}{\# \Sc_{\rho, L}(n)}\sum_{x: \ |\log \|\rho(x)\|-n|\le L}\exp\(s \log \|\rho_\ast(x)\|\)=\th(-s)-\th(0).
\]
Since $\th(-s)$ is continuously differentiable (in fact, real analytic) by Theorem \ref{Thm:Manhattan_rep},
the G\"artner-Ellis theorem (e.g., \cite[Theorem 2.3.6]{DemboZeitouni}) implies the first claim.

Concerning the second claim, by Proposition \ref{Prop:Knieper}, 
we have that for each $s \in \R$,
\[
C_{s, L, R}^{-1}\,n^{-2}e^{\th(s)n} \le \sum_{[x]: \ |\log \lambda(\rho(x))-n| \le L}\exp\(-s \log \lambda(\rho_\ast(x))\) \le C_{s, L, R}\,e^{\th(s)n},
\]
for all large enough $n$.
Hence for all large enough $L>0$ and for each $s \in \R$,
\[
\lim_{n \to \infty}\frac{1}{n}\log \frac{1}{\#\Cc_{\rho, L}(n)}\sum_{x: \ |\log \lambda(\rho(x))-n|\le L}\exp\(s \log \lambda(\rho_\ast(x))\)=\th(-s)-\th(0).
\]
The rest follows as in the proof of the first claim.
\qed

\begin{theorem}\label{Thm:intersection_norm}
Let $(\rho, \rho_\ast)$ be a pair of dominated representations from a non-elementary hyperbolic group $\G$ to $\GL(m, \R)$ for $m \ge 2$.
There exists a finite positive constant $\t(\rho_\ast/\rho)$ as the following limit:
\[
\t(\rho_\ast/\rho)=\lim_{T \to \infty}\frac{1}{\#\Bc_{\rho}(T)}\sum_{x \in \Bc_{\rho}(T)\setminus \Bc_\rho(0)}\frac{\log \|\rho_\ast(x)\|}{\log \|\rho(x)\|},
\]
where $\Bc_{\rho}(T):=\{x \in \G \ : \ \log \|\rho(x)\| \le T\}$ for positive real numbers $T$.

Moreover, we have that
\[
v_\rho \le \t(\rho_\ast/\rho)v_{\rho_\ast},
\]
where $v_\rho$ and $v_{\rho_\ast}$ are the exponential growth rates for $\rho$ and $\rho_\ast$, respectively.
\end{theorem}

We call the limit $\t(\rho_\ast/\rho)$ in Theorem \ref{Thm:intersection_norm} the {\it intersection number} for a pair of dominated representations $(\rho, \rho_\ast)$.

\proof[Proof of Theorem \ref{Thm:intersection_norm}]
Applying the cocycle $v_\rho c_\rho(x, \x)$ to Proposition \ref{Prop:PSpotential},
we obtain a finite Borel measure $\m_\rho$ satisfying \eqref{Eq:QC} and the exponent is $0$ relative to any metric $d$ in $\Dc_\G$ and for positive constants $C, L>0$,
\begin{equation}\label{Eq:mrho}
C^{-1} \le \exp(v_\rho \log \|\rho(x)\|)\cdot \m_\rho(\Oc(x, R)) \le C \quad \text{for all $x \in \G$},
\end{equation}
and 
\begin{equation}\label{Eq:count}
C^{-1}\exp(v_\rho n) \le \#\{x \in \G \ : \ |\log \|\rho(x)\|-n| \le L\} \le C \exp(v_\rho n),
\end{equation}
for all positive integers $n$ by letting $s=0$ in Lemma \ref{Lem:growth}.
Fix a metric $d$ in $\Dc_\G$.
By Lemma \ref{Lem:localint} applied to the $\G$-invariant tempered potentials 
\[
\psi_\ast(x, y):=\log \|\rho_\ast(x^{-1}y)\| \quad \text{and} \quad \psi(x, y):=\log \|\rho(x^{-1}y)\|,
\]
there exist constants $\t_\ast$ and $\t$ such that
\[
\t^{\psi_\ast}(\x)=\t_\ast \quad \text{and} \quad \t^{\psi}(\x)=\t \quad\text{for $\m_\rho$-almost every $\x \in \partial \G$}.
\]
Namely, for $\m_\rho$-almost every $\x \in \partial \G$,
we have that $\psi_\ast(o, \g_\x(n))/n \to \t_\ast$ as $n \to \infty$
where $\g_\x$ is any rough geodesic ray relative to $d$ issuing from $o$,
and similarly for $\psi(o, \g_\x(n))$.
This implies that
\[
\frac{\psi_\ast(o, \g_\x(n))}{\psi(o, \g_\x(n))} \to \frac{\t_\ast}{\t} \quad \text{for $\m_\rho$-almost every $\x \in \partial \G$ as $n \to \infty$}.
\]
Let
\[
S_{n, \e}:=\left\{x \in \G \ : \ |\psi_\ast(o, x)-n (\t_\ast/\t)|>\e n \right\} \quad \text{for $n \ge 0$ and $\e>0$},
\]
and $S(n, L):=\{x \in \G \ : \ |\psi(o, x)-n|\le L\}$.
For a large enough $R>0$,
the shadows $\Oc(x, R)$ for $x \in S(n, L)$ cover the boundary $\partial \G$ with a uniformly bounded multiplicity.
Therefore we have by \eqref{Eq:mrho} and \eqref{Eq:count},
\begin{equation*}
\frac{\#\(S_{n, \e}\cap S(n, L)\)}{\# S(n, L)} \le C'\sum_{x \in S_{n, \e}\cap S(n,L)}\m_\rho(\Oc(x, R))\le C''\m_\rho\left(\bigcup_{x \in S_{n,\e}\cap S(n, L)}\Oc(x, R)\right) \to 0,
\end{equation*}
as $n \to \infty$,
where the last term tends to $0$ as $n$ goes to infinity by the Fatou lemma; indeed
for $\x \in \Oc(x, R)$ and $x \in S_{n, \e}\cap S(n, L)$, 
if $d(\g_\x(k), x) \le R$,
then
\[
|\psi_\ast(o, \g_\x(k))- (\t_\ast/\t)\psi(o, \g_\x(k))| \ge \e \psi(o, \g_\x(k))-R_0,
\]
where $R_0$ is a constant independent of $n$; it depends only on the constants in the condition \eqref{Eq:RG} for $\psi_\ast$ and $\psi$, the metric $d$, $\e$, $L$ and $R$.
This implies that for any given $\e>0$, for a uniformly distributed random element in $S(n, L)$,
we have $|\psi_\ast(o, x)-(\t_\ast/\t)\psi(o, x)| \le \e \psi(o, x)$ with probability at least $1-\e$ for all large enough $n$, and thus
\[
\frac{\t_\ast}{\t}=\lim_{n \to \infty}\frac{1}{\#S(n, L)}\sum_{x \in S(n, L)}\frac{\psi_\ast(o, x)}{n}.
\] 
For any real $T>0$, let $n:=\lfloor T\rfloor$ be the largest integer at most $T$.
If $x_n$ is distributed according to the uniform counting measure on the set 
$\Bc_{\rho}(T):=\{x \in \G \ : \ \log \|\rho(x)\| \le T\}$,
then since $x_n \notin S(n, L)$ with probability at most $O(e^{-v_\rho L})$ by \eqref{Eq:count} and the probability that $x_n \in S_{n, \e}\cap S(n, L)$ tends to $0$ as $n$ goes to infinity,
the event that $x_n \in S_{n, \e}$ has the probability at most $O(e^{-v_\rho L})$ for all large enough $n$.
Hence letting $T \to \infty$ and then $L \to \infty$,
we have that
\[
\frac{\t_\ast}{\t}=\lim_{T \to \infty}\frac{1}{\#\Bc_{\rho}(T)}\sum_{x \in \Bc_{\rho}(T)\setminus \Bc_\rho(0)}\frac{\psi_\ast(o, x)}{\psi(o, x)},
\]
since $\Bc_\rho(0)$ is finite by \eqref{Eq:dom_QI}.
Defining $\t(\rho_\ast/\rho):=\t_\ast/\t$, we obtain the first claim.

Moreover, the above discussion shows that for all small enough $\e>0$ and all large enough $L$ and $T$,
\[
 (1-C e^{-v_\rho L})\# \Bc_\rho(T) \le \#\Bc_{\rho_\ast}((\t(\rho_\ast/\rho)+\e)T),
\]
where $C$ is independent of $L$ or $T$.
Since 
\[
c_1 \exp(v_\rho T) \le \#\Bc_\rho(T)\le c_2 \exp(v_\rho T),
\]
and similarly for $\#\Bc_{\rho_\ast}(T)$,
by taking logarithms, dividing by $T$ and letting $T \to \infty$,
we obtain the second claim.
\qed

\begin{corollary}\label{Cor:zero}
The rate function $I$ for the large deviation principles in Proposition \ref{Prop:LDP} has a unique zero at $\t(\rho_\ast/\rho)$.
\end{corollary}

\proof
We have $\th'(0)=-\t(\rho_\ast/\rho)$ by Theorem \ref{Thm:Manhattan_rep} and by the second claim of Theorem \ref{Thm:Manhattan} adapted to the pair of representations $(\rho, \rho_\ast)$ 
(cf.\ the proof of Theorem \ref{Thm:intersection_norm}).
Thus $I$ has a zero at $\t(\rho_\ast/\rho)$.
If $t \neq \t(\rho_\ast/\rho)$, then $I(t)>0$ 
since $\th$ is convex and differentiable at $0$.
\qed

\medskip

The following is shown in \cite[Theore 3.7 and (12) in p.1113]{BCLS}.
We provide another proof by using the large deviation principles in Proposition \ref{Prop:LDP}.

\begin{theorem}\label{Thm:intersection_rep}
Let $(\rho, \rho_\ast)$ be a pair of dominated representations from a non-elementary hyperbolic group $\G$ to $\GL(m, \R)$ for $m \ge 2$.
There exists a finite positive constant $\t(\rho_\ast/\rho)$ as the following limit:
\[
\t(\rho_\ast/\rho)=\lim_{T \to \infty}\frac{1}{\#\Cc_{\rho}(T)}\sum_{[x] \in \Cc_{\rho}(T)\setminus \Cc_\rho(0)}\frac{\log \lambda(\rho_\ast(x))}{\log \lambda(\rho(x))},
\]
where $\Cc_{\rho}(T):=\{[x] \in \conj_\G \ : \ \log \lambda(\rho(x)) \le T\}$ for positive real numbers $T$.
\end{theorem}

\proof
We recall that for a fixed large enough $L>0$ and for positive integers $n$,
\[
\Cc_{\rho, L}(n)=\{[x] \in \conj_\G \ : \ |\log \lambda(\rho(x))-n|\le L\},
\]
in the second claim in Proposition \ref{Prop:LDP}; it shows that 
since the rate function $I$ has a unique zero at $\t(\rho_\ast/\rho)$ by Corollary \ref{Cor:zero},
for each $\e>0$ there exists $0<c_\e<v_\rho$ such that for all large enough $n$,
\[
\#\{[x] \in \Cc_{\rho, L}(n) \ : \ |\log \lambda(\rho_\ast(x))-n\t(\rho_\ast/\rho)| \ge \e n\} \le \#\Cc_{\rho, L}(n)\cdot e^{-c_\e n} \le C e^{(v_\rho-c_\e)n},
\]
where in the last inequality we have used $\#\Cc_{\rho, L}(n) \le C e^{v_\rho n}$ by Proposition \ref{Prop:Knieper} at $s=0$.
Hence
\begin{align*}
\sum_{n=0}^{\lceil T\rceil}\#\{[x] \in \Cc_{\rho, L}(n) \ : \ |\log \lambda(\rho_\ast(x))-n\t(\rho_\ast/\rho)| \ge \e n\} 
\le \sum_{n=0}^{\lceil T\rceil}C e^{(v_\rho-c_\e)n} \le C' e^{(v_\rho-c_\e)T},
\end{align*}
for all $T>0$,
where $\lceil T\rceil$ is the smallest integer greater than or equal to $T$.
Moreover, 
we have a positive constant $C>0$ such that for all large enough $T$ and $n:=\lfloor T\rfloor$,
\[
\#\Cc_\rho(T) \ge \#\Cc_{\rho, L}(n-\lceil L \rceil) \ge \frac{1}{Cn^2}e^{v_\rho n},
\]
which follows by letting $s=0$ in Proposition \ref{Prop:Knieper}.
Therefore we obtain for each $\e>0$,
\begin{align*}
&\frac{1}{\#\Cc_\rho(T)}\#\{[x] \in \Cc_\rho(T) \ : \ |\log \lambda(\rho_\ast(x))-\t(\rho_\ast/\rho)\log \lambda(\rho(x))| \ge \e \log \lambda(\rho(x))\}\\
&\qquad \qquad \qquad \qquad \qquad \qquad \qquad \le \frac{C n^2}{e^{v_\rho n}}C' e^{(v_\rho-c_\e)n}
\le C'' n^2 e^{-c_\e n} \to 0 \quad \text{as $n \to \infty$},
\end{align*}
implying that 
\[
\frac{1}{\#\Cc_\rho(T)}\sum_{[x] \in \Cc_\rho(T)\setminus\Cc_\rho(0)}\frac{\log \lambda(\rho_\ast(x))}{\log \lambda(\rho(x))} \to \t(\rho_\ast/\rho) \quad \text{as $T \to \infty$},
\]
since $\Cc_\rho(0)$ is finite; \eqref{Eq:dom_QI_sp} and there are only finitely many conjugacy class of torsion elements (e.g., \cite[Lemma 2.3.6]{Calegari}),
as claimed.
\qed

\subsection{Real analytic variations of dominated representations}

Let $\Uc$ be a real analytic manifold.
We say that a family of representation $\rho_u: \G \to \GL(m, \R)$ parameterized by $u \in \Uc$ is a \textit{real analytic family} if for each $x \in \G$, the map $u \mapsto \rho_u(x)$ yields a real analytic map with values in $\GL(m, \R)$.
Recall that associated with any dominated representation $\rho$,
there exists a pair of H\"older continuous $\rho$-equivariant maps,
\[
\check E_\rho: \partial \G\to \Pb({\R^m}^\ast) \quad \text{and} \quad E_\rho:\partial \G \to \Pb(\R^m).
\]
Let us consider the spaces $C^\a(\partial \G, \Pb(\R^m))$ and $C^\a(\partial \G, \Pb({\R^m}^\ast))$ of $\a$-H\"older maps for some $\a>0$ as Banach manifolds.

\begin{theorem}[Theorem 6.1 in \cite{BCLS} and Theorem 6.1 in \cite{BochiPotrieSambarino}]\label{Thm:RA}
Let $\rho_u:\G \to \GL(m, \R)$ for $u \in \Uc$ be a real analytic family of representations.
If $\rho_{u_0}$ is dominated at $u_0 \in \Uc$,
then there exists a neighborhood $\Uc'$ of $u_0$ in $\Uc$ such that
for every $u \in \Uc'$,
the representations $\rho_u$ are dominated.
Moreover, for some $\a>0$, the map from $\Uc'$ to $C^\a(\partial \G, \Pb(\R^m))$, $u \mapsto E_{\rho_u}$ and the map from $\Uc'$ to $C^\a(\partial \G, \Pb({\R^m}^\ast))$, $u \mapsto \check E_{\rho_u}$ are real analytic, respectively.
\end{theorem}

The following is shown in \cite[Propositions 3.12 and 8.1]{BCLS}.
We give a proof by using Theorem \ref{Thm:enhanced}.

\begin{theorem}
Let $\rho_{u}$ for $u \in \Uc$ and $\rho_{\ast v}$ for $v \in \Vc$ be two real analytic families of dominated representations from $\G$ to $\GL(m, \R)$ for $m \ge 2$.
The exponential growth rate $v_{\rho_u}$ is a real analytic function on $\Uc$, and
the intersection number $\t(\rho_{\ast v}/\rho_u)$ is a real analytic function on $\Uc \times \Vc$.
\end{theorem}

\proof
We consider $(\SS_0, \s)$ in Theorem \ref{Thm:enhanced}.
Letting $c_{\rho_u}(x, \x)$ be the cocycle associated with $\rho_u$ for $u \in \Uc$, 
we define
\[
\Psi_u(\o):=c_{\rho_u}(w(\o_0), w_\ast(\o)) \quad \text{for $\o \in \SS_0$}.
\]
Note that $\Psi_u$ is a H\"older continuous potential for each $u \in \Uc$.
Since $v_{\rho_u}c_{\rho_u}(x, \x)$ has the exponent $0$ with respect to any metric in $\Dc_\G$,
by Lemma \ref{Lem:Pr}, we have $\Pr(v_{\rho_u}\Psi_u)=0$, and thus
$\Pr_{\Cc_0}(v_{\rho_u}\Psi_u)\le 0$.
For each $u \in \Uc$, there exists a $\G$-invariant Radon measure $\L_u$ on $\partial^2 \G$ equivalent to $\m_{\check \rho_u}\otimes \m_{\rho_u}$.
By Proposition \ref{Prop:code-measure}(2),
there exists a $\s$-invariant Borel probability measure $\lambda_u$ on $\SS_0$ such that
\[
w_\ast \lambda_u \otimes \Leb_{[0, T_0)}\le c \L_u \otimes dt \qquad \text{on $\partial^2 \G \times \R$},
\]
for some positive constants $T_0, c>0$, further $w_\ast \lambda_u \le C \m_{\check \rho_u}\otimes \m_{\rho_u}$ on $\partial^2 \G$ for some positive constant $C>0$.
For the fixed strongly hyperbolic metric $\hat d$ used to define $\Sus(\SS_0, r_0)$ in Theorem \ref{Thm:enhanced},
let $\hat q$ be the corresponding quasi-metric in $\partial \G$.
We then have that
\[
D_{\hat q}(\x_-, \x_+; \L_u)=\hat \t_u \qquad \text{for $\L_u$-almost all $(\x_-, \x_+) \in \partial^2 \G$},
\]
where $\hat \t_u$ is the local intersection number of $\log \|\rho_u(x^{-1}y)\|$ relative to $\hat d$ for $\m_{\rho_u}$-almost every point in the boundary by Lemma \ref{Lem:localint}.
Furthermore, Lemma \ref{Lem:localint_potential} shows that
\[
\int_{\SS_0}\Psi_u\,d\lambda_u=\hat \t_u\int_{\SS_0}\Psi_B\,d\lambda_u.
\]
Lemma \ref{Lem:dim-ent} shows that
\[
h(\s, \lambda_u)+\int_{\SS_0}\Psi_u\,d\lambda_u=h(\s, \lambda_u)+\hat \t_u\int_{\SS_0}\Psi_B\,d\lambda_u \ge 0,
\]
and by the variational principle (Proposition \ref{Prop:VP}),
we have $\Pr_{\Cc_0}(v_{\rho_u}\Psi_u)\ge 0$.
Hence $\Pr_{\Cc_0}(v_{\rho_u}\Psi_u)=0$.
Proposition \ref{Prop:VP} shows that $\lambda_u$ is the unique equilibrium state for $\Psi_u$.
Note that $\Uc\times \R \to \R$, $(u, s)\mapsto \Pr_{\Cc_0}(s\Psi_u)$ is real analytic by Theorems \ref{Thm:TF}(2) and \ref{Thm:RA}.
Furthermore, we have that
\[
\frac{\partial}{\partial s}\Big|_{s=v_{\rho_u}}\Pr_{\Cc_0}(s\Psi_u)=\int_{\SS_0}\Psi_u\,d\lambda_u <0,
\]
where the last strict inequality holds since $S_n \Psi_u<0$ for all large enough $n$ by Lemmas \ref{Lem:rep_RG} and \ref{Lem:rep_QE} together with \eqref{Eq:dom_QI}.
The implicit function theorem shows that on an open neighborhood $\Uc_0$ of $u$ in $\Uc$, the function $w \mapsto v_{\rho_w}$ for $w \in \Uc_0$ is real analytic.
Since this holds for each $u \in \Uc$, the function $u \mapsto v_{\rho_u}$ for $u \in \Uc$ is real analytic,
concluding the first claim.

Let $\Psi_{\ast v}(\o):=c_{\rho_{\ast v}}(w(\o_0), w_\ast(\o))$ for $\o \in \SS_0$ a family of H\"older continuous potentials associated with $\rho_{\ast v}$ for $v \in \Vc$.
Since 
the function $\Uc \times \R \to \R$ and the function $\Uc\times \Vc \times \R \to \R$ defined by
\[
(u, s) \mapsto \Pr_{\Cc_0}(s\Psi_u), \quad \text{and} \quad (u, v, s) \mapsto \Pr_{\Cc_0}(v_{\rho_u}\Psi_u+s \Psi_{\ast v}),
\]
respectively, are real analytic by Theorems \ref{Thm:TF}(2) and \ref{Thm:RA} combined with the real analyticity of $v_{\rho_u}$ in $u \in \Uc$ in the first claim.
Hence for $u \in \Uc$ and for $(u, v) \in \Uc \times \Vc$, the following values depend real analytically in $u$ and in $(u, v)$, respectively,
\[
\frac{\partial}{\partial s}\Big|_{s=v_{\rho_u}}\Pr_{\Cc_0}(s \Psi_u)=\int_{\SS_0}\Psi_u\,d\lambda_u 
\quad \text{and} \quad
\frac{\partial}{\partial s}\Big|_{s=0}\Pr_{\Cc_0}(v_{\rho_u}\Psi_u+s\Psi_{\ast v})=\int_{\SS_0}\Psi_{\ast v}\,d\lambda_u.
\]
Noting that
\[
\t(\rho_{\ast v}/\rho_u)=\frac{\int_{\SS_0}\Psi_{\ast v}\,d\lambda_u}{\int_{\SS_0}\Psi_u\,d\lambda_u},
\]
by Lemma \ref{Lem:localint_potential} (cf.\ the proof of Theorem \ref{Thm:intersection_norm}).
Therefore, we conclude that $\t(\rho_{\rho_{\ast v}}/\rho_u)$ is a real analytic function on $\Uc \times \Vc$ and thus the second claim.
\qed

\section{Examples}\label{Sec:example}

\subsection{A free product group}\label{Sec:free}

Consider the following group with presentation
\[
\G:=\langle s_1, s_2 \ | \ s_1^4, s_2^6\rangle,
\]
which is isomorphic to the free product $(\Z/4)\ast(\Z/6)$.
Let $S:=\{s_1, s_1^{-1}, s_2, s_2^{-1}\}$.
We define a family of presentations $\rho_u$ of $\G$ into $\SL(2, \R)$ by setting
\[
\rho_u(s_1):=A_{-u}R(\pi/2)A_u \quad \text{and} \quad \rho_u(s_2):=R(\pi/3),
\]
where
\[
A_u:=
\begin{pmatrix}
e^u & 0\\
0 & e^{-u}
\end{pmatrix}
\quad \text{and} \quad
R(t):=
\begin{pmatrix}
\cos {t} & -\sin{t}\\
\sin {t} & \cos {t}
\end{pmatrix}
\quad \text{where $u \in \R$ and $t\in [0, \pi)$}.
\]
It indeed defines a representation since $\rho_u(s_1)^2=\rho_u(s_2)^3=-I$,
and further it descends to $\PSL(2, \R)$.
Note that $\{\rho_u\}_{u \in \R}$ is a real analytic family.
If $|u|$ is large enough, then $\rho_u$ defines a dominated representation (cf.\ \cite[Section 5.5]{BochiPotrieSambarino}).
For $u \in \R$ such that $\rho_u$ is dominated,
let $\m_{\rho_u}$ be an associated finite Borel measure on $\partial \G$ (see the end of Section \ref{Sec:dom_rep}), and $\t(\rho_u/d_S)$ be the local intersection number of the associated tempered potential with $\rho_u$ over the word metric $d_S$.
By construction, if $u$ is large, then $\t(\rho_u/d_S)$ is large.
Since the local intersection number is constant $\m_{\rho_u}$-almost everywhere on $\partial \G$ (see Lemma \ref{Lem:localint}),
if $\t(\rho_u/d_S) \neq \t(\rho_{u'}/d_S)$,
then $\m_{\rho_u}$ and $\m_{\rho_{u'}}$ are mutually singular for $u, u' \in \R$.

\subsection{A surface group with the translation generators}\label{Sec:surface}
\begin{figure}
\centering
\includegraphics[width=120mm, bb=0 0 595 842]{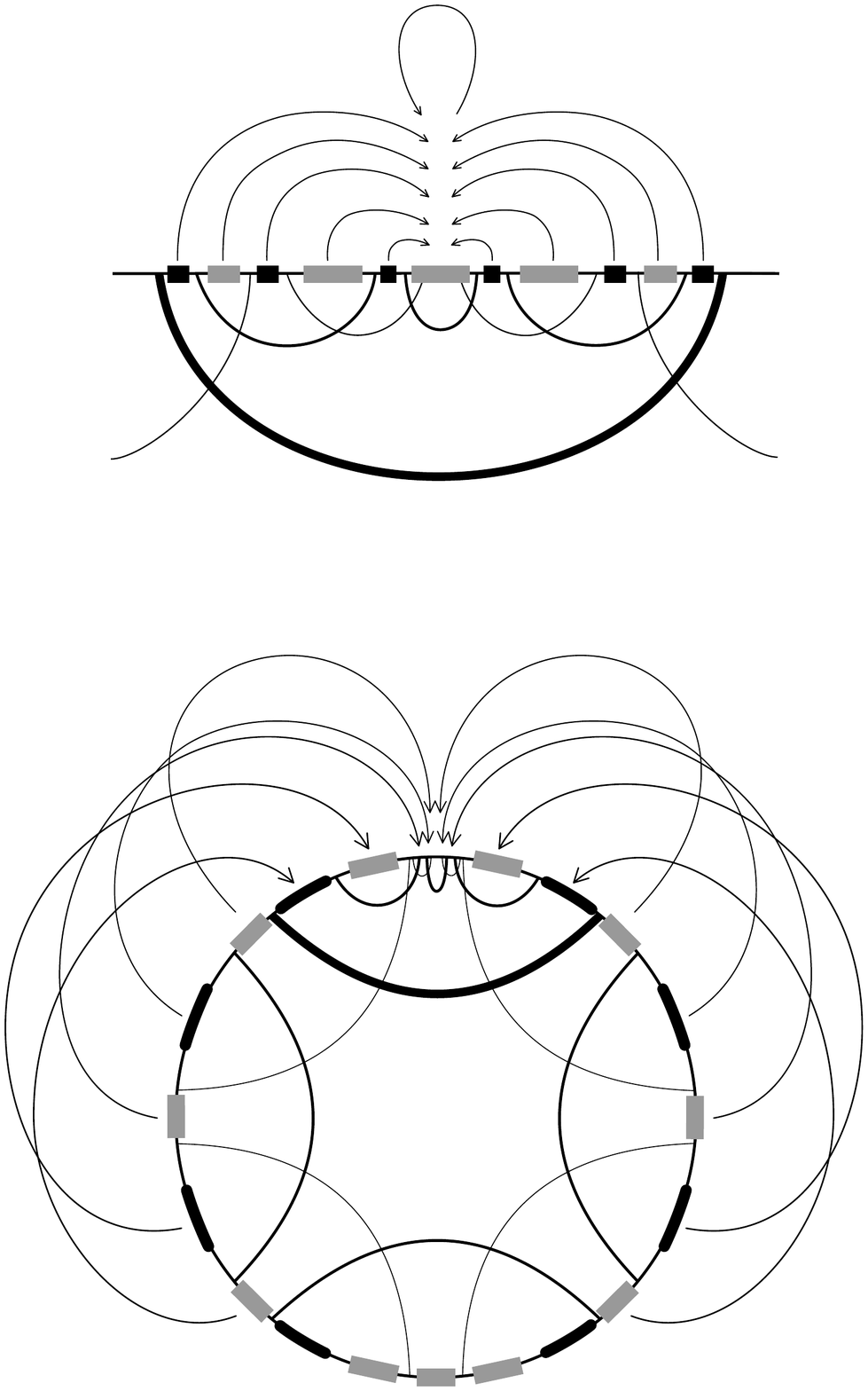}
\caption{A part of a strictly invariant family of multicones of index $1$ for $(\G, S)$, where $\G$ is a genus two surface group and $S$ is the translation generators.
The arrows indicate the destinations of intervals after translation by a single generator $s \in S$ (bottom), and an expanded picture shows a small portion of the disk bounded by the thick arc (top).}
\label{Fig:multicones}
\end{figure}

Let us consider the fundamental group of closed Riemann surface of genus $2$ with the set of generators $S:=\{s_i, s_i^{-1}, i=1,2,3,4\}$ arising from identifying pairs of opposite sides of a regular octagon in the hyperbolic plane, i.e., 
\[
\G:=\langle s_1, s_2, s_3, s_4 \ | \ s_4^{-1}s_3^{-1}s_2^{-1}s_1^{-1}s_4 s_3 s_2 s_1\rangle.
\]
We will refer to $S$ as the set of {\it translation generators}.
We define a representation of $\G$ into $\SL(2, \R)$ which acts via $\PSL(2, \R)$ and admits a hyperbolic regular octagon as a fundamental polygon.
Letting
\[
t:=\frac{\pi}{8} \quad \text{and} \quad u:=\log\frac{\sqrt{\cos 2t}+\cos t}{\sin t},
\]
we define a homomorphism
\[
\rho: \G \to \SL(2, \R)
\]
by
\[
\rho(s_1):=A_u,
\quad
\rho(s_2):=R(3t)A_u R(-3t),
\]
\[
\rho(s_3):=R(6t)A_u R(-6t),
\quad \text{and} \quad
\rho(s_4):=R(t)A_u R(-t),
\]
where we have used the same notations in the previous Section \ref{Sec:free}.
One can check that $\rho$ is indeed a homomorphism since the group generated by $\rho(s_i)$, $i \in \{1, 2, 3, 4\}$ admits a regular octagon as a Dirichlet fundamental domain and they satisfy the required relation.
Moreover, the representation $\rho$ is dominated;
this follows since the associated topological flow is H\"older orbit equivalent to the classical geodesic flow on $\rho(\G)\backslash \SL(2, \R)/{\{\pm I\}}$ (cf.\ \cite[Theorem A]{BochiGourmelon}).
More concretely, one can describe a strictly invariant family of multicones of index $1$ (Figure \ref{Fig:multicones}; see \cite[Section 5]{BochiPotrieSambarino} for the details).

\begin{figure}
\begin{minipage}[t]{0.5\columnwidth}
\centering
\includegraphics[clip, width=65mm, bb=0 0 360 214]{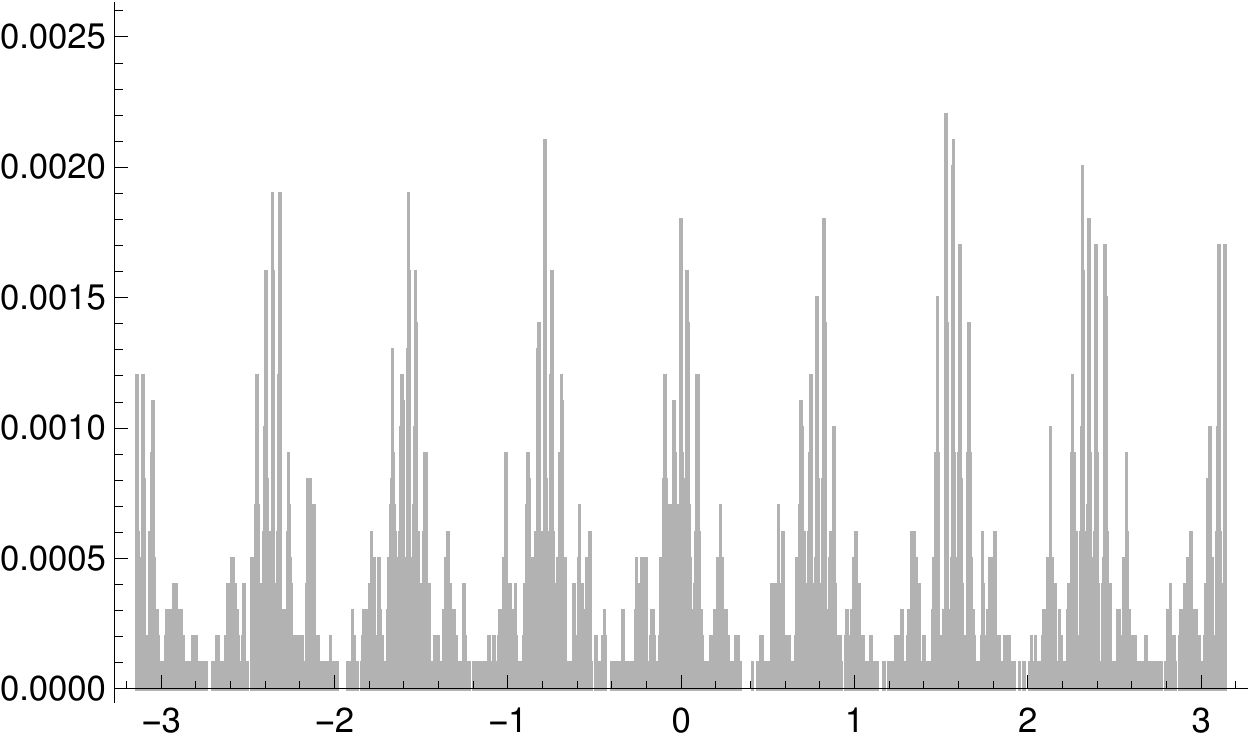}
\subcaption{A Patterson-Sullivan measure}
\end{minipage}%
\begin{minipage}[t]{0.5\columnwidth}
\centering
\includegraphics[clip, width=65mm, bb=0 0 360 214]{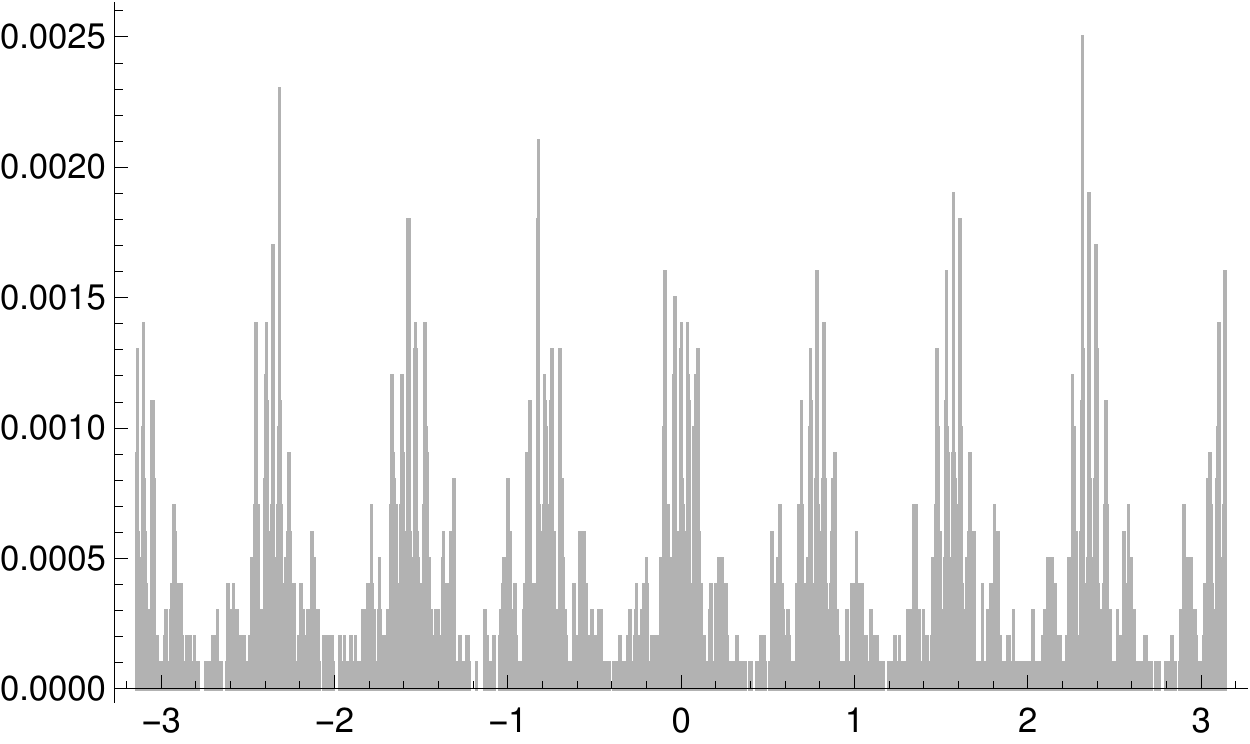}
\subcaption{The harmonic measure}
\end{minipage}
\caption{The histograms (the number of steps $n=10^3$ and the number of samples $10^4$) for
the distributions of angles in $[-3.14\dots, 3.14\dots]$ of $\rho(w_n).o$ where $o$ is the center of the unit disk:
 $\{w_n\}_{n=0, 1, \dots}$ is associated with the Markov chain on the strongly connected component in the underlying directed graph in an automatic structure for $(\G, S)$  (left), and $\{w_n\}_{n=0, 1, \dots}$ is the simple random walk on the Cayley graph for $(\G, S)$ (right).}
\label{Fig:histogram_translation}
\end{figure}

The group $\G$ with the set of generators $S$ admits a deterministic automaton on $8$ letters with $37$ states \cite{GAP4}.
The underlying directed graph has a unique strongly connected component with $36$ states.
We have performed a Markov chain $\{w_n\}_{n=0, 1, \dots}$ on the strongly connected component where the transition probability is given by the Parry measure; the extreme point of sample path is distributed according to a Patterson-Sullivan measure $\m_S$ defined in terms of the word metric for $S$ (cf.\ \cite[Section 4.2]{CalegariFujiwara2010} and Figure \ref{Fig:histogram_translation} (left)).
A numerical experiment shows that
\[
\t(\rho/d_S)=\lim_{n \to \infty}\frac{1}{n}\Eb\log \s_1(\rho(w_n)) \quad \text{is approximately}\quad
1.13837...,
\]
(where $\Eb$ stands for the expectation-----we compute an average over the number of samples $10^4$ with the number of steps $n=10^3$).
The exponential growth rate of $(\G, S)$ is
\[
v_{S}=1.94303...,
\]
and $v_\rho=2$ (since the exponential growth rate for the hyperbolic metric is $1$ and $2\log \s_1(\rho(w))$ coincides with the displacement by $\rho(w)$ in the hyperbolic plane).
Computations suggest that
\[
\frac{v_S}{v_\rho}<1.13837....
\]
In fact, we can actually show that a finite Borel measure $\m_\rho$ associated with $\rho$ (where we note that $\m_\rho$ is in the Lebesgue measure class if it is realized on the boundary of the hyperbolic plane) and a Patterson-Sullivan measure $\m_S$ associated with the word metric for $S$ are mutually singular.
The following is a sketch of the proof:
Suppose that two measures $\m_\rho$ and $\m_S$ are mutually absolutely continuous,
then the logarithm of spectral radius for $\rho$ and the stable length function are proportional on the set of conjugacy classes of $\G$; this follows from \cite[Theorem 1.2]{CT} for the pair of the hyperbolic metric and the word metric.
However, the values for the one associated with $\rho$ is not arithmetic, i.e., not included in a discrete subgroup in $\R$ \cite[Proposition 2.1]{Dalbo} whereas the values for the one associated with the word metric are contained in a discrete set of rational numbers $\Z[1/N]$ for some $N \ge 1$ \cite[Theorem 3.17, Chapter III.$\Gamma$]{BridsonHaefliger}, yielding a contradiction.

Concerning the harmonic measure, we have performed the simple random walks $\{w_n\}_{n=0, 1, \dots}$ on the Cayley graph for $(\G, S)$ (Figure \ref{Fig:histogram_translation} (right)).
Let $\n$ be the associated harmonic measure on the boundary $\partial \G$, and $\L_{\rm harm}$ be a $\G$-invariant Radon measure in the class of $\check \n\otimes \n$ on $\partial^2 \G$ (see Section \ref{Sec:RW}).
A numerical experiment provides that
\[
\text{$\t(\rho/d_S: \L_{\rm harm})$ is approximately $1.12909...$.}
\]
Note that the local intersection number for $\t(\rho/d_S:\L)$ relative to an ergodic measure $\L$ is constant $\L$-almost everywhere,
and $\t(\rho/d_S)=\t(\rho/d_S: \L_S)$ where $\L_S$ is a $\G$-invariant Radon measure in the class of $\m_S\otimes \m_S$ on $\partial^2 \G$.
Therefore if $\t(\rho/d_S)\neq \t(\rho/d_S: \L_{\rm harm})$,
then a Patterson-Sullivan measure $\m_S$ and the harmonic measure $\n$ are mutually singular on $\partial \G$;
the fact itself has been shown as a special case of \cite[Theorems 1.2 and 1.3]{GMM2018}.
Moreover, we note that $\m_\rho$ and $\n$ are mutually singular \cite[Theorem 1]{KosenkoTiozzo};
see also the discussion in Remark \ref{Rem:multifractal}.

\subsection*{Acknowledgment}
We thank Eduardo Oreg\'on-Reyes for his helpful comments and the referee for a careful review that helped us to improve the presentation of the paper.
R.T.\ is partially supported by JSPS Grant-in-Aid for Scientific Research JP20K03602. 
This work was supported by the Research Institute for Mathematical
Sciences, an International Joint Usage/Research Center located
in Kyoto University.

\bibliographystyle{alpha}
\bibliography{coding}

\end{document}